\documentclass[a4paper,11pt]{amsart}

\usepackage{amsthm,amsmath,amssymb,stmaryrd,enumerate,euscript,mathrsfs,a4,array}
\usepackage{geometry}

\usepackage{pictexwd,dcpic,graphicx,cite}
\usepackage[export]{adjustbox}
\usepackage{subcaption}
\usepackage{url}
\usepackage[latin1]{inputenc} 
\usepackage{indentfirst}
\usepackage{calligra} 
\usepackage{cleveref}
\usepackage{graphicx}
\usepackage{amssymb}
\usepackage{fancyhdr}
\usepackage{amsmath}
\usepackage{latexsym}
\usepackage{amsthm}
\usepackage[all]{xy} 
\usepackage{amscd}
\usepackage{tikz-cd}
\usepackage{tikz}
\usetikzlibrary{arrows.meta}
\numberwithin{equation}{section} 

\graphicspath{ {images/} }

\newtheorem{theorem}{Theorem}[section]
\newtheorem*{theorem*}{Theorem}
\newtheorem{prop}[theorem]{Proposition}

\newtheorem{lemma}[theorem]{Lemma}
\newtheorem{cor}[theorem]{Corollary}

\theoremstyle{definition} 
\newtheorem{defn}[theorem]{Definition}

\theoremstyle{remark}
\newtheorem{rmk}[theorem]{Remark}

\newtheorem*{remark*}{Remark}
\newtheorem*{remarks*}{Remarks}
\newtheorem*{notation*}{Notation}
\newtheorem*{convention*}{Convention}

\newcommand{\op}{\mathit{op}}

\newcommand{\co}{\,}

\newcommand{\cat}{\mathcal}

\newcommand{\catK}{\mathcal{K}}

\newcommand{\tcat}{\mathbf}

\newcommand{\bcat}{}
\newcommand{\bx}{\bcat{X}}
\newcommand{\by}{\bcat{Y}}

\newcommand{\Gray}{\tcat{Gray}}

\newcommand{\Mnd}{\mathbf{Mnd}}
\newcommand{\Psm}{\mathbf{Psm}}
\newcommand{\Lift}{\mathbf{Lift}}

\newcommand{\algfont}{\mathrm}

\newcommand{\psSalg}{\algfont{Ps}\text{-}\ms\text{-}\algfont{Alg}}

\newcommand{\psTalg}{\algfont{Ps}\text{-}\mt\text{-}\algfont{Alg}}

\newcommand{\ms}{S}
\newcommand{\mss}{\ms^2}
\newcommand{\msss}{\ms^3}
\newcommand{\mssss}{\ms^4}

\newcommand{\es}{s}
\newcommand{\mus}{ {m} }

\newcommand{\mt}{T}
\newcommand{\mtt}{\mt^2}
\newcommand{\mttt}{\mt^3}

\newcommand{\et}{t}
\newcommand{\mut}{ {n} }

\newcommand{\mts}{\mt \ms}
\newcommand{\mst}{\ms \mt}  
\newcommand{\msst}{\mss \mt}  
\newcommand{\mtss}{\mt \mss}  
\newcommand{\mtts}{\mtt \ms}  
\newcommand{\msts}{\ms \mt \ms}  
\newcommand{\mtst}{\mt \ms \mt}  
\newcommand{\mstt}{\ms \mtt}

\newcommand{\mh}{F}
\newcommand{\mk}{F'}

\newcommand{\mth}{\mt \mh}
\newcommand{\mhs}{\mh \ms}
\newcommand{\mtk}{\mt \mk}
\newcommand{\mks}{\mk \ms}

\newcommand{\mths}{\mt \mh \ms}
\newcommand{\mhss}{\mh \mss}
\newcommand{\mtth}{\mtt \mh}

\newcommand{\mthss}{\mth \mss}
\newcommand{\mhsss}{\mh \msss}
\newcommand{\mtths}{\mtt \mhs}

\newcommand{\mssst}{\msss \mt}
\newcommand{\mssts}{\mss \mts}
\newcommand{\mstss}{\ms \mt \mss}
\newcommand{\mtsss}{\mt \msss}

\newcommand{\msstt}{\mss \mtt} 
\newcommand{\mtsst}{\mt \mss \mt}
\newcommand{\mstst}{\mst \mst}
\newcommand{\mstts}{\ms \mtt \ms}
\newcommand{\mtsts}{\mts \mts}

\newcommand{\msttt}{\ms \mttt} 
\newcommand{\mttst}{\mtt \ms \mt}
\newcommand{\mttts}{\mttt \ms} 
\newcommand{\mtstt}{\mt \ms \mtt}

\newcommand{\mttth}{\mttt \mh}
\newcommand{\mkss}{\mk \mss}
\newcommand{\mttk}{\mtt \mk}
\newcommand{\mtks}{\mt \mks}

\newcommand{\lift}[1]{\hat{#1}}

\newcommand{\di}{d}

\newcommand{\myemph}{\textit}

\title[]{On the formal theory of pseudomonads and pseudodistributive laws}
\author[]{Nicola Gambino}
\author[]{Gabriele Lobbia}

\date{January 20th, 2021}

\begin{document}

\begin{abstract}
We contribute to the  formal theory of pseudomonads, {\em i.e.} the analogue for pseudomonads of the formal theory of monads. In particular, we solve a problem posed by Lack by proving that, for every Gray-category~$\catK$, there is a  Gray-category~$\Psm(\catK)$ of pseudomonads, pseudomonad morphisms, pseudomonad transformations and pseudomonad modifications in~$\catK$. We then establish a triequivalence between~$\Psm(\catK)$ and the Gray-category of pseudomonads introduced by Marmolejo and give a simpler proof of the equivalence between pseudodistributive laws and liftings of pseudomonads to 2-categories of
pseudoalgebras. 
\end{abstract}

\maketitle

\section*{Introduction}

\subsection*{Context and motivation} 

Monads are one of the fundamental notions of category theory~\cite[Chapter~VI]{MacLaneS:catwm}. For example, they provide a homogeneous approach to
the study of categories of sets equipped with algebraic structure, such as groups and monoids~\cite{BarrM:toptt}. Furthermore, Beck's theorem on  distributive laws
between monads~\cite{BeckJ:disl} describes concisely the structure that is necessary and sufficient in order to combine two algebraic structures, so that 
the operations of one distribute over those of the other. For example, the monads for groups and for monoids can be combined via a distributive law to define the monad for rings. Subsequently, the formal theory of monads, originally introduced by Street~\cite{StreetR:fortm} and later developed further by Lack and Street~\cite{LackS:fortm}, has offered an elegant and mathematically efficient account of the theory of monads, starting from the observation that the notion of a monad can be defined within any 2-category 
(so that the usual notion is recovered by considering the 2-category of categories, functors and natural transformations). Among many other results, it provides
a characterisation of the existence of categories of Eilenberg-Moore algebras as a completeness property and, importantly for our purposes, a simple account of Beck's theorem on distributive laws.

In recent years, motivation from pure mathematics, {\em e.g.}~in the theory of operads~\cite{FioreM:carcbg,FioreM:relpsm,GambinoN:opebaf,GarnerR:polvpd}, and theoretical computer science, {\em e.g.}~in the study of variable binding~\cite{CattaniG:proomb,CurienP:opecdl,TanakaM:unicsb}, led to significant interest in pseudomonads~\cite{LackS:cohapm,MarmolejoF:docwsf,MarmolejoF:dislpII,WalkerC:dislva}, which are the counterparts of monads in 2-dimensional category theory, obtained by requiring the axioms for a monad to hold only up to coherent isomorphism rather than strictly~\cite{BungeM:cohera}. Here, one of the key issues has been the proof of a counterpart of Beck's theorem on distributive laws, which requires a satisfactory axiomatisation of the notion of a pseudostributive law
\cite{ChengE:psedl,MarmolejoF:dislp,MarmolejoF:cohplr,TanakaM:psedl,TanakaM:psedlav}, building on early work of Kelly~\cite{KellyG:cohtla} on
semi-strict distributive laws. This is a difficult question because such a notion necessarily involves complex coherence conditions. 

Given how the formal theory of monads offers a simple proof of Beck's theorem on distributive laws,   it seems natural to attack this problem by developing a formal theory of pseudomonads. In order to do this, however, one needs to face the challenge that, just as the formal theory of monads is formulated within 2-dimensional category theory~\cite{KellyG:revetc}, the formal theory of pseudomonads is formulated within 3-dimensional category theory~\cite{GordonR:coht,GurskiN:algtt}, which is notoriously hard.  In this setting, it is convenient to work with Gray-categories, \emph{i.e.}~semistrict tricategories~\cite[Section 4.8]{GordonR:coht}, which are easier to handle than tricategories, but sufficiently general for many purposes, since every tricategory is triequivalent to a Gray-category~\cite[Theorem 8.1]{GordonR:coht}. 

In spite of significant advances in the creation of a formal theory of pseudomonads in the works cited above, there are
still fundamental questions to be addressed. In particular, there is not yet a direct counterpart of
 the 2-category~$\Mnd(\catK)$ of monads, monad morphisms and monad transformations in a 
 2-category~$\catK$,  which is the the starting point of the formal theory of monads~\cite{StreetR:fortm}.
Filling this gap would involve the definition, for a  Gray-category $\catK$, of a 3-dimensional category 
$\Psm(\catK)$ having pseudomonads in $\catK$ as 0-cells, pseudomonad morphisms as 1-cells,
and appropriately defined pseudomonad transformations and pseudomonad modifications
as  2-cells and 3-cells, respectively. This issue was raised by Lack
in~\cite[Section 6]{LackS:cohapm}, who suspected that defining 
$\Psm(\catK)$ in this way would give rise only to a tricategory, not a Gray-category, and hence require lengthy verifications of the coherence 
conditions. For this reason, Lack preferred to define a Gray-category of
pseudomonads in $\cat{K}$ using the description of pseudomonads in $\catK$ as suitable lax
functors and developing parts of the theory using enriched category theory.

Yet another approach was taken earlier by Marmolejo in~\cite{MarmolejoF:dislp}, who introduced, for a Gray-category $\catK$, a Gray-category that we denote here $\Lift(\catK)$ to avoid confusion, that has pseudomonads in $\catK$ as 0-cells and liftings  of 1-cells, 2-cells and 3-cells of $\catK$ 
to 2-categories of pseudoalgebras as 1-cells, 2-cells and 3-cells, respectively. He then used~$\Lift(\catK)$ to introduce the notion of a lifting of a pseudomonad
and of a pseudodistributive law, proving  the fundamental result that pseudodistributive 
laws are equivalent liftings of pseudomonads are equivalent to ~\cite[Theorems 6.2, 9.3 and~10.2]{MarmolejoF:dislp}, thus obtaining an analogue of Beck's  result on distributive laws. Here, Marmolejo defined pseudodistributive laws
explicitly, giving nine coherence conditions for them~\cite{MarmolejoF:dislp}. Later, Marmolejo and Wood~\cite{MarmolejoF:cohplr} showed not only that   an additional tenth coherence condition, introduced by Tanaka~\cite{TanakaM:psedl},
can be derived from Marmolejo's conditions, but also that one of the original nine conditions 
 introduced by Marmolejo is derivable from the others, thus reducing the number of 
coherence axioms for a pseudodistributive law  to eight.

\subsection*{Main results.} The aim of this paper is take some further steps in the development of the formal theory of pseudomonads. In particular, our main contributions are the following:

\begin{itemize} 
\item Theorem~\ref{thm:psm-gray}, which answers the question raised in~\cite{LackS:cohapm} by showing that for every Gray-category $\catK$, there is a Gray-category $\Psm(\catK)$ of pseudomonads, pseudomonad morphisms, pseudomonad transformations and pseudomonad modifications in $\catK$; 
 \item Theorem~\ref{thm:equiv-lift-psm},  the analogue of a fundamental result of the formal theory of monads, asserting that $\Psm(\catK)$ is equivalent, in a suitable 3-categorical sense, to the Gray-category $\Lift(\catK)$; 
 \item Proposition \ref{prop:distr-law-as-obj}, recording that an object of~$\Psm(\Psm(\catK))$  is the same thing a pseudodistributive law in $\catK$;
\item a new, simpler proof of Marmolejo's theorem equivalence between pseudodistributive laws and liftings of pseudomonads to 2-categories of pseudoalgebras, given as the proof of Theorem~\ref{thm:distr-iff-lift}.
\end{itemize}

Theorem~\ref{thm:psm-gray} supports the definition of a pseudodistributive law of~\cite{MarmolejoF:dislp,MarmolejoF:cohplr}, since it allows us to show
that a pseudodistributive law is the same 
 thing as a pseudomonad in~$\Psm(\catK)$ (Proposition~\ref{prop:distr-law-as-obj}), as one would expect by analogy with the
situation in the formal theory of monads. Thanks to this observation,
we provide an interpretation of the complex coherence conditions for a pseudodistributive law in terms of the simpler ones, namely
those for a pseudomonad morphism, a pseudomonad transformation and a pseudomonad modification (see Table~\ref{tab:coherence-axioms} for details). 
This point of view allows us to give 
a principled presentation of the conditions for pseudodistributive laws of~\cite{MarmolejoF:dislp,TanakaM:psedl}, included in
Appendix~\ref{app:coh}, which hopefully provides a useful reference for future work in this area.
For the convenience of readers, we also describe how our formulation relates to the
ones of Marmolejo and of Tanaka (see Table~\ref{tab:comp-coh-axioms}).

Theorem~\ref{thm:equiv-lift-psm}, which establishes the equivalence between $\Psm(\catK)$ and $\Lift(\catK)$, does not seem to be part of the 
literature
(in part because its very statement requires the introduction of the 3-dimensional category~$\Psm(\catK)$, which is defined here for the first time), 
but extends existing results. In particular, 
the equivalence between pseudomonad morphisms and liftings of morphisms to categories of pseudoalgebras is proved in~\cite{MarmolejoF:cohplr}. 
 Related results appear also in~\cite{TanakaM:psedl}, but with important differences. First, the work carried out
therein is developed for the particular tricategory $\mathbf{2\text{-}Cat}_{\mathrm{psd}}$ of 2-categories, pseudofunctors,
pseudonatural transformations and modifications, rather than for a general tricategory or Gray-category. While that is
an important example ({\em cf.}~Remark~\ref{rem:kleisli}), restricting to a particular tricategory does not allow us to
exploit the various dualities that are essential to derive results in the formal theory. Secondly, the results obtained
therein focus on hom-2-categories of pseudomonad endomorphisms, {\em i.e.}~of the form $\Psm(\cat{K})((X,S), (X,S))$, rather than on general hom-2-categories of pseudomonad morphisms.

Our proof Theorem~\ref{thm:distr-iff-lift}, which
asserts the equivalence between pseudodistributive laws and liftings of pseudomonads to 2-categories of pseudoalgebras 
established in~\cite{MarmolejoF:dislp}, follows naturally combining Theorem~\ref{thm:psm-gray} and Theorem~\ref{thm:equiv-lift-psm}. 
More specifically, combining our identification of pseudodistributive laws with pseudomonads in $\Psm(\catK)$ of
Proposition~\ref{prop:distr-law-as-obj} with the fact that a pseudomonad in~$\Lift(\catK)$ is a lifting of a pseudomonad~$T$ to the 2-categories of pseudoalgebras of another pseudomonad~$S$, we obtain the desired equivalence between a pseudodistributive law of~$S$ over~$T$ and a lifting of~$T$ to pseudo-$S$-algebras.
This proof is simpler than that in~\cite{MarmolejoF:dislp} since it takes a modular, more abstract, approach to the verification of the coherence conditions
and avoids completely the notion of a composite of pseudomonads with compatible structure.

 As the proofs of our main results involve lengthy, subtle calculations with pasting diagrams,
 we tried to strike a reasonable compromise between rigour and conciseness by giving what we hope are
 the key diagrams of the proofs, and describing the additional steps in the text. When in doubt, we preferred
 to err on the side of rigour, since one of our initial goals was to answer the question raised in~\cite{LackS:cohapm} about whether $\Psm(\catK)$ is a Gray-category or not. We hope that this
 did not make the paper too long. For the convenience of the readers, some of the diagrams are confined to the Appendices.
  
 \subsection*{Outline of the paper} Section~\ref{sec:grayc} provides background on Gray-categories,
 pseudomonads and 2-categories of pseudoalgebras. Section~\ref{sec:psem} defines the
 Gray-category~$\Psm(\catK)$. We  prove the equivalence of~$\Psm(\catK)$ and~$\Lift(\catK)$
 in Section~\ref{sec:lift}. We conclude the paper in Section~\ref{sec:psel} by discussing pseudodistributive laws.
 
 \subsection*{Acknowledgements} We are very grateful to Steve Lack for patiently discussing his
work on pseudomonads with us, especially during the 2019 International Conference on Category Theory at the University of Edinburgh.
 We wish to thank also John Bourke, Emily Riehl and Nick Gurski for fruitful conversations,
 Martin Hyland for encouraging us to consider carefully the directions of 3-cells, and the anonymous referee for helpful comments.
Nicola Gambino gratefully acknowledges the support of EPSRC under grant EP/M01729X/1, of the US Air Force Office for Scientific Research under agreement FA8655-13-1-3038 and the hospitality of the School of Mathematics of the University of Manchester while on study leave from the University of Leeds in 2018/19. 
Gabriele Lobbia gratefully acknowledges the support of an EPSRC PhD Scholarship.

\section{Preliminaries}
\label{sec:grayc}

\subsection*{Gray-categories} 

We begin by reviewing the notion of a Gray-category and fixing some notation. A Gray-category can be defined very
succinctly in terms of enriched category theory (see Remark~\ref{thm:gray-enriched}). For our purposes, however, 
it is useful to give an explicit definition, which we recall from~\cite[Section~2]{MarmolejoF:dislp} in Definition~\ref{defn:gray-cat} below. 
The explicit definition makes it easier to see that Gray-categories are special tricategories~\cite[Proposition~3.1]{GordonR:coht} in
which the only non-strict operation is horizontal composition of 2-cells~\cite[Section~5.2]{GordonR:coht}.
Throughout this paper, for a Gray-category $\catK$, we use $X, Y, Z, \ldots$ to denote its 0-cells, 
$F \colon X \to Y$, $G \colon Y \to Z$, \ldots for its 1-cells, $f \colon F \to F' \, , g \colon G \to G' \, \ldots$ for 2-cells, and $\alpha \colon f \to f' \, , \beta \colon g \to g' \, \ldots$ for 3-cells. 

When stating the definition of a Gray-category below, we make use of the notion of a cubical functor from~\cite{GordonR:coht}, which we unfold in Remark~\ref{rmk:cub-comp}.

\begin{defn}\label{defn:gray-cat} A Gray-category $\catK$ consists of the the data in (G1)-(G4),
subject to axioms (G5) and (G6), as given below.
\begin{enumerate}[(G1)]
\item A class of objects $\catK_0$. We call the elements of $\catK_0$ the 0-cells of $\catK$.
\item\label{ax:strict-hom} For every $X \, , Y \in \catK_0$, a 2-category $\catK(X \, , Y)$. We refer
to the $n$-cells of these 2-categories as the $(n+1)$-cells of~$\catK$. 
\item For every $X,\,Y,\,Z\in \catK_0$, a cubical functor 
\[
\catK(Y,\,Z)\times \catK(X,\,Y)\rightarrow \catK(X,\,Z) \, , 
\]
whose action on $F \colon X \to Y$ and $G \colon Y \to Z$
is written $GF \colon X \to Z$, and whose action on $f \colon F \to F'$ and $g \colon G \to G'$ gives rise
to an invertible 3-cell
\begin{displaymath}
\begin{tikzcd}
GF \arrow[r, "Gf", ""{name=U, inner sep=6pt,below}] \arrow[d, "gF"{left}] & GF' \arrow[d, "gF'"] \\
G'F \arrow[r, "G'f"{below}, ""{name=D, inner sep=6pt, above}] & G'F' 
\arrow[Rightarrow, from=U, to=D, "g_f"]
\end{tikzcd}
\end{displaymath}
called the \emph{interchange maps} of $\catK$. 
\item For any $X\in\catK_0$, a 1-cell $1_X \colon X \to X$. We call these the identity 1-cells of $\catK$.
\item\label{equ:gray} For every 
\[
\begin{tikzcd}
F \arrow[r, bend left, "f", ""{name=U, inner sep=1pt,below}] \arrow[r, bend right, "f'"{below}, ""{name=D, inner sep=1pt, above}] & F'
\arrow[Rightarrow, from=U, to=D, "\alpha"]         
\end{tikzcd} \qquad 
\begin{tikzcd}
G \arrow[r, bend left, "g", ""{name=U, inner sep=1pt,below}] \arrow[r, bend right, "g'"{below}, ""{name=D, inner sep=1pt, above}] & G'
\arrow[Rightarrow, from=U, to=D, "\beta"]         
\end{tikzcd} \qquad 
\begin{tikzcd}
K \arrow[r, bend left, "k", ""{name=U, inner sep=1pt,below}] \arrow[r, bend right, "k'"{below}, ""{name=D, inner sep=1pt, above}] & K'
\arrow[Rightarrow, from=U, to=D, "\gamma"]         
\end{tikzcd} 
\]
in $\catK(X,\,Y)$, $\catK(Y,\,Z)$ and $\catK(Z,\,W)$, respectively,  
\begin{center}
$(KG)F=K(GF) \, , $  \\
\vspace{0.15cm}
$(KG)f=K(Gf) \, , \quad (Kg)F=K(gF) \, , \quad (kG)F=k(GF) \, ,$ \ \\
\vspace{0.15cm}
$(KG)\alpha =K(G\alpha) \, , \quad (K\beta)F=K(\beta F) \, , \quad (\gamma G)F=\gamma(GF) \, , $ \\
\vspace{0.15cm}
$(Kg)_f=K(g_f) \, , \quad (kG)_f=k_{Gf} \, , \quad (k_g)F=k_{(gF)} \, .$  \\
\end{center}
\item\label{ax-gray-id} For every $X$, the 2-functors 
\[
1_X(-) \colon \catK(X,Y) \to \catK(X,Y) \, , \quad 
(-)1_X \colon \catK(X,Y) \to \catK(X,Y)
\] 
defined by composition with $1_X \colon X \to X$, are identities.
\end{enumerate}
\end{defn}

\begin{rmk}\label{rmk:cub-comp} Asserting that composition in a Gray-category $\catK$ is a cubical functor means that the properties in (i)-(v) below hold,
for every $F \, , F' \, , F'' \colon X \to Y$, $G \, , G'  \, , G'' \colon Y \to Z$ and
\[
\begin{tikzcd}
F \arrow[r, bend left, "f", ""{name=U, inner sep=1pt,below}] \arrow[r, bend right, "f'"{below}, ""{name=D, inner sep=1pt, above}] & F' \arrow[r, "f''"] & F'' \, , 
\arrow[Rightarrow, from=U, to=D, "\phi"]         
\end{tikzcd}
\qquad
\begin{tikzcd}
G \arrow[r, bend left, "g", ""{name=U, inner sep=1pt,below}] \arrow[r, bend right, "g'"{below}, ""{name=D, inner sep=1pt, above}] & G' \arrow[r, "g''"] & G'' \, .
\arrow[Rightarrow, from=U, to=D, "\psi"]         
\end{tikzcd}
\]
\begin{enumerate}[(i)] 
\item Composition with 1-cells on either side, 
\[
(-)F \colon \catK(Y,\,Z)\rightarrow\catK(X,\,Z) \,  \quad G(-):\catK(X,\,Y)\rightarrow\catK(X,\,Z) \, , 
\]
is a strict 2-functor.
\item Composition with 2-cells, 
\[
(-)f:(-)F\rightarrow(-)F' \, , \quad g(-):G(-)\rightarrow G'(-) \, , 
\]
 is a  pseudo-natural transformation. 
\item Composition with 3-cells, 
\[
(-)\varphi \colon (-)f\rightarrow(-)f' \, , \quad \psi (-) \colon g(-)\rightarrow g'(-) \, , 
\]
is a modification, 
\item The following coherence equations hold: 
\begin{equation}
\label{equ:cc1}
\begin{gathered}
\begin{tikzpicture}
\node (a) at (0.5,3) {$GF$};
\node (b) at (3.5,3) {$GF'$};
\node (c) at (0.5,0) {$G'F$};
\node (d) at (3.5,0) {$G'F'$};
\draw[->] (a) to [bend left] node[scale=.7] (f) [above] {$Gf$} (b);
\draw[->] (a) to [bend right] node[scale=.7] (g) [below] {$Gf'$} (b);
\draw[->] (a) to [bend left] node[scale=.7] (Gf) [right] {$gF$} (c);
\draw[->] (a) to [bend right] node[scale=.7] (Gf') [left] {$g'F$} (c);
\draw[->] (b) to node[scale=.7] (G'f) [right]{$gF'$} (d);
\draw[->] (c) to node[scale=.7] (g'F) [below]{$G'f'$} (d);
\draw[-{Implies},double distance=1.5pt,shorten >=10pt,shorten <=10pt] (f) to node[scale=.7] [right] {$G\varphi$} (g);
\draw[-{Implies},double distance=1.5pt,shorten >=40pt,shorten <=40pt, pos=0.7] (b) to node[scale=.7] [right, xshift=0.8cm, yshift=0.4cm] {$g_{f'}$} (c);
\draw[-{Implies},double distance=1.5pt,shorten >=10pt,shorten <=10pt] (Gf) to node[scale=.7] [below] {$\gamma F$} (Gf');
\node at (5,1.5) {$=$};
\node (a2) at (6.5,3) {$GF$};
\node (b2) at (9.5,3) {$GF'$};
\node (c2) at (6.5,0) {$G'F$};
\node (d2) at (9.5,0) {$G'F'$};
\draw[->] (c2) to [bend left] node[scale=.7] (f2) [above] {$Gf$} (d2);
\draw[->] (c2) to [bend right] node[scale=.7] (g2) [below] {$G'f'$} (d2);
\draw[->] (b2) to [bend left] node[scale=.7] (Gf2) [right] {$gF'$} (d2);
\draw[->] (b2) to [bend right] node[scale=.7] (Gf'2) [left] {$g'F'$} (d2);
\draw[->] (a2) to node[scale=.7] (G'f2) [left]{$g'F$} (c2);
\draw[->] (a2) to node[scale=.7] (g'F2) [above]{$G'f$} (b2);
\draw[-{Implies},double distance=1.5pt,shorten >=10pt,shorten <=10pt] (f2) to node[scale=.7] [right] {$G'\varphi$} (g2);
\draw[-{Implies},double distance=1.5pt,shorten >=40pt,shorten <=40pt, pos=0.7] (b2) to node[scale=.7] [right, yshift=1cm] {$g'_f$} (c2);
\draw[-{Implies},double distance=1.5pt,shorten >=10pt,shorten <=10pt] (Gf2) to node[scale=.7] [above] {$\gamma F'$} (Gf'2);
\end{tikzpicture} 
\end{gathered}
\end{equation}

\begin{equation}
\label{equ:cc2}
\begin{gathered} 
\begin{tikzpicture}
\node (a) at (3,4) {$GF$};
\node (b) at (5,4) {$GF'$};
\node (c) at (3,2) {$G'F$};
\node (d) at (5,2) {$G'F'$};
\node (e) at (3,0) {$G''F'$};
\node (f) at (5,0) {$G''F'$};
\draw[->] (a) to node[scale=.7] (gF) [above]{$Gf$} (b);
\draw[->] (a) to node[scale=.7] (Gf) [left]{$gF$} (c);
\draw[->] (c) to node[scale=.7] (gF') [below] {$G'f$} (d);
\draw[->] (b) to node[scale=.7] (G'f) [right] {$gF'$} (d);
\draw[->] (c) to node[scale=.7] (Gf'') [left] {$g''F$} (e);
\draw[->] (d) to node[scale=.7] (G'f'') [right] {$g''F'$} (f);
\draw[->] (e) to node[scale=.7] (gF'') [below] {$G''f$} (f);
\draw[-{Implies},double distance=1.5pt,shorten >=20pt,shorten <=20pt] (gF) to node[scale=.7] [right] {$g_f$} (gF');
\draw[-{Implies},double distance=1.5pt,shorten >=20pt,shorten <=10pt] (gF') to node[scale=.7] [right] {${g''}_f$} (gF'');
\node at (7,2) {$=$};
\node (a2) at (8.8,4) {$GF$};
\node (b2) at (10.8,4) {$GF'$};
\node (e2) at (8.8,0) {$G''F'$};
\node (f2) at (10.8,0) {$G''F'$};
\draw[->] (a2) to node[scale=.7] (gF2) [above]{$Gf$} (b2);
\draw[->] (a2) to node[scale=.7] (Gf''f2) [left]{$(g''g)F$} (e2);
\draw[->] (b2) to node[scale=.7] (G'f2) [right] {$(g''g)F'$} (f2);
\draw[->] (e2) to node[scale=.7] (gF''2) [below] {$G''f$} (f2);
\draw[-{Implies},double distance=1.5pt,shorten >=45pt,shorten <=45pt] (gF2) to node[scale=.7] [right] {${(g''g)}_f$} (gF''2);
\end{tikzpicture} 
\end{gathered}
\end{equation}

\begin{equation}
\label{equ:cc3}
\begin{gathered} 
\begin{tikzpicture}
\node (a) at (0,2) {$GF$};
\node (b) at (2,2) {$GF'$};
\node (c) at (0,0) {$G'F$};
\node (d) at (2,0) {$G'F'$};
\node (e) at (4,2) {$GF''$};
\node (f) at (4,0) {$G'F''$};
\draw[->] (a) to node[scale=.7] (gF) [above]{$Gf$} (b);
\draw[->] (a) to node[scale=.7] (Gf) [left]{$gF$} (c);
\draw[->] (c) to node[scale=.7] (gF') [below] {$G'f$} (d);
\draw[->] (b) to node[scale=.7] (G'f) [right] {$gF'$} (d);
\draw[->] (b) to node[scale=.7] (g''F) [above] {$Gf''$} (e);
\draw[->] (d) to node[scale=.7] (g''F') [below] {$G'f''$} (f);
\draw[->] (e) to node[scale=.7] (G''f) [right] {$gF''$} (f);
\draw[-{Implies},double distance=1.5pt,shorten >=20pt,shorten <=20pt] (gF) to node[scale=.7] [right, xshift=0.2cm] {$g_{f}$} (gF');
\draw[-{Implies},double distance=1.5pt,shorten >=20pt,shorten <=20pt] (g''F) to node[scale=.7] [right, xshift=0.2cm] {$g_{f''}$} (g''F');
\node at (5,1) {$=$};
\node (a2) at (6,2) {$GF$};
\node (b2) at (10,2) {$GF''$};
\node (e2) at (6,0) {$G'F$};
\node (f2) at (10,0) {$G'F''$};
\draw[->] (a2) to node[scale=.7] (gF2) [above]{$G(f''f)$} (b2);
\draw[->] (a2) to node[scale=.7] (Gf''f2) [left]{$gF$} (e2);
\draw[->] (b2) to node[scale=.7] (G'f2) [right] {$gF''$} (f2);
\draw[->] (e2) to node[scale=.7] (gF''2) [below] {$G'(f''f)$} (f2);
\draw[-{Implies},double distance=1.5pt,shorten >=20pt,shorten <=20pt] (gF2) to node[scale=.7] [right, xshift=0.2cm] {$g_{(f''f)}$} (gF''2);
\end{tikzpicture} 
\end{gathered}
\end{equation}
\item The interchange map $g_f$ is the identity 3-cell when either $f$ or $g$ is the identity.
\end{enumerate}
\end{rmk} 

\begin{rmk} \label{gray-notation}
When working with a Gray-category,  we sometimes write $G \circ F$ instead of $GF$ for 
cubical composition of $1$-cells. For 2-cells, we write 
$g' \cdot g$ (or $g' \co g$)
for the vertical composition and $g\circ f$ for cubical composition. For 3-cells, we write $\beta \circ\alpha$ for cubical composition, $\alpha' \ast \alpha$ for  vertical composition 
in $\catK(X,\,Y)$ and $\bar{\alpha}\cdot\alpha$ for horizontal composition in $\catK(X,\,Y)$, where $\alpha' \colon f'\rightarrow f''$ and $\bar{\alpha}\in\catK(X,\,Y)[F',\,F'']$. 
\end{rmk}

\begin{rmk} \label{thm:gray-enriched}
We write $\Gray$ for the  category of 2-categories and 2-functors. For 2-categories $X$ and $Y$,
let $[X,Y]$ be the 2-category of 2-functors from $X$ to $Y$, pseudonatural transformations, and modifications~\cite{KellyG:revetc}. This definition equips the category~$\Gray$ with 
the structure of a closed category~\cite{EilenbergS:cloc}. The closed structure of~$\Gray$  is part of symmetric monoidal structure, 
whose the tensor product is known as the \myemph{Gray tensor product}~\cite[Section~4.8]{GordonR:coht}. We will write~$X \otimes Y$ for the Gray tensor product of 2-categories~$X$ and~$Y$.
 A Gray-category can then be defined equivalently as a   $\Gray$-enriched category~\cite[Section~5.1]{GordonR:coht}. Since $\Gray$ is a monoidal closed category, it is enriched over itself. 
Therefore, it can be viewed as a Gray-category, as we will do from now on. More explicitly, 
$\Gray$ is the Gray-category having 2-categories as 0-cells, 2-functors as 1-cells, pseudonatural transformations as 2-cells, and modifications as 3-cells. 
\end{rmk}

The notions of a {\em Gray-functor} and of a {\em Gray-natural transformation} are instances of the general notions of enriched functor and enriched natural transformation~\cite[Section~1.2]{KellyG:bascec}.  We will use the terminology of {\em Gray-modification} and {\em Gray-perturbation} to denote the strict counterparts of the corresponding tricategorical notions~\cite[Section~3.3]{GordonR:coht}. 

When working with the Yoneda embedding for Gray-categories, which is just an instance of the Yoneda embedding for enriched categories~\cite[Section~2.4]{KellyG:bascec}, we often identify an object $X \in \catK$ with the representable Gray-functor $\catK(-, X) \colon \catK^\op \to \Gray$ associated to it. Analogous conventions will be used also for the $n$-cells of~$\catK$, where~$n = 1, 2, 3$.  For further information on Gray-categories and tricategories, we invite the reader to refer to~\cite{GarnerR:lowdsf,GordonR:coht,GurskiN:algtt,LackS:bicntg}.

\subsection*{Pseudomonads and their pseudoalgebras} Let $\catK$ be a Gray-category, to be considered fixed for the rest of this section. 
We recall the definition of a  pseudomonad.

\begin{defn}
\label{defn:0-cells-psm} Let $X\in\catK$. A \emph{pseudomonad} on $X$ in $\catK$ consists of:
\begin{itemize}
\item a $1$-cell $S:X\rightarrow X$ in $\catK$;
\item two 2-cells $\mus:S^2\rightarrow$ and $s:1_X\rightarrow S$ in $\catK$;
\item three invertible $3$-cells in $\catK$ of the following form:
\begin{center}
\begin{tikzpicture}
\node (a) at (0,2) {$S^3$};
\node (b) at (2,2) {$S^2$};
\node (e) at (0,0) {$S^2$};
\node (f) at (2,0) {$S$};
\draw[->] (a) to node[scale=.7] (gF) [above]{$S\mus$} (b);
\draw[->] (a) to node[scale=.7] (Gf) [left]{$\mus S$} (e);
\draw[->] (b) to node[scale=.7] (G'f) [right] {$\mus$} (f);
\draw[->] (e) to node[scale=.7] (gF'') [below] {$\mus$} (f);
\draw[-{Implies},double distance=1.5pt,shorten >=20pt,shorten <=20pt] (gF) to node[scale=.7] [right] {$\mu$} (gF'');

\node (a2) at (5,2) {$S$};
\node (b2) at (7,2) {$S^2$};
\node (f2) at (7,0) {$S$};
\node (g2) at (9,2) {$S$};
\node (x) at (5,1) {};
\node (x1) at (7,1) {};
\node (x2) at (9,1) {};
\draw[->] (a2) to node[scale=.7] (gF2) [above]{$Ss$} (b2);
\draw[->] (a2) to [bend right] node[scale=.7] (Gf''f2) [left]{$1_S$} (f2);
\draw[->] (b2) to node[scale=.7] (G'f2) [right] {$\mus $} (f2);
\draw[->] (g2) to [bend left] node[scale=.7] (Gf''f2) [right]{$1_S$} (f2);
\draw[->] (g2) to node[scale=.7] (gF3) [above]{$sS$} (b2);
\draw[-{Implies},double distance=1.5pt,shorten >=20pt,shorten <=20pt] (x1) to node[scale=.7] [above] {$\lambda$} (x);
\draw[-{Implies},double distance=1.5pt,shorten >=20pt,shorten <=20pt] (x2) to node[scale=.7] [above] {$\rho$} (x1);
\end{tikzpicture}
\end{center} 
satisfying the coherence axioms in~\eqref{equ:psm1} and~\eqref{equ:psm2} below:

\begin{equation}
\label{equ:psm1}
\begin{xy}  0;/r.28pc/:
(0,15)*+{\mssss}="1";
(20,15)*+{\msss}="2";
(15,0)*+{\msss}="5";
(35,0)*+{\mss}="6";
(0,-5)*+{\msss}="8";
(15,-20)*+{\mss}="9"; 
(35,-20)*+{\ms}="10";
(16,12)*+{}="11";
(20,6)*+{}="12";
(25,-7)*+{}="15";
(25,-14)*+{}="16";
(7,0)*+{}="17";
(7,-6)*+{}="18";
{\ar^{\mss \mus} "1";"2"};
{\ar_{\mus \mss} "1";"8"};
{\ar^{\ms \mus} "5";"6"};
{\ar^{\mus} "6";"10"};
{\ar^{\ms \mus \ms} "1";"5"};
{\ar^{\ms \mus} "2";"6"};
{\ar^{\mus \ms} "5";"9"};
{\ar_{\mus \ms} "8";"9"};
{\ar_{\mus} "9";"10"};
{\ar@{=>}^(.6){\ms \mu} "11";"12"};
{\ar@{=>}^{\mu} "15";"16"};
{\ar@{=>}^{\mu \ms} "17";"18"}
\end{xy}
 = 
\begin{xy} 0;/r.28pc/:
(0,15)*+{\mssss}="1";
(20,15)*+{\msss}="2";
(35,0)*+{\mss}="6";
(8,8)*+{}="15";
(8,4)*+{}="16";
{\ar@{=>}^{{\mus}_{\mus}} "15";"16"};
(0,-5)*+{\msss}="8";
(15,-20)*+{\mss}="9"; 
(35,-20)*+{\ms}="10";
(20,-5)*+{\mss}="19";
(7,0)*+{}="17";
(7,-6)*+{}="18";
(16,-8)*+{}="30";
(21,-14)*+{}="31";
(26,0)*+{}="32";
(26,-5)*+{}="33";
{\ar@{=>}_{\mu} "30";"31"};
{\ar^{\mss \mus} "1";"2"};
{\ar^{\mus \ms} "2";"19"};
{\ar^{\mus} "19";"10"};
{\ar_{\mus \mss} "1";"8"};
{\ar^{\ms \mus} "8";"19"};
{\ar^{\mus} "6";"10"};
{\ar^{\ms \mus} "2";"6"};
{\ar_{\mus \ms} "8";"9"};
{\ar^{\mus} "9";"10"};
{\ar@{=>}^{\mu} "32";"33"};
\end{xy}
\end{equation}

\begin{equation}
\label{equ:psm2}
\begin{xy}
(0,15)*+{\mss}="1";
(15,0)*+{\msss}="5";
(40,0)*+{\mss}="6";
(15,-20)*+{\mss}="9"; 
(40,-20)*+{\ms}="10";
(28,-7)*+{}="15";
(28,-14)*+{}="21";
{\ar@/^1pc/^{1_{\mss}} "1";"6"};
{\ar@/_1pc/_{1_{\mss}} "1";"9"};
{\ar^{\mus} "6";"10"};
{\ar^{\ms \mus} "5";"9"};
{\ar^{\ms \es \ms} "1";"5"};
{\ar^{\ms \mus} "5";"6"};
{\ar_{\mus} "9";"10"};
{\ar@{=>}^{\mu} "15";"21"};
(20,9)*+{}="30";
(20,4)*+{}="31";
{\ar@{=>}^{\ms \rho} "30";"31"};
(9,-2)*+{}="36";
(9,-7)*+{}="37";
{\ar@{=>}^{\lambda \ms} "36";"37"};
\end{xy}
\qquad = \qquad  
\begin{xy} 
(-20,0)*+{\mss}="1";
(0,0)*+{\ms}="10";
{\ar^{\mus} "1";"10"};
\end{xy}
\end{equation} 
\end{itemize}
\end{defn}

For brevity, we will refer to an object $X \in \catK$ and a pseudomonad~$(S, m, s, \mu, \lambda, \rho)$ on~$X$ simply as a pseudomonad in $\catK$ and
write simply $(X, S)$ to denote it.  

Note that the notion of a pseudomonad is self-dual, in the sense that a pseudomonad
in~$\catK$ is the same thing as a pseudomonad in~$\catK^\op$, where~$\catK^\op$ is
the Gray-category  obtained from~$\catK$ by reversing the direction of the 1-cells, but not that of the 2-cells and 3-cells. 
As in the formal theory of monads, this is important to obtain results by duality.

Let $(X,S)$ be a pseudomonad in $\catK$. For $I \in \catK$, there is a 2-category $\psSalg(I)$ 
of $I$-indexed pseudo-$S$-algebras, pseudoalgebra morphisms, and pseudoalgebra 2-cells,
whose definitions  we recall below. An {\em $I$-indexed pseudoalgebra} for $S$  
consists of a $1$-cell~$A \colon I \rightarrow X$, called the underlying 1-cell of the pseudoalgebra,
a 2-cell $a \colon \ms A \rightarrow A$, called the {\em structure map} of the pseudoalgebra, and invertible 3-cells 
\begin{equation*}
\begin{xy}
(-12,8)*+{\mss A}="4";
(-12,-8)*+{\ms A }="3";
(12,8)*+{\ms A}="2";
(12,-8)*+{A \, , }="1"; 
(-1,2)="5";
(-1,-3)="6";
{\ar^{\ms a} "4";"2"};
{\ar_{a} "3";"1"};
{\ar_{m_A} "4";"3"};
{\ar^{a} "2";"1"};
{\ar@{=>}^{\ \scriptstyle{\bar{a}}} "5";"6"};
\end{xy}
\hspace{1.5cm}
\begin{xy}
(-8,8)*+{A}="4";
(8,8)*+{\ms A}="2";
(8,-8)*+{A \, ,}="1";
(5,1)="5";
(0,1)="6";
{\ar@/_/_{1_A} "4";"1"};
{\ar^(.45){\es_A} "4";"2"};
{\ar^{a} "2";"1"};
{\ar@{=>}^(.4){{\scriptstyle{\tilde{a}}}} "6";"5"};
\end{xy} 
\end{equation*}
called the {\em associativity} and {\em unit}  of the pseudoalgebra, satisfying the coherence 
axioms~\eqref{equ:alg1} and~\eqref{equ:alg2} stated below.

\begin{equation}
\label{equ:alg1}
\begin{xy}
(0,15)*+{\msss A}="1";
(25,15)*+{\mss A}="2";
(15,0)*+{\mss A}="5";
(40,0)*+{\ms A }="6";
(0,-5)*+{\mss A}="8";
(15,-20)*+{\ms A}="9"; 
(40,-20)*+{A}="10";
(17,12)*+{}="11";
(23,6)*+{}="12";
(28,-7)*+{}="15";
(28,-14)*+{}="16";
(6,0)*+{}="17";
(6,-6)*+{}="18";
{\ar^{\mss a} "1";"2"};
{\ar_{\mus_{\ms A}} "1";"8"};
{\ar^{\ms a} "5";"6"};
{\ar^{a} "6";"10"};
{\ar^{\ms \mus_A} "1";"5"};
{\ar^{\ms a} "2";"6"};
{\ar^{\mus_A} "5";"9"};
{\ar_{\mus_A} "8";"9"};
{\ar_{a} "9";"10"};
{\ar@{=>}^(.65){ \ms \bar{a}} "11";"12"};
{\ar@{=>}^{\bar{a}} "15";"16"};
{\ar@{=>}^{\alpha_A} "17";"18"}
\end{xy}
 =  
\begin{xy}
(0,15)*+{\msss A }="1";
(25,15)*+{\mss A}="2";
(40,0)*+{\ms A}="6";
(12,9)*+{}="15";
(12,4)*+{}="16";
{\ar@{=>}^>>>>{\mus_a} "15";"16"};
(0,-5)*+{\mss A}="8";
(15,-20)*+{\ms A}="9"; 
(40,-20)*+{A \, .}="10";
(25,-5)*+{\ms A}="19";
(7,0)*+{}="17";
(7,-6)*+{}="18";
(15,-8)*+{}="30";
(20,-14)*+{}="31";
(33,-2)*+{}="32";
(33,-7)*+{}="33";
{\ar@{=>}^{\bar{a}} "30";"31"};
{\ar^{\mss a} "1";"2"};
{\ar_{\mus_A} "2";"19"};
{\ar_{a} "19";"10"};
{\ar_{\mus_{\ms A}} "1";"8"};
{\ar^{\ms a} "8";"19"};
{\ar^{a} "6";"10"};
{\ar^{\ms a} "2";"6"};
{\ar_{\mus_A} "8";"9"};
{\ar_{a} "9";"10"};
{\ar@{=>}_{\bar{a}} "32";"33"};
\end{xy}   
\end{equation} 

\begin{equation}
\label{equ:alg2}
\begin{xy}  
(0,15)*+{\ms A}="1";
(15,0)*+{\mss A}="5";
(40,0)*+{\ms A}="6";
(15,-20)*+{\ms A}="9"; 
(40,-20)*+{A}="10";
(28,-7)*+{}="15";
(28,-14)*+{}="21";
{\ar@/^1pc/^{1_{\ms A}} "1";"6"};
{\ar@/_1pc/_{1_{\ms A}} "1";"9"};
{\ar^{a} "6";"10"};
{\ar^{\mus_A} "5";"9"};
{\ar^{\ms \es_A} "1";"5"};
{\ar^{\ms a} "5";"6"};
{\ar_{a} "9";"10"};
{\ar@{=>}^{\bar{a}} "15";"21"};
(20,9)*+{}="30";
(20,4)*+{}="31";
{\ar@{=>}^{\; \tilde{a}} "30";"31"};
(9,-2)*+{}="36";
(9,-7)*+{}="37";
{\ar@{=>}^{\lambda_A} "36";"37"};
\end{xy}
=  
\begin{xy} 
(0,15)*+{};
(15,0)*+{\ms A}="1";
(40,0)*+{A \, .}="10";
{\ar^{a} "1";"10"};
(40,-20)*+{};
\end{xy} 
\end{equation} 

As usual, we refer to a pseudoalgebra by the name of its underlying 1-cell, leaving the rest of its data implicit. Similar conventions will be
implicitly assumed for other kinds of structures.

\begin{prop}[Marmolejo] Let  $(X,S)$ be a pseudomonad in $\catK$, $I \in \cat{K}$ and $A$ an~$I$-indexed pseudoalgebra for $S$. Then,
the coherence condition
\begin{equation}
\label{equ:alg3} 
\begin{xy}
(0,15)*+{\ms A}="1";
(25,15)*+{ A}="2";
(0,-5)*+{\mss A}="8";
 (15,-20)*+{\ms A}="9"; 
(12,9)*+{}="15";
(12,4)*+{}="16";
{\ar@{=>}^{\es_a} "15";"16"};
(40,-20)*+{A}="10";
(25,-5)*+{\ms A}="19";
(7,0)*+{}="17";
(7,-6)*+{}="18";
(15,-8)*+{}="30";
(20,-14)*+{}="31";
(33,-2)*+{}="34";
(33,-7)*+{}="35";
{\ar@{=>}^{\bar{a}} "30";"31"};
{\ar^{a} "1";"2"};
{\ar^{\es_A} "2";"19"};
{\ar_{a} "19";"10"};
{\ar_{\es_{\ms A}} "1";"8"};
{\ar^{\ms a} "8";"19"};
{\ar@/^1.5pc/^{1_{\ms A}} "2";"10"};
{\ar_{\mus_A} "8";"9"};
{\ar_{a} "9";"10"};
{\ar@{=>}_{\tilde{a}\, } "34";"35"};
\end{xy}
\   = \  
\begin{xy}
(0,15)*+{\ms A}="1";
(0,-5)*+{\mss A}="8";
(15,-20)*+{\ms A}="9";
(40,-20)*+{A}="10";
(7,-2)*+{}="36";
(7,-7)*+{}="37";
{\ar_{\es_{\ms A}} "1";"8"};
{\ar_{\mus_A} "8";"9"};
{\ar_{a}  "9";"10"};
{\ar@/^1.5pc/^{1_{\ms A}} "1";"9"};
{\ar@{=>}^{\, \rho_A} "36";"37"};
\end{xy}
\end{equation} 
is derivable.
\end{prop}

\begin{proof} See~\cite[Lemma~9.1]{MarmolejoF:docwsf}.
\end{proof} 

Given pseudoalgebras $A$ and $B$, a \myemph{pseudoalgebra morphism}  $f : A \rightarrow B$ consists of a 2-cell  $f :  A \rightarrow B$ and an  invertible 3-cell \begin{equation*}
\begin{xy}
(-12,8)*+{\ms A}="4";
(-12,-8)*+{A }="3";
(12,8)*+{\ms B}="2";
(12,-8)*+{B}="1"; 
(-1,2)="5";
(-1,-3)="6";
{\ar^{\ms f} "4";"2"};
{\ar_{f} "3";"1"};
{\ar_{a} "4";"3"};
{\ar^{b} "2";"1"};
{\ar@{=>}^{\ \scriptstyle{\bar{f}}} "5";"6"};
\end{xy}
\end{equation*}
satisfying the coherence conditions~\eqref{equ:mor1} and~\eqref{equ:mor2} stated below.
\begin{equation}
\label{equ:mor1}
\begin{xy}
(0,15)*+{\mss A}="1";
(25,15)*+{\mss B}="2";
(15,0)*+{\ms A}="5";
(40,0)*+{\ms B}="6";
(0,-5)*+{\ms A}="8";
(15,-20)*+{A}="9"; 
(40,-20)*+{B}="10";
(16,12)*+{}="11";
(22,6)*+{}="12";
(28,-7)*+{}="15";
(28,-14)*+{}="16";
(7,0)*+{}="17";
(7,-6)*+{}="18";
{\ar^{\mss f} "1";"2"};
{\ar_{\mus_A} "1";"8"};
{\ar^{\ms f} "5";"6"};
{\ar^{b} "6";"10"};
{\ar^{\ms a} "1";"5"};
{\ar^{\ms b} "2";"6"};
{\ar^{a} "5";"9"};
{\ar_{a} "8";"9"};
{\ar_{f} "9";"10"};
{\ar@{=>}^(.7){\ms \bar{f}} "11";"12"};
{\ar@{=>}^{\bar{f}} "15";"16"};
{\ar@{=>}^{\bar{a}} "17";"18"}
\end{xy}
 =  
\begin{xy}
(0,15)*+{\mss A}="1";
(25,15)*+{\ms B}="2";
(40,0)*+{\ms B}="6";
(0,-5)*+{\ms A}="8";
(15,-20)*+{A}="9"; 
(40,-20)*+{B \, .}="10";
(12,9)*+{}="15";
(12,4)*+{}="16";
{\ar@{=>}^>>>>>{\mus_f} "15";"16"};
(25,-5)*+{\ms B}="19";
(7,0)*+{}="17";
(7,-6)*+{}="18";
(15,-8)*+{}="30";
(20,-14)*+{}="31";
(33,0)*+{}="32";
(33,-5)*+{}="33";
{\ar@{=>}^{\, \bar{f}} "30";"31"};
{\ar^{\mss f} "1";"2"};
{\ar_{\mus_B} "2";"19"};
{\ar^(.45){b} "19";"10"};
{\ar_{\mus_A} "1";"8"};
{\ar^{\ms f} "8";"19"};
{\ar^{b} "6";"10"};
{\ar^{\ms b} "2";"6"};
{\ar_{a} "8";"9"};
{\ar_{f} "9";"10"};
{\ar@{=>}^{\, \bar{b}} "32";"33"};
\end{xy} 
\end{equation} 

\begin{equation} 
\label{equ:mor2}
\begin{xy}
(0,15)*+{A}="1";
(0,-5)*+{\ms A}="8";
(15,-20)*+{A}="9";
(40,-20)*+{B}="10";
(7,-2)*+{}="36";
(7,-7)*+{}="37";
{\ar_{\es_A} "1";"8"};
{\ar_{a} "8";"9"};
{\ar_{f}  "9";"10"};
{\ar@/^1pc/^{1_{A}} "1";"9"};
{\ar@{=>}_{\tilde{a}} "36";"37"};
\end{xy}
 =  
\begin{xy}
(0,15)*+{A}="1";
(25,15)*+{B}="2";
(0,-5)*+{\ms A}="8";
 (15,-20)*+{A}="9"; 
(12,9)*+{}="15";
(12,4)*+{}="16";
{\ar@{=>}^{\es_f} "15";"16"};
(40,-20)*+{B \, .}="10";
(25,-5)*+{\ms B}="19";
(7,0)*+{}="17";
(7,-6)*+{}="18";
(15,-8)*+{}="30";
(20,-14)*+{}="31";
(33,-2)*+{}="34";
(33,-7)*+{}="35";
{\ar@{=>}^{\bar{f}} "30";"31"};
{\ar^{f} "1";"2"};
{\ar^{\es_B} "2";"19"};
{\ar_{b} "19";"10"};
{\ar_{\es_A} "1";"8"};
{\ar^{\ms f} "8";"19"};
{\ar@/^1pc/^{1_{B}} "2";"10"};
{\ar_{a} "8";"9"};
{\ar_{f} "9";"10"};
{\ar@{=>}_{\tilde{b} \, } "34";"35"};
\end{xy}
\end{equation} 
Given pseudoalgebra morphisms $f : A \rightarrow B$ and $g : A \rightarrow B$, 
a  \myemph{pseudoalgebra $2$-cell} consists of a 3-cell $\alpha : f \rightarrow g$  
satisfying the coherence condition~\eqref{equ:2cell}. 

\begin{equation}
\label{equ:2cell} 
\begin{xy}
(-12,8)*+{\ms A}="4";
(-12,-12)*+{A }="3";
(12,8)*+{\ms B}="2";
(12,-12)*+{B}="1"; 
(0,-4)="5";
(0,-9)="6";
(0,10)="7";
(0,6)="8";
{\ar@/^1pc/^{\ms f} "4";"2"};
{\ar@/_1pc/_{\ms g} "4";"2"};
{\ar@/_1pc/_{g} "3";"1"};
{\ar_{a} "4";"3"};
{\ar^{b} "2";"1"};
{\ar@{=>}^{\ \scriptstyle{\bar{g}}} "5";"6"};
{\ar@{=>}^{\; \scriptstyle{S \alpha}} "7";"8"};
\end{xy}
\qquad = \qquad  
\begin{xy}
(-12,8)*+{\ms A}="4";
(-12,-12)*+{A }="3";
(12,8)*+{\ms B}="2";
(12,-12)*+{B \, .}="1"; 
(0,6)="5";
(0,0)="6";
(0,-10)="7";
(0,-14)="8";
{\ar@/^1pc/^{\ms f} "4";"2"};
{\ar@/^1pc/^{f} "3";"1"};
{\ar@/_1pc/_{g} "3";"1"};
{\ar_{a} "4";"3"};
{\ar^{b} "2";"1"};
{\ar@{=>}^{\ \scriptstyle{\bar{f}}} "5";"6"};
{\ar@{=>}^{\; \scriptstyle{\alpha}} "7";"8"};
\end{xy}
\end{equation} 
There is a forgetful 2-functor $U_I : \psSalg(I) \rightarrow \catK(I,X)$, defined by mapping a
pseudo-$S$-algebra to its underlying 1-cell, which has a left pseudoadjoint, defined
by mapping a 1-cell $A : I \rightarrow X$ to the free pseudoalgebra on it, given 
by the composite 1-cell~$SA : I \rightarrow X$. Attentive readers will have observed that the
directions of the structural 3-cells $\mu$ and~$\lambda$ for a pseudomonad  as in
Definition~\ref{defn:0-cells-psm} match those of the 3-cells necessary to make $SA$ into a pseudoalgebra. 

The function mapping an object $I \in \catK$ to the 2-category $\psSalg(I)$ extends to a Gray-functor
$\psSalg : \catK^\op \rightarrow \Gray$.  We also have a Gray-transformation 
\begin{equation}
\label{equ:forgetful} 
U : \psSalg \rightarrow X \, ,
\end{equation}
with components given by the forgetful 2-functors $U_I : \psSalg(I) \rightarrow \catK(I,X)$, for~$I \in$~$\catK$. 
Note the use of our convention on the Yoneda lemma in~\eqref{equ:forgetful}. 
Note that  the structure of pseudo-$S$-algebra on a 1-cell $A \colon I \to X$ can be viewed as a 
left $S$-action on $A$, associative and unital up to coherent isomorphism. For this reason, we
sometimes refer to pseudoalgebras as {\em left pseudomodules}. This terminology is convenient
when we discuss dualities in Section~\ref{sec:psel}.


\section{The Gray-category of pseudomonads}
\label{sec:psem} 

The aim of this section is to introduce the 3-dimensional category $\Psm(\catK)$ of pseudomonads in a Gray-category $\catK$ and prove that it is a Gray-category. In order to do so, we review the notion of a pseudomonad morphism from~\cite{MarmolejoF:cohplr} and introduce the notions of a pseudomonad transformation and modification. Again, we fix a Gray-category $\catK$. When working with two pseudomonads $(X, S)$ and $(Y, T)$, we use $m$ and $s$ for the multiplication and
unit of $S$, $n$ and $t$ for the multiplication and unit of $T$, but we use the same letters $\mu$, $\lambda$, $\rho$ for the structural 3-cells of both monads to
simplify notation, as the
context makes it always clear to which we are referring.

\begin{defn} Let $(X, S)$ and $(Y, T)$ be pseudomonads in $\catK$.
A \emph{pseudomonad morphism} $(F, \phi) \colon (X,\,S) \rightarrow (Y,\,T)$  consists of a 1-cell $F:X\rightarrow Y$, a $2$-cell $\phi:TF\rightarrow FS$ and two invertible $3$-cells
\begin{center}
\begin{tikzpicture}
\node (a) at (0,4) {$T^2F$};
\node (b) at (2,4) {$TFS$};
\node (d) at (2,2) {$FS^2$};
\node (e) at (0,0) {$TF$};
\node (f) at (2,0) {$FS$};
\draw[->] (a) to node[scale=.7] (gF) [above]{$T\phi$} (b);
\draw[->] (a) to node[scale=.7] (Gf) [left]{$\mut F$} (e);
\draw[->] (b) to node[scale=.7] (G'f) [right] {$\phi S$} (d);
\draw[->] (d) to node[scale=.7] (G'f) [right] {$F\mus $} (f);
\draw[->] (e) to node[scale=.7] (gF'') [below] {$\phi$} (f);
\draw[-{Implies},double distance=1.5pt,shorten >=45pt,shorten <=45pt] (gF) to node[scale=.7] [right] {$\bar{\phi}$} (gF'');

\node (a2) at (5,4) {$F$};
\node (b2) at (7,4) {$TF$};
\node (f2) at (7,0) {$FS$.};
\node (Y) at (5,2) {};
\node (Y') at (7,2) {};
\draw[->] (a2) to node[scale=.7] (gF2) [above]{$tF$} (b2);
\draw[->] (a2) to [bend right] node[scale=.7] (Gf''f2) [left]{$Fs$} (f2);
\draw[->] (b2) to node[scale=.7] (G'f2) [right][xshift=0.1cm] {$\phi$} (f2);
\draw[-{Implies},double distance=1.5pt,shorten >=15pt,shorten <=15pt] (Y) to node[scale=.7] [above] [yshift=0.2cm] {$\widetilde{\phi}$} (Y');
\end{tikzpicture}
\end{center} 
These data are required to satisfy the coherence axioms in~\eqref{equ:m1} and~\eqref{equ:m2}. 
\begin{multline}
\label{equ:m1}
\begin{xy}   0;/r.20pc/:
(0,25)*+{\mttth}="1";
(30,25)*+{\mtths}="2";
(45,12)*+{\mthss}="6";
(0,0)*+{\mtth}="3";
(30,0)*+{\mtth}="9";
(60,0)*+{\mths}="10";
(30,-30)*+{\mth}="11";
(60,-30)*+{\mhs}="12";
(60,-12)*+{\mhss}="13";
{\ar^{\mtt \phi} "1";"2"};
{\ar_{\mut \mth} "1";"3"};
{\ar^{\mt \mut \mh} "1";"9"};
{\ar^{\mt \phi \ms} "2";"6"};
{\ar_{\mut \mh} "3";"11"};
{\ar^{\mth \mus} "6";"10"};
{\ar_{\mt \phi} "9";"10"};
{\ar_{\mut \mh} "9";"11"};
{\ar^{\phi \ms} "10";"13"};
{\ar^{\mh \mus} "13";"12"};
{\ar_{h} "11";"12"};
{\ar@{}|{\Downarrow \; \mt \bar{\phi}} "1";"10"};
{\ar@{}|{\Downarrow \; \bar{\phi}} "9";"12"};
{\ar@{}|{\Downarrow \;  \mu \mh} "1";"11"};
\end{xy} = \qquad \\
\begin{xy}    0;/r.20pc/:
(0,30)*+{\mttth}="1";
(30,30)*+{\mtths}="2";
(0,-4)*+{\mtth}="3";
(30,-4)*+{\mths}="4";
(50,14)*+{\mthss}="6";
(12,6)*+{};
(50,-20)*+{\mhss}="8";
(50,-4)*+{\mhsss}="14";
(70,-2)*+{\mths}="10";
(40,-37)*+{\mth}="11";
(70,-37)*+{\mhs}="12";
(70,-20)*+{\mhss}="13";
(7,0)*+{}="17";
(7,-6)*+{}="18";
(15,-8)*+{}="30";
(20,-14)*+{}="31";
(33,-2)*+{}="32";
(33,-7)*+{}="33";
{\ar^{\mtt \phi} "1";"2"};
{\ar_{\mut \mth} "1";"3"};
{\ar^{\mut \mhs} "2";"4"};
{\ar^{\mt \phi \ms} "2";"6"};
{\ar^{\mt \phi} "3";"4"};
{\ar_{\mut \mh} "3";"11"};
{\ar_{\phi \ms} "4";"8"};
{\ar_(.4){\phi \mss} "6";"14"};
{\ar_(.4){\mh \mus \ms} "14";"8"};
{\ar^(.4){\mhs \mus} "14";"13"};
{\ar^{\mth \mus} "6";"10"};
{\ar_{\mh \mus} "8";"12"};
{\ar^{\phi \ms} "10";"13"};
{\ar^{\mh \mus} "13";"12"};
{\ar_{\phi} "11";"12"};
{\ar@{}|{\Downarrow \; \mut _\phi}  "1";"4"};
{\ar@{}|{\Downarrow \, \bar{\phi}S \quad }  "2";"8"};
{\ar@{}|{\Downarrow \; \phi_m }  "6";"13"};
{\ar@{}|{\Downarrow \; \bar{\phi}  }  "3";"12"};
{\ar@{}|{\Downarrow \; F \mu  }  "14";"12"};
\end{xy}
\end{multline}

\begin{equation}
\label{equ:m2}
\begin{xy}
(0,15)*+{\mth}="1";
(15,0)*+{\mtth}="5";
(40,0)*+{\mths}="6";
(40,-18)*+{\mhss}="7";
(15,-35)*+{\mth}="9"; 
(40,-35)*+{\mhs}="10";
(28,-7)*+{}="15";
(28,-14)*+{}="21";
(10,-10)*+{{}_{\Downarrow \lambda H}};
{\ar@/^1pc/^{\mth \es} "1";"6"};
{\ar@/_1pc/_{1_{\mth}} "1";"9"};
{\ar^{\phi \ms} "6";"7"};
{\ar^{\mh \mus} "7";"10"};
{\ar^{\mut \mh} "5";"9"};
{\ar^{\mt \et \mh} "1";"5"};
{\ar^{\mt \phi} "5";"6"};
{\ar_{\phi} "9";"10"};
{\ar@{}|{\Downarrow \, \mt \tilde{\phi}} "1";"6" };
{\ar@{}|{\Downarrow \,  \bar{\phi}}  "5";"10" };
(20,9)*+{}="30";
(20,4)*+{}="31";
(9,-2)*+{}="36";
(9,-7)*+{}="37";
\end{xy}
\quad = \quad 
\begin{xy} 
(0,15)*+{\mth}="1";
(0,-3)*+{\mhs}="8";
(40,0)*+{\mths}="6";
(40,-18)*+{\mhss}="7";
(40,-35)*+{\mhs}="10";
{\ar@/^1pc/^{\mth \es} "1";"6"};
{\ar_{\phi} "1";"8"};
{\ar^{\mhs \es} "8";"7"};
{\ar^{\phi \ms} "6";"7"};
{\ar^{\mh \mus} "7";"10"};
{\ar@/_1pc/_{1_{\mhs}} "8";"10"};
{\ar@{}|{\Downarrow \, \phi_s} "1";"7" };
{\ar@{}|{\Downarrow \, \mh \lambda} "8";"10" };
\end{xy}
\end{equation}

\end{defn}

\begin{prop}[Marmolejo and Wood]  \label{thm:marwood}
Let $(F, \phi) : (X,S) \rightarrow (Y,T)$ be a 
pseudomonad morphism. The coherence condition
\begin{equation*}
\label{equ:m3} 
\begin{xy} 
(45,-20)*+{\mth}="1";
(65,-20)*+{\mhs}="2";
(5,20)*+{\mth}="3";
(30,20)*+{\mhs}="4";
(5,0)*+{\mtth}="5";
(30,0)*+{\mths}="6";
(45,-8)*+{\mhss}="7";
(15,14)*+{}="8";
(15,8)*+{}="9";
(36,-5)*+{}="10";
(44,-11)*+{}="11";
(38,8)*+{}="12";
(38,2)*+{}="13";
(62,0)*+{}="14";
(62,-6)*+{}="15";
{\ar_{\phi} "1";"2"};
{\ar^{\phi} "3";"4"};
{\ar_{\et \mth} "3";"5"};
{\ar_{\mut \mh} "5";"1"};
{\ar^(.35){\mh \mus} "7";"2"};
{\ar^(.4){\phi \ms} "6";"7"};
{\ar_{\et \mhs} "4";"6"};
{\ar^{\mt \phi} "5";"6"};
{\ar@/^/^(.55){\mh \es \ms}  "4";"7"};
{\ar@/^2pc/^{1_{\mhs}} "4";"2"};
{\ar@{}|{\Downarrow \, t_\phi}"5";"4"};
{\ar@{}|{\Downarrow \, \bar{\phi}}"1";"6"};
{\ar@{}|{\quad \Downarrow \, \mh \rho_S }"2";"4"};
(34,6)*+{{}_{\Downarrow \, \tilde{\phi} \ms} };
\end{xy} = 
\begin{xy} 0;/r.22pc/: 
(25,-20)*+{\mth}="1";
(45,-20)*+{\mhs}="2";
(5,20)*+{\mth}="3";
(5,0)*+{\mtth}="5";
(12,5)*+{}="6";
(12,0)*+{}="7";
{\ar_{\phi} "1";"2"};
{\ar_{\et \mth} "3";"5"};
{\ar_{\mut \mh} "5";"1"};
{\ar@/^1.8pc/^{1_{\mth}} "3";"1"};
{\ar@{=>}^{\rho_T \mh} "6";"7"}
\end{xy} 
\end{equation*} 
is derivable. 
\end{prop} 

\begin{proof} See~\cite[Theorem 2.3]{MarmolejoF:cohplr}.
\end{proof} 

\begin{defn} Let $(F,\,\phi),\,(F',\,\phi'): (X,\,S)\rightarrow(Y,T)$ be  pseudomonad morphisms. A \emph{pseudomonad transformation} $(p,\,\bar{p}):(F,\,\phi)\rightarrow(F',\,\phi')$ consists of a 2-cell $p:F\rightarrow F'$ and an invertible $3$-cell 
\begin{center}
\begin{tikzpicture}
\node (a) at (0,2) {$TF$};
\node (b) at (2,2) {$TF'$};
\node (c) at (0,0) {$FS$};
\node (d) at (2,0) {$F'S$};
\draw[->] (a) to node[scale=.7] (Tp) [above]{$Tp$} (b);
\draw[->] (a) to node[scale=.7] (phisx) [left]{$\phi$} (c);
\draw[->] (b) to node[scale=.7] (phidx) [right]{$\phi'$} (d);
\draw[->] (c) to node[scale=.7] (pS) [below]{$pS$} (d);
\draw[-{Implies},double distance=1.5pt,shorten >=20pt,shorten <=20pt] (Tp) to node[scale=.7] [right] [xshift=0.2cm] {$\bar{p}$} (pS);
\end{tikzpicture}
\end{center}
satisfying the coherence conditions in~\eqref{equ:t1} and~\eqref{equ:t2} below:
\begin{equation}
\label{equ:t1}
\begin{xy}
(-20,25)*+{\mtth}="1";
(0,25)*+{\mttk}="2";
(-20,-5)*+{\mth}="3";
(0,5)*+{\mths}="5";
(20,5)*+{\mtks}="6";
(0,-10)*+{\mhss}="7";
(20,-10)*+{\mkss}="8";
(0,-25)*+{\mhs}="9";
(20,-25)*+{\mks}="10";
{\ar^{\mtt p} "1";"2"};
{\ar_{\mut \mh} "1";"3"};
{\ar^{\mt \phi} "1";"5"};
{\ar^{\mt \phi'} "2";"6"};
{\ar_{\phi} "3";"9"};
{\ar^{\mt p \ms} "5";"6"};
{\ar^{\phi \ms} "5";"7"};
{\ar^{\phi' \ms} "6";"8"};
{\ar^{p \mss} "7";"8"};
{\ar_<<<{\mh \mus} "7";"9"};
{\ar^{\mk \mus} "8";"10"};
{\ar_{p \ms} "9";"10"};
{\ar@{}|{\ \Downarrow \, \mt \bar{p} } "1";"6"};
{\ar@{}|{\ \Downarrow \, \bar{\phi}} "1";"9"};
{\ar@{}|{\ \Downarrow \, \bar{p} \ms } "5";"8"};
{\ar@{}|{\ \ \Downarrow \, p_{\mus }^{-1}} "7";"10"};
\end{xy}
\quad = 
\begin{xy}
(-20,25)*+{\mtth}="1";
(0,25)*+{\mttk}="2";
(-20,-5)*+{\mth}="3";
(0,-5)*+{\mtk}="4";
(20,5)*+{\mtks}="6";
(20,-10)*+{\mkss}="8";
(0,-25)*+{\mhs}="9";
(20,-25)*+{\mks}="10";
{\ar^{\mtt p} "1";"2"};
{\ar_{\mut \mh} "1";"3"};
{\ar^{\mut \mk} "2";"4"};
{\ar^{\mt \phi'} "2";"6"};
{\ar^{\mt p} "3";"4"};
{\ar_{\phi} "3";"9"};
{\ar^{\phi'} "4";"10"};
{\ar^{\phi' \ms} "6";"8"};
{\ar^{\mk \mus} "8";"10"};
{\ar_{p \ms} "9";"10"};
{\ar@{}|{\Downarrow \, \mut_p} "1";"4"};
{\ar@{}|{\Downarrow\, \bar{\phi'} } "2";"10"};
{\ar@{}|{\Downarrow \, \bar{p} } "3";"10"};
\end{xy}
\end{equation}

\begin{equation}
\label{equ:t2}
\begin{xy}  
(0,15)*+{\mh}="1";
(20,15)*+{\mk}="2";
(0,-5)*+{\mth}="8";
(20,-20)*+{\mhs}="9";
(40,-20)*+{\mks}="10";
(7,0)*+{}="36";
(7,-5)*+{}="37";
(25,-2)*+{}="38";
(25,-7)*+{}="39";
{\ar_{\et \mh} "1";"8"};
{\ar^{p} "1";"2"};
{\ar@/^1pc/^{\mk \es} "2";"10"};
{\ar_{\phi} "8";"9"};
{\ar_{p \ms}  "9";"10"};
{\ar@/^1pc/^{\mh \es} "1";"9"};
{\ar@{=>}^{\; \tilde{\phi}} "36";"37"};
{\ar@{=>}^{\; p_s^{-1}} "38";"39"};
\end{xy}
\quad = \quad 
\begin{xy}   
(0,15)*+{\mh}="1";
(20,15)*+{\mk}="2";
(0,-5)*+{\mth}="8";
 (20,-20)*+{\mhs}="9"; 
(40,-20)*+{\mks}="10";
(20,-5)*+{\mtk}="19";
(12,8)*+{}="17";
(12,4)*+{}="18";
(13,-8)*+{}="30";
(20,-14)*+{}="31";
(30,-2)*+{}="34";
(30,-7)*+{}="35";
{\ar@{=>}^(.65){\; \bar{p}} "30";"31"};
{\ar^{p} "1";"2"};
{\ar^{\et \mk} "2";"19"};
{\ar_{\phi'} "19";"10"};
{\ar_{\et \mh} "1";"8"};
{\ar^{\mt p} "8";"19"};
{\ar@/^1pc/^{\mk \es} "2";"10"};
{\ar_{\phi} "8";"9"};
{\ar_{p \ms} "9";"10"};
  {\ar@{=>}^{\; \tilde{\phi'}} "34";"35"};
 {\ar@{=>}^{\; t_p} "17";"18"}
\end{xy}
\end{equation}

\end{defn}

\begin{defn} Let $(p,\,\widetilde{p}),\,(p',\,\widetilde{p'}):(F,\,\phi)\rightarrow (F',\,\phi')$ be pseudomonad transformations. A \emph{pseudomonad modification} $\alpha:(p,\,\widetilde{p})\rightarrow(p',\,\widetilde{p'})$ is a $3$-cell $\alpha:p\rightarrow p'$ satisfying the  coherence condition
\begin{equation}
\label{equ:md}
\begin{xy}
(-12,8)*+{\mth}="4";
(-12,-12)*+{\mhs}="3";
(12,8)*+{\mtk}="2";
(12,-12)*+{\mks}="1"; 
(-1,-4)="5";
(-1,-9)="6";
(0,10)="7";
(0,6)="8";
{\ar@/^1pc/^{\mt p} "4";"2"};
{\ar@/_1pc/_{\mt p'} "4";"2"};
{\ar@/_1pc/_{p' \ms} "3";"1"};
{\ar_{\phi} "4";"3"};
{\ar^{\phi'} "2";"1"};
 {\ar@{=>}^{\ \scriptstyle{\bar{p'}}} "5";"6"};
{\ar@{=>}^{\; \scriptstyle{\mt \alpha}} "7";"8"};
\end{xy}
\qquad = \qquad  
\begin{xy}
(-12,8)*+{\mth}="4";
(-12,-12)*+{\mhs }="3";
(12,8)*+{\mtk}="2";
(12,-12)*+{\mks}="1"; 
(-1,6)="5";
(-1,0)="6";
(0,-10)="7";
(0,-14)="8";
{\ar@/^1pc/^{\mt p} "4";"2"};
{\ar@/^1pc/^{p \ms} "3";"1"};
{\ar@/_1pc/_{p' \ms} "3";"1"};
{\ar_{\phi} "4";"3"};
{\ar^{\phi'} "2";"1"};
{\ar@{=>}^{\ \scriptstyle{\bar{p}}} "5";"6"};
{\ar@{=>}^{\; \scriptstyle{\alpha}} "7";"8"};
\end{xy}
\end{equation}
\end{defn}

The following is our first main result, which solves the problem raised in~\cite[Section 6]{LackS:cohapm}.

\begin{theorem}
\label{thm:psm-gray}
Let $\catK$ be a Gray-category. Then there is a Gray-category $\Psm(\catK)$, called the Gray-category of pseudomonads in $\catK$,
having pseudomonads in $\catK$ as 0-cells, pseudomonad morphisms as 1-cells, pseudomonad
transformations as 2-cells, and pseudomonad modifications as 3-cells.
\end{theorem}

The rest of this section is devoted to the proof of Theorem~\ref{thm:psm-gray}, which will be obtained
by combining Lemmas~\ref{lemma:strict-hom-cats}, ~\ref{lemma:comp-ass}, ~\ref{lemma:cub-comp} and~\ref{lemma:ass-2-3-cells} below.  We begin by giving the
definition of the hom-2-categories of $\Psm(\catK)$.

\begin{lemma}
\label{lemma:strict-hom-cats}
Let $(X,\,S)$ and $(Y,\,T)$ be two pseudomonads in $\catK$. Then there is a 2-category 
$\Psm(\catK)(\,(X,\,S),\,(Y,\,T)\,)$ having 
pseudomonad morphisms from~$(X,S)$ to~$(Y,T)$ as 0-cells, pseudomonad transformations as 1-cells and and pseudomonad modifications as 2-cells. 
\end{lemma}

\begin{proof}
First of all, for any pair of composable 1-cells
$(p_0,\,\widetilde{p_0})\colon (F_0,\phi_0)\rightarrow(F_1,\phi_1)$ and $(p_1,\,\widetilde{p_1})\colon (F_1,\,\phi_1)\rightarrow(F_2,\,\phi_2)$
we define their composition as $(p_1p_0,\,\widetilde{p_1p_0})$ where~$\widetilde{p_1p_0}$ is defined as the pasting of 
\begin{center}
\begin{tikzpicture}
\node (a) at (0,2) {$TF_0$};
\node (b) at (2,2) {$TF_1$};
\node (c) at (0,0) {$F_0S$};
\node (d) at (2,0) {$F_1S$};
\node (f) at (4,0) {$F_2S$.};
\node (e) at (4,2) {$TF_2$};
\draw[->] (a) to node[scale=.7] (Tp) [above]{$Tp_0$} (b);
\draw[->] (a) to node[scale=.7] (phisx) [left]{$\phi_0$} (c);
\draw[->] (b) to node[scale=.7] (phidx) [right]{$\phi_1$} (d);
\draw[->] (c) to node[scale=.7] (pS) [below]{$p_0S$} (d);

\draw[->] (b) to node[scale=.7] (Tp1) [above]{$Tp_1$} (e);
\draw[->] (e) to node[scale=.7] (phidx) [right]{$\phi_2$} (f);
\draw[->] (d) to node[scale=.7] (p1S) [below]{$p_1S$} (f);
\draw[-{Implies},double distance=1.5pt,shorten >=20pt,shorten <=20pt] (Tp) to node[scale=.7] [right] [xshift=0.2cm] {$\widetilde{p_0}$} (pS);
\draw[-{Implies},double distance=1.5pt,shorten >=20pt,shorten <=20pt] (Tp1) to node[scale=.7] [right] [xshift=0.2cm] {$\widetilde{p_1}$} (p1S);
\end{tikzpicture}
\end{center}
We want to show that composition is strictly associative. So let us consider three composable 1-cells
$$(F_0,\,\phi_0)\xrightarrow{(p_0,\,\widetilde{p_0})}(F_1,\,\phi_1)\xrightarrow{(p_1,\,\widetilde{p_1})}(F_2,\,\phi_2)\xrightarrow{(p_2,\,\widetilde{p_2})}(F_3,\,\phi_3)$$
By definition, the two possible composites are
\begin{align*} 
(p_2,\,\widetilde{p_2})\cdot \big( (p_1,\,\widetilde{p_1})\cdot(p_0,\,\widetilde{p_0}) \big) & =(p_2(p_1p_0),\,\widetilde{p_2(p_1p_0)}) \, , \\ 
\big ( (p_2,\,\widetilde{p_2})\cdot(p_1,\,\widetilde{p_1}) \big)\cdot(p_0,\,\widetilde{p_0}) & =((p_2p_1)p_0,\,\widetilde{(p_2p_1)p_0}) \, .
\end{align*}
We want to show that these are equal. Since $\catK$ is a Gray-category, $p_2(p_1p_0)=(p_2p_1)p_0$. Moreover,  $\widetilde{p_2(p_1p_0)}=\widetilde{(p_2p_1)p_0}$ since they are both the pasting of 
\begin{center}
\begin{tikzpicture}
\node (a) at (0,2) {$TF_0$};
\node (b) at (2,2) {$TF_1$};
\node (c) at (0,0) {$F_0S$};
\node (d) at (2,0) {$F_1S$};
\node (f) at (4,0) {$F_2S$};
\node (e) at (4,2) {$TF_2$};
\node (g) at (6,2) {$TF_3$};
\node (h) at (6,0) {$F_3S$.};
\draw[->] (a) to node[scale=.7] (Tp) [above]{$Tp_0$} (b);
\draw[->] (a) to node[scale=.7] (phisx) [left]{$\phi_0$} (c);
\draw[->] (b) to node[scale=.7] (phidx) [right]{$\phi_1$} (d);
\draw[->] (c) to node[scale=.7] (pS) [below]{$p_0S$} (d);

\draw[->] (b) to node[scale=.7] (Tp1) [above]{$Tp_1$} (e);
\draw[->] (e) to node[scale=.7] (phidx) [right]{$\phi_2$} (f);
\draw[->] (d) to node[scale=.7] (p1S) [below]{$p_1S$} (f);

\draw[->] (e) to node[scale=.7] (Tp2) [above]{$Tp_2$} (g);
\draw[->] (g) to node[scale=.7] (phidx) [right]{$\phi_3$} (h);
\draw[->] (f) to node[scale=.7] (p2S) [below]{$p_2S$} (h);
\draw[-{Implies},double distance=1.5pt,shorten >=20pt,shorten <=20pt] (Tp) to node[scale=.7] [right] [xshift=0.2cm] {$\widetilde{p_0}$} (pS);
\draw[-{Implies},double distance=1.5pt,shorten >=20pt,shorten <=20pt] (Tp1) to node[scale=.7] [right] [xshift=0.2cm] {$\widetilde{p_1}$} (p1S);
\draw[-{Implies},double distance=1.5pt,shorten >=20pt,shorten <=20pt] (Tp2) to node[scale=.7] [right] [xshift=0.2cm] {$\widetilde{p_2}$} (p2S);
\end{tikzpicture}
\end{center}

It remains to define the identity 1-cells of $\Psm(\catK)(\,(X,\,S),\,(Y,\,T)\,)$. For a pseudomonad morphism $(F, \phi) \colon (X, S) \to (Y, T)$,
we define the identity on it to be 
\[
(1_F,\,1_\phi):(F,\,\phi)\rightarrow(F,\,\phi) \, .
\] 
This is allowed since $T1_F=1_{TF}$ and $1_FS=1_{FS}$. These can be shown to be a strict identities, using  that $\catK$ is a Gray-category
and in particular Axiom (G\ref{ax-gray-id}).
\end{proof} 

We proceed by defining the composition of $1$-cells in $\Psm(\catK)$ and proving that is {\em strictly} associative, as required to have a Gray-category.
Since a pseudomonad morphism is a tuple of the form $(F,  \phi, \bar{\phi}, \tilde{\phi})$, where
$F$ is a 1-cell, $\phi$ is a 2-cell while $\bar{\phi}$ and $\tilde{\phi}$ are 3-cells, we will need to check
equalities at three levels. The key level of the verification is that of 2-cells. Indeed, 
strict associativity at the level of 1-cells will follow easily from the strict associativity of
 composition of 1-cells in $\catK$. The key issue are the equalities at the level of 2-cells, since 2-cells could be isomorphic (by means of an invertible 3-cell), 
but not equal.  Instead, equalities of 3-cells will be quite straightforward. In fact,  the required equations for 3-cells either hold strictly or they fail completely, 
since there are no 4-cells that could make these equations hold only up to isomorphism.

 In the following, for a pseudomonad morphism $\underline{F}=(F,\,\phi,\,\bar{\phi},\,\widetilde{\phi})$, we define
\[
\underline{F}\;\bar{}:=\bar{\phi} \, , \quad\quad\quad\underline{F}\;\widetilde{}:=\widetilde{\phi} \, .
\]

Let $(F,\,\phi):(X,\,S)\rightarrow(Y,\,T)$ and $(G,\,\psi):(Y,\,T)\rightarrow(Z,\,Q)$ be two pseudomonad morphisms. We define their composition as 
\begin{equation}
\label{equ:comp-in-psm}
(G,\,\psi, \underline{G}\;\bar{}, \underline{G}\;\widetilde{} \; )
\circ
(F,\,\phi,  \underline{F}\;\bar{}, \underline{F}\;\widetilde{} \; )
:= \big(GF,\,G\phi\cdot\psi F \, , \underline{G}\circ\underline{F}\;\bar{} \, , \underline{G}\circ\underline{F}\;\widetilde{} \; \big) 
\end{equation}
where the invertible $3$-cells are defined by the following pasting diagrams:
\begin{center}
\begin{tikzpicture}
\node (a) at (0,6) {$Q^2GF$};
\node (b) at (3,6) {$QGTF$};
\node (d) at (3,4) {$GT^2F$};
\node (e) at (0,0) {$QGF$};
\node (f) at (3,0) {$GTF$};
\node (QGFS) at (6,6) {$QGFS$};
\node (GTFS) at (6,4) {$GTFS$};
\node (GFS2) at (6,2) {$GFS^2$};
\node (GFS) at (6,0) {$GFS$};
\node at (-3,3) {$\underline{G}\circ\underline{F}\;\bar{}:=$};
\draw[->] (a) to node[scale=.7] (gF) [above]{$Q\psi F$} (b);
\draw[->] (a) to node[scale=.7] (Gf) [left]{$m_QGF$} (e);
\draw[->] (b) to node[scale=.7] (G'f) [left] {$\psi TF$} (d);
\draw[->] (e) to node[scale=.7] (gF'') [below] {$\psi F$} (f);
\draw[->] (d) to node[scale=.7] (Gf2) [left]{$G\mut F$} (f);
\draw[->] (b) to node[scale=.7] (QGphi) [above]{$QG\phi$} (QGFS);
\draw[->] (QGFS) to node[scale=.7] (psiFS) [right]{$\psi FS$} (GTFS);

\draw[->] (d) to node[scale=.7] (gF2) [below]{$GT\phi$} (GTFS);
\draw[->] (GTFS) to node[scale=.7] (G'f2) [right] {$G\phi S$} (GFS2);
\draw[->] (GFS2) to node[scale=.7] (G'f2) [right] {$GF\mus $} (GFS);
\draw[->] (f) to node[scale=.7] (gF''2) [below] {$G\phi$} (GFS);
\draw[-{Implies},double distance=1.5pt,shorten >=75pt,shorten <=75pt] (gF) to node[scale=.7] [left] {$\bar{\psi}F$} (gF'');
\draw[-{Implies},double distance=1.5pt,shorten >=40pt,shorten <=40pt] (gF2) to node[scale=.7] [left] {$G\bar{\phi}$} (gF''2);
\draw[-{Implies},double distance=1.5pt,shorten >=20pt,shorten <=20pt] (QGphi) to node[scale=.7] [left] {${\psi_\phi}$} (gF2);
\end{tikzpicture}\\
\end{center}

\begin{center}
\begin{tikzpicture}
\node (GF) at (0,4) {$GF$};
\node (GFS) at (6,4) {$GFS$};
\node (QGF) at (2,0) {$QGF$};
\node (GTF) at (4,2) {$GTF$};
\node (a) at (4,4) {};
\node at (-3,2) {$\underline{G}\circ\underline{F}\;\widetilde{}:=$};
\draw[->] (GF) to node[scale=.7] (GFs) [above]{$GFs$} (GFS);
\draw[->] (GF) to node[scale=.7] (GtF) [above, xshift=0.2cm, yshift=0.1cm]{$GtF$} (GTF);
\draw[->] (GTF) to node[scale=.7] (Gphi) [right, xshift=0.2cm, yshift=-0.1cm] {$G\phi$} (GFS);
\draw[->] (GF) to node[scale=.7] (qGF) [left] {$qGF$} (QGF);
\draw[->] (QGF) to node[scale=.7] (psiF) [right, xshift=0.2cm, yshift=-0.1cm]{$\psi F$} (GTF);
\draw[-{Implies},double distance=1.5pt,shorten >=10pt,shorten <=10pt] (a) to node[scale=.7] [left, xshift=-0.2cm] {$G\widetilde{\phi}$} (GTF);
\draw[-{Implies},double distance=1.5pt,shorten >=25pt,shorten <=25pt] (GtF) to node[scale=.7] [left, xshift=-0.1cm] {$\widetilde{\psi}F$} (QGF);
\end{tikzpicture}
\end{center}

The proof that this definition gives a pseudomonad morphism is in Appendix~\ref{app:coherence-comp}.

\begin{rmk}
We did not consider any parenthesis in the diagrams above thanks to axiom~(G\ref{equ:gray}) for a Gray-category. Moreover since $Q(-)$ is a strict $2$-functor we have $Q(G\phi\cdot\psi F)=QG\phi\cdot Q\psi F$ (and similarly for other compositions in the diagrams).
\end{rmk}

\begin{lemma}
\label{lemma:comp-ass}
The composition of pseudomonad morphisms defined in~\eqref{equ:comp-in-psm} is strictly associative. 
\end{lemma}

\begin{proof}
From now on, let us consider three pseudomonad morphisms in $\catK$: 
$$(X,\,S)\xrightarrow{(F,\,\phi)}(Y,\,T)\xrightarrow{(G,\,\psi)}(Z,\,Q)\xrightarrow{(H,\,\xi)}(V,\,R)$$
In order to prove this statement we have to prove that the equation for associativity holds for the respective $1$-, $2$- and $3$-cell components. For 1-cells, since $\catK$ is a Gray-category, then $H(GF)=(HG)F$. 

For 2-cells, the idea is to reduce both composites to $HG\phi\cdot H\psi F\cdot\xi GF$. 
%
%
%
%
%
%
On the one hand, 
\begin{align*}
H(G\phi\cdot\psi F)\cdot\xi GF 
&  = \left[ H(G\phi)\cdot H(\psi F)\right] \cdot \xi GF 
& (\textrm{because}\;H(-)\;\textrm{is strict}) 
\\
& = \left[ HG\phi\cdot H\psi F\right] \cdot \xi GF 
& (\textrm{by}\;(\textrm{G}\ref{equ:gray})) 
\\
& =  HG\phi\cdot H\psi F \cdot \xi GF 
& (\textrm{since}\;\catK(X,\,V)\;\textrm{is a}\;2\textrm{-category}).
\end{align*}
On the other hand, 
\begin{align*}
HG\phi\cdot (H\psi\cdot\xi G)F & = HG\phi\cdot\left[ (H\psi)F\cdot(\xi G)F\right] & (\textrm{because}\;(-)F\;\textrm{is strict}) \\
& =  HG\phi\cdot \left[H\psi F \cdot \xi GF\right]
& (\textrm{by}\;(\textrm{G}\ref{equ:gray})) 
\\
& =  HG\phi\cdot H\psi F \cdot \xi GF 
& (\textrm{since}\;\catK(X,\,V)\;\textrm{is a}\;2\textrm{-category}). 
\end{align*}
For 3-cells, to prove that $\big( \underline{H}(\underline{G}\underline{F}) \big) \; \widetilde{}= \big( ( \underline{H}\underline{G})\underline{F} \big)\;\widetilde{}$ we just need to notice that, using the fact that $H(-)$ and $(-)F$ are strict $2$-functors, both of them are pasting of: 
\begin{center}
\begin{tikzpicture}
\node (HGF) at (-0.5,4.5) {$HGF$};
\node (HGFS) at (6.5,4.5) {$HGFS$};
\node (a) at (5,4.5) {};
\node (HGTF) at (5,3) {$HGTF$};
\node (b) at (3.7,3.3) {};
\node (e) at (3.7,1.8) {};
\node (HQGF) at (3.5,1.5) {$HQGF$};
\node (c) at (2.3,2) {};
\node (d) at (2.3,0.5) {};
\node (RHGF) at (2,0) {$RHGF$.};
\draw[->] (HGF) to node[scale=.7] (HGFs) [above]{$HGFs$} (HGFS);
\draw[->] (HGF) to node[scale=.7] (HGFs) [left]{$HqGF$} (HQGF);
\draw[->] (HGF) to node[scale=.7] (GtF) [above]{$HGtF$} (HGTF);
\draw[->] (HGTF) to node[scale=.7] (HGphi) [right] {$HG\phi$} (HGFS);
\draw[->] (HGF) to node[scale=.7] (rHGF) [left] {$rHGF$} (RHGF);
\draw[->] (RHGF) to node[scale=.7] (xiGF) [right]{$\xi GF$} (HQGF);
\draw[->] (HQGF) to node[scale=.7] (HpsiF) [right]{$H\psi F$} (HGTF);
\draw[-{Implies},double distance=1.5pt,shorten >=7pt,shorten <=7pt] (a) to node[scale=.7] [left] {$HG\widetilde{\phi}$} (HGTF);
\draw[-{Implies},double distance=1.5pt,shorten >=7pt,shorten <=7pt] (b) to node[scale=.7] [left] {$H\widetilde{\psi}F$} (e);
\draw[-{Implies},double distance=1.5pt,shorten >=7pt,shorten <=7pt] (c) to node[scale=.7] [left] {$\widetilde{\xi}GF$} (d);
\end{tikzpicture}
\end{center}
We get the required equality by the pasting theorem for $2$-categories\cite{AJPow90}. Finally, let us prove the equality on the other $3$-cell component. By definition,

\[
\begin{tikzpicture}
\node (a) at (0,4) {$R^2HGF$};
\node (b) at (3,4) {$RHQGF$};
\node (c) at (6,4) {$RHGFS$};
\node (d) at (3,2) {$HQ^2GF$};
\node (e) at (6,2) {$HQGFS$};
\node (f) at (0,0) {$RHGF$};
\node (g) at (3,0) {$HQGF$};
\node (h) at (6,0) {$HGFS$\,.};
\node at (-1.5,2) {$\underline{H}(\underline{G}\underline{F})\;\bar{}\,=$};
\draw[->] (a) to node[scale=.7](F) {} (b);
\draw[->] (b) to node[scale=.7](G) {} (c);
\draw[->] (d) to node[scale=.7](H) {} (e);
\draw[->] (f) to node[scale=.7](I) {} (g);
\draw[->] (g) to node[scale=.7](J) {} (h);
\draw[->] (b) to (d);
\draw[->] (c) to (e); 
\draw[->] (d) to (g);
\draw[->] (e) to (h); 
\draw[->] (a) to (f);
\2-cells 
\draw[-{Implies},double distance=1.5pt,shorten >=45pt,shorten <=45pt] (F) to node[scale=.7] [left] {$\bar{\xi}GF$} (I);
\draw[-{Implies},double distance=1.5pt,shorten >=20pt,shorten <=20pt] (G) to node[scale=.7] [left] {$\xi_{(G\phi\cdot\psi F)}$} (H);
\draw[-{Implies},double distance=1.5pt,shorten >=20pt,shorten <=20pt] (H) to node[scale=.7] [left] {$H(\underline{G}\underline{F}\;\bar{}\,)$} (J);
\end{tikzpicture} 
\]
Using the definition of $\underline{G}\underline{F}\;\bar{}$ and (\ref{equ:cc3}), the right-hand side pasting becomes
\[
\begin{tikzpicture}
\node (A) at (-2.5,6) {$R^2HGF$}; 
\node (B) at (0,6) {$RHQGF$};
\node (C) at (2.5,6) {$RHGTF$};
\node (D) at (5,6) {$RHGFS$};
\node (E) at (-2.5,0) {$RHGF$}; 
\node (a) at (0,4) {$HQ^2GF$};
\node (b) at (2.5,4) {$HGTFS$};
\node (c) at (5,4) {$HQGFS$};
\node (d) at (2.5,2) {$HGT^2F$};
\node (e) at (5,2) {$HGTFS$};
\node (f) at (0,0) {$HQGF$};
\node (g) at (2.5,0) {$HGTF$};
\node (h) at (5,0) {$HGFS$\,.};
\draw[->] (a) to node[scale=.7](F) {} (b);
\draw[->] (b) to node[scale=.7](G) {} (c);
\draw[->] (d) to node[scale=.7](H) {} (e);
\draw[->] (f) to node[scale=.7](I) {} (g);
\draw[->] (g) to node[scale=.7](J) {} (h);
\draw[->] (B) to node[scale=.7](F1) {} (C);
\draw[->] (C) to node[scale=.7](G1) {} (D);
\draw[->] (A) to node[scale=.7](S1) {} (B);
\draw[->] (E) to node[scale=.7](S) {} (f);
\draw[->] (A) to (E); 
\draw[->] (B) to (a); 
\draw[->] (C) to (b); 
\draw[->] (D) to (c); 
\draw[->] (b) to (d);
\draw[->] (c) to (e); 
\draw[->] (d) to (g);
\draw[->] (e) to (h); 
\draw[->] (a) to (f);
\draw[-{Implies},double distance=1.5pt,shorten >=45pt,shorten <=45pt] (F) to node[scale=.7] [left] {$H\bar{\psi}F$} (I);
\draw[-{Implies},double distance=1.5pt,shorten >=20pt,shorten <=20pt] (G) to node[scale=.7] [right] {$H(\psi_\phi)$} (H);
\draw[-{Implies},double distance=1.5pt,shorten >=20pt,shorten <=20pt] (H) to node[scale=.7] [right] {$HG\bar{\phi}$} (J);
\draw[-{Implies},double distance=1.5pt,shorten >=70pt,shorten <=70pt] (S1) to node[scale=.7] [right] {$\bar{\xi}GF$} (S);
\draw[-{Implies},double distance=1.5pt,shorten >=20pt,shorten <=20pt] (F1) to node[scale=.7] [right] {$\xi_{\psi F}$} (F);
\draw[-{Implies},double distance=1.5pt,shorten >=20pt,shorten <=20pt] (G1) to node[scale=.7] [right] {$\xi_{G\phi}$} (G);
\end{tikzpicture} 
\]

Let us notice that, by (G\ref{equ:gray}), $H(\psi_\phi)={H\psi}_\phi$, $\xi_{\psi F}=(\xi_\psi)F$ and $\xi_{G\phi}=\xi G_\phi$. Moreover, using the definition of $\underline{H}\underline{G}\;\bar{}$ and (\ref{equ:cc2}), the diagram above is equal to
\begin{center}
\begin{tikzpicture}
\node (a) at (0,4) {$R^2HGF$};
\node (b) at (3,4) {$RHGTF$};
\node (c) at (6,4) {$RHGFS$};
\node (d) at (3,2) {$HGT^2FS$};
\node (e) at (6,2) {$HGTFS$};
\node (f) at (0,0) {$RHGF$};
\node (g) at (3,0) {$HGTF$};
\node (h) at (6,0) {$HGFS$,};
\node (G') at (4, 4) {};
\node (H') at (4, 2) {};
\node (A) at (1.5,4) {};
\node (B) at (1.5,0) {};

\draw[->] (a) to node[scale=.7](F) {} (b);
\draw[->] (b) to node[scale=.7](G) {} (c);
\draw[->] (d) to node[scale=.7](H) {} (e);
\draw[->] (f) to node[scale=.7](I) {} (g);
\draw[->] (g) to node[scale=.7](J) {} (h);
\draw[->] (b) to (d);
\draw[->] (c) to (e); 
\draw[->] (d) to (g);
\draw[->] (e) to (h); 
\draw[->] (a) to (f);

\2-cells 
\draw[-{Implies},double distance=1.5pt,shorten >=45pt,shorten <=45pt] (A) to node[scale=.7] [left] {$(\underline{H}\underline{G}\;\bar{}\,)F$} (B);
\draw[-{Implies},double distance=1.5pt,shorten >=18pt,shorten <=18pt] (G') to node[scale=.7] [right] {$(H\psi\cdot\xi G)_\phi$} (H');
\draw[-{Implies},double distance=1.5pt,shorten >=18pt,shorten <=18pt] (H) to node[scale=.7] [right] {$HG\bar{\phi}$} (J);
\end{tikzpicture}
\end{center} 
which is exactly the definition of $\underline{H}(\underline{G}\underline{F})\;\bar{}$.
\end{proof}

For brevity, we sometimes write $P_\catK$ instead of $\Psm(\catK)$, so for any pair of pseudomonads $(X,\,S)$ and $(Y,\,T)$ the $2$-category of pseudomonads morphisms from $(X,\,S)$ to $(Y,\,T)$ can be written
as $P_\catK\big((X,\,S),\,(Y,\,T)\big)$.

\begin{lemma}
\label{lemma:cub-comp}
The definition of composition of pseudomonad morphisms extends to a cubical functor 
\[
-\circ -:P_\catK(\,(Y,\,T),\,(Z,\,Q)\,)\times P_\catK(\,(X,\,S),\,(Y,\,T)\,)\longrightarrow P_\catK(\,(X,\,S),\,(Z,\,Q)\,) 
\]
for $(X, S) \, , (Y, T) \, , (Z, Q) \in \Psm(\cat{K})$. 
\end{lemma}

\begin{proof} This is a just a long verification, but we spell  it out in some detail. 
By the definition of a cubical functor,  for $(F,\,\phi):(X,\,S)\rightarrow(Y,\,T)$ 
and~$(G,\,\psi):(Y,\,T)\rightarrow(Z,\,Q)$ in~$\Psm(\cat{K})$, we need to define strict $2$-functors 
\begin{align}
F_\phi \colon & P_\catK(\,(Y,\,T),\,(Z,\,Q)\,)\rightarrow P_\catK(\,(X,\,S),\,(Z,\,Q)\,) \, , \label{equ:F-phi} \\
G_\psi \colon & P_\catK(\,(X,\,S),\,(Y,\,T)\,)\rightarrow P_\catK(\,(X,\,S),\,(Z,\,Q)\,) \label{equ:G-psi}
\end{align}
such that 
\begin{equation}
\label{equ:F-phi-on-obj}
F_\phi(\,(G,\,\psi)\,)=G_\psi(\,(F,\,\phi)\,)= (G,\,\psi)\circ (F,\,\phi) \, , 
\end{equation}
plus, for 2-cells $(p,\,\widetilde{p}):(F,\,\phi)\rightarrow(F',\,\phi')$ and $(q,\,\widetilde{q}):(G,\,\psi)\rightarrow(G',\,\psi')$, an invertible 3-cell in $\Psm(\catK)$
\begin{equation}
\label{equ:Big-Sigma}
\begin{gathered} 
\begin{tikzpicture}
\node (a) at (0,2) {$(G,\,\psi)\circ (F,\,\phi)$};
\node (c) at (5.5,2) {$(G,\,\psi)\circ (F',\,\phi')$};
\node (b) at (0,0) {$(G',\,\psi')\circ (F,\,\phi)$};
\node (d) at (5.5,0) {$(G',\,\psi')\circ (F',\,\phi')$};

\draw[->] (a) to node[scale=.7] (Tp) [left]{$(q,\,\widetilde{q})\circ (F,\,\phi)$} (b);
\draw[->] (a) to node[scale=.7] (phisx) [above]{$(G,\,\psi)\circ (p,\,\widetilde{p})$} (c);
\draw[->] (b) to node[scale=.7] (phidx) [below]{$(G',\,\psi')\circ (p,\,\widetilde{p})$} (d);
\draw[->] (c) to node[scale=.7] (pS) [right]{$(q,\,\widetilde{q})\circ (F',\,\phi')$} (d);
\draw[-{Implies},double distance=1.5pt,shorten >=20pt,shorten <=20pt] (phisx) to node[scale=.7] [right] [xshift=0.2cm] {$\Sigma_{(p,\,\widetilde{p}),\,(q,\,\widetilde{q})}$} (phidx);
\end{tikzpicture}
\end{gathered}
\end{equation}
satisfying  axioms (\ref{equ:cc1}), (\ref{equ:cc2}) and (\ref{equ:cc3}).

We begin by defining $F_\phi$ in~\eqref{equ:F-phi}. Its action on objects is determined by~\eqref{equ:F-phi-on-obj}. For
its action on 1-cells, we send $(q,\,\widetilde{q}) \colon (G, \psi) \to (G', \psi)$ to the pseudomonad modification~$(qF, \widetilde{qF}) \colon (GF, G\phi \cdot \psi F) \to (G'F, G' \phi \cdot \psi' F)$, 
where $\widetilde{qF}$ is defined as the following pasting:
\begin{center}
\begin{tikzpicture}
\label{cubical comp 2}
\node (a) at (0,4) {$QGF$};
\node (b) at (2,4) {$QG'F$};
\node (c) at (0,2) {$GTF$};
\node (d) at (2,2) {$G'TF$};
\node (e) at (0,0) {$GFS$};
\node (f) at (2,0) {$G'FS$.};
\draw[->] (a) to node[scale=.7] (gF) [above]{$QqF$} (b);
\draw[->] (a) to node[scale=.7] (Gf) [left]{$\psi F$} (c);
\draw[->] (c) to node[scale=.7] (gF') [below] {$qTF$} (d);
\draw[->] (b) to node[scale=.7] (G'f) [right] {$\psi'F$} (d);
\draw[->] (c) to node[scale=.7] (Gf'') [left] {$G\phi$} (e);
\draw[->] (d) to node[scale=.7] (G'f'') [right] {$G'\phi$} (f);
\draw[->] (e) to node[scale=.7] (gF'') [below] {$qFS$} (f);
\draw[-{Implies},double distance=1.5pt,shorten >=20pt,shorten <=20pt] (gF) to node[scale=.7] [right, xshift=0.2cm] {$\widetilde{q}F$} (gF');
\draw[-{Implies},double distance=1.5pt,shorten >=20pt,shorten <=10pt] (gF') to node[scale=.7] [right, xshift=0.2cm, yshift=0.2cm] {${q_\phi}^{-1}$} (gF'');
\end{tikzpicture}
\end{center}
The action of $F_\phi$ on 3-cells $\beta:(q,\,\widetilde{q})\rightarrow(q',\,\widetilde{q'})$ is defined 
by letting $\beta\circ(F,\,\phi):=\beta F$ in $\mathbb{\catK}$. The proof that this is a pseudomonad modification, and therefore a 3-cells in $\Psm(\catK)$, is  in Appendix~\ref{app:well-defineness-Fphi}.

We now  show that $F_\phi$ is a 2-functor. For this, we use extensively the notation of Remark~\ref{gray-notation} to avoid writing some diagrams. To prove that composition is preserved strictly, we show that
\begin{equation}
\label{equ:F-phi-strict-comp}
F_\phi(\,(q',\,\widetilde{q'})\cdot(q,\,\widetilde{q})\,)=F_\phi(q',\,\widetilde{q'})\cdot F_\phi(q,\,\widetilde{q})
\end{equation}
for any  
\[
(G,\,\psi)\xrightarrow{(q,\,\widetilde{q})}(G',\,\psi')\xrightarrow{(q',\,\widetilde{q'})}(G',\,\psi')
\]
in $P_\catK(\,(Y,\,T),\,(Z,\,Q)\,)$. The composition 
$(q',\,\widetilde{q'})\cdot(q,\,\widetilde{q})$
is defined as $(q'q,\,\widetilde{q'q})$
where $\widetilde{q'q}$ is defined as the pasting of 
\begin{center}
\begin{tikzpicture}
\node (a) at (0,1.5) {$QG$};
\node (b) at (1.5,1.5) {$QG'$};
\node (c) at (0,0) {$GT$};
\node (d) at (1.5,0) {$G'T$};
\node (f) at (3,0) {$G'T$.};
\node (e) at (3,1.5) {$QG'$};
\draw[->] (a) to node[scale=.7] (Tp) [above]{} (b);
\draw[->] (a) to node[scale=.7] (phisx) [left]{} (c);
\draw[->] (b) to node[scale=.7] (phidx) [right]{} (d);
\draw[->] (c) to node[scale=.7] (pS) [below]{} (d);

\draw[->] (b) to node[scale=.7] (Tp1) [above]{} (e);
\draw[->] (e) to node[scale=.7] (phidx) [right]{} (f);
\draw[->] (d) to node[scale=.7] (p1S) [below]{} (f);
\draw[-{Implies},double distance=1.5pt,shorten >=15pt,shorten <=15pt] (Tp) to node[scale=.7] [right] [xshift=0.2cm] {$\widetilde{q}$} (pS);
\draw[-{Implies},double distance=1.5pt,shorten >=15pt,shorten <=15pt] (Tp1) to node[scale=.7] [right] [xshift=0.2cm] {$\widetilde{q'}$} (p1S);
\end{tikzpicture}
\end{center}
Using the equation (\ref{equ:cc3}) we can see that the 3-cells components of $F_\phi(\,(q',\widetilde{q'})\cdot(q,\widetilde{q})\,)$ and~$F_\phi(q',\,\widetilde{q'})\cdot$~$F_\phi(q,\,\widetilde{q})$ are two pasting of the diagram below
\begin{center}
\begin{tikzpicture}
\node (a) at (0,3) {$QG$};
\node (b) at (1.5,3) {$QG'$};
\node (c) at (0,1.5) {$GT$};
\node (d) at (1.5,1.5) {$G'T$};
\node (f) at (3,1.5) {$G'T$};
\node (e) at (3,3) {$QG'$};
\node (g) at (0,0) {$GF$};
\node (h) at (1.5,0) {$G'F$};
\node (i) at (3,0) {$G'F$.};

\draw[->] (a) to node (Tp) [above]{} (b);
\draw[->] (a) to node (phisx) [left]{} (c);
\draw[->] (b) to node (phidx) [right]{} (d);
\draw[->] (c) to node (pS) [below]{} (d);

\draw[->] (b) to node (Tp1) [above]{} (e);
\draw[->] (e) to node (phidx) [right]{} (f);
\draw[->] (d) to node (p1S) [below]{} (f);

\draw[->] (c) to node [below]{} (g);
\draw[->] (d) to node [below]{} (h);
\draw[->] (f) to node [below]{} (i);

\draw[->] (g) to node (pS') [below]{} (h);
\draw[->] (h) to node (p1S') [below]{} (i);

\draw[-{Implies},double distance=1.5pt,shorten >=15pt,shorten <=15pt] (Tp) to node[scale=.7] [right] [xshift=0.2cm] {$\widetilde{q}$} (pS);
\draw[-{Implies},double distance=1.5pt,shorten >=15pt,shorten <=15pt] (Tp1) to node[scale=.7] [right] [xshift=0.2cm] {$\widetilde{q'}$} (p1S);

\draw[-{Implies},double distance=1.5pt,shorten >=10pt,shorten <=10pt] (pS) to node[scale=.7] [right] [xshift=0.1cm] {${q_\phi}^{-1}$} (pS');
\draw[-{Implies},double distance=1.5pt,shorten >=10pt,shorten <=10pt] (p1S) to node[scale=.7] [right] [xshift=0.1cm] {${q'_\phi}^{-1}$} (p1S');
\end{tikzpicture}
\end{center}
Moreover, $(q'\cdot q)F=q'F\cdot qF$ since $(-)F$ is a strict 2-functor (since $\catK$ is a Gray category). Hence, the required equality in~\eqref{equ:F-phi-strict-comp} holds. Let us also verify that $F_\phi$  
preserves identities strictly. Recall from Lemma~\ref{lemma:strict-hom-cats} that $1_{(G,\,\psi)}:=(1_G,\,1_\psi)$ in $P_\catK(\,(Y,\,T),\,(Z,\,Q)\,)$.
Therefore,  
\[
F_\phi(1_G,\,1_\psi) = (1_GF,\,\widetilde{1_GF})
\]
and moreover
\begin{align*}
(1_GF,\,\widetilde{1_GF})&=(1_{GF},\,(\,(1_G)_\phi\cdot1_{\psi F})\ast(1_{G\phi}\cdot1_\psi F)\,) & (\textrm{by definition of}\,F_\phi)\\
&=(1_{GF},\,(1_{G\phi}\cdot1_{\psi F})\ast(1_{G\phi}\cdot1_{\psi F})\,) & (\textrm{by Remark \ref{rmk:cub-comp}})\\
&=(1_{GF},\,(1_{G\phi\cdot\psi F})\ast(1_{G\phi\cdot\psi F})\,) & (\textrm{since }\cdot\textrm{ preserves identities})\\
&=(1_{GF},\,1_{G\phi\cdot\psi F}) & (\textrm{by (G\ref{ax:strict-hom})})\\
&=1_{(GF,\,G\phi\cdot\psi F)} & \\
& =1_{F_\phi(G,\,\psi)}\, ,
\end{align*}
as required.

We now define the 2-functor $G_\phi$ of~\eqref{equ:G-psi}. Again, its action on objects is determined by~\eqref{equ:F-phi-on-obj}.  On
morphisms, it sends $(p, \widetilde{p}) \colon (F, \phi) \to (F', \phi')$ to 
\[
(Gp, \widetilde{Gp}) \colon (GF, G \phi \cdot \psi F) \to (GF', G \phi' \cdot \psi F') \, ,
\]
where $\widetilde{Gp}$ is defined as the following pasting:
\begin{center}
\begin{tikzpicture}
\label{cubical comp 2}
\node (a) at (0,4) {$QGF$};
\node (b) at (2,4) {$QGF'$};
\node (c) at (0,2) {$GTF$};
\node (d) at (2,2) {$GTF'$};
\node (e) at (0,0) {$GFS$};
\node (f) at (2,0) {$GF'S \, .$};
\draw[->] (a) to node[scale=.7] (gF) [above]{$QGp$} (b);
\draw[->] (a) to node[scale=.7] (Gf) [left]{$\psi F$} (c);
\draw[->] (c) to node[scale=.7] (gF') [below] {$GTp$} (d);
\draw[->] (b) to node[scale=.7] (G'f) [right] {$\psi F'$} (d);
\draw[->] (c) to node[scale=.7] (Gf'') [left] {$G\phi$} (e);
\draw[->] (d) to node[scale=.7] (G'f'') [right] {$G\phi'$} (f);
\draw[->] (e) to node[scale=.7] (gF'') [below] {$GpS$} (f);
\draw[-{Implies},double distance=1.5pt,shorten >=20pt,shorten <=20pt] (gF) to node[scale=.7] [right, xshift=0.2cm] {$\psi_p$} (gF');
\draw[-{Implies},double distance=1.5pt,shorten >=20pt,shorten <=10pt] (gF') to node[scale=.7] [right, xshift=0.2cm, yshift=0.2cm] {$G\widetilde{p}$} (gF'');
\end{tikzpicture}
\end{center}
On 3-cells $\alpha:(p,\,\widetilde{p})\rightarrow(p',\,\widetilde{p'})$ we let $(G,\,\phi)\circ\alpha :=G\alpha$, which is a  3-cell in $\Psm(\catK)$  by a similar argument to the one used for $F_\phi$. The proof that this is a 2-functor is completely analogous to the one for $F_\phi$ and hence omitted.


To conclude the proof, we need to define the 3-cell $\Sigma_{(p,\,\widetilde{p}),\,(q,\,\widetilde{q})}$  in~\eqref{equ:Big-Sigma}. We take 
this to be $q_p$, which is shown to be a pseudomonad modification in Appendix~\ref{app:coherence-interch}. 
The required axioms for $\Sigma_{(p,\,\widetilde{p}),\,(q,\,\widetilde{q})}$, as in (\ref{equ:cc1}), (\ref{equ:cc2}) and (\ref{equ:cc3}), hold as they are instances of the ones for $q_p$ for $\catK$. 
\end{proof}

\begin{lemma}
\label{lemma:ass-2-3-cells}
The cubical functor providing composition in $\Psm(\cat{K})$ satisfies the coherence
conditions of Axiom (G\ref{equ:gray}).
 \end{lemma}

\begin{proof}
The first one is just Lemma~\ref{lemma:comp-ass}. Since the definitions on 3-cells coincide with the ones in $\catK$, all the equations regarding them hold directly. Therefore, we only need to prove the ones for 2-cells. Let us consider the following diagram in $\Psm(\catK)$:
\begin{center}
\begin{tikzpicture}
\node (a) at (0,0) {$(X,\,S)$};
\node (b) at (4,0) {$(Y,\,T)$};
\node (c) at (8,0) {$(Z,\,Q)$};
\node (d) at (12,0) {$(V,\,R)$}; 

\draw[->] (a) to [bend left] node[scale=.7] (f) [above] {$(F,\,\phi)$} (b);
\draw[->] (a) to [bend right] node[scale=.7] (f') [below] {$(F',\,\phi')$} (b);
\draw[-{Implies},double distance=1.5pt,shorten >=10pt,shorten <=10pt] (f) to node[scale=.7] [right, xshift=0.2cm] {$(p,\,\widetilde{p})$} (f');

\draw[->] (b) to [bend left] node[scale=.7] (g) [above] {$(G,\,\psi)$} (c);
\draw[->] (b) to [bend right] node[scale=.7] (g') [below] {$(G',\,\psi')$} (c);
\draw[-{Implies},double distance=1.5pt,shorten >=10pt,shorten <=10pt] (g) to node[scale=.7] [right, xshift=0.2cm] {$(q,\,\widetilde{q})$} (g');

\draw[->] (c) to [bend left] node[scale=.7] (h) [above] {$(H,\,\xi)$} (d);
\draw[->] (c) to [bend right] node[scale=.7] (h') [below] {$(H',\,\xi')$} (d);
\draw[-{Implies},double distance=1.5pt,shorten >=10pt,shorten <=10pt] (h) to node[scale=.7] [right, xshift=0.2cm] {$(r,\,\widetilde{r})$} (h');
\end{tikzpicture}
\end{center}
We need to prove: 
\begin{enumerate}[(i)]
\item $(H_\xi\circ G_\psi)\,(p,\,\widetilde{p})=H_\xi(G_\psi\,(p,\,\widetilde{p})\,)$, 
\item $(\,H_\xi\,(q,\,\widetilde{q})\,)F_\phi=H_\xi(\,(q,\,\widetilde{q})F_\phi\,)$, 
\item $(\,(r,\,\widetilde{r})\,G_\psi\,)F_\phi=(r,\,\widetilde{r})\,(G_\psi\circ F_\phi\,)$. 
\end{enumerate}
At the 2-cells level we have $(HG)p=H(Gp)$ because $\catK$ is a Gray-category. The same happens in (ii) and (iii) so we will just prove that the associated 3-cells are equal in each case. 
Let us start with part (i). On the one hand, by definition,  
\[
\widetilde{(HG)p}=\widetilde{H\widetilde{Gp}} \, , 
\] 
and therefore
\begin{center}
\begin{tikzpicture}[scale=0.9]
\node (A) at (-2.5,5) {$RHGF$}; 
\node (B) at (0,5) {$RHGF'$};
\node (C) at (2.5,6) {$RHGF$};
\node (D) at (5,6) {$RHGF'$};
\node (E) at (-2.5,1) {$HGFS$};
\node (F') at (-2.5,3) {$HQGF$}; 
\node (a) at (0,3) {$HQGF'$};
\node (b) at (2.5,4) {$HQGF$};
\node (c) at (5,4) {$HQGF'$};
\node (d) at (2.5,2) {$HGTF$};
\node (e) at (5,2) {$HGTF'$};
\node (f) at (0,1) {$HGF'S$};
\node (g) at (2.5,0) {$HGFS$};
\node (h) at (5,0) {$HGF'S$.};
\node at (1.25,3) {$=$};

\draw[->] (b) to node[scale=.7](G) {} (c);
\draw[->] (d) to node[scale=.7](H) {} (e);
\draw[->] (g) to node[scale=.7](J) {} (h);
\draw[->] (C) to node[scale=.7](G1) {} (D);
\draw[->] (A) to node[scale=.7](S1) {} (B);
\draw[->] (E) to node[scale=.7](S) {} (f);
\draw[->] (F') to node[scale=.7](S') {} (a);
\draw[->] (A) to (F');
\draw[->] (F') to (E); 
\draw[->] (B) to (a); 
\draw[->] (C) to (b); 
\draw[->] (D) to (c); 
\draw[->] (b) to (d);
\draw[->] (c) to (e); 
\draw[->] (d) to (g);
\draw[->] (e) to (h); 
\draw[->] (a) to (f);

\draw[-{Implies},double distance=1.5pt,shorten >=15pt,shorten <=15pt] (G) to node[scale=.7] [right] {$H\psi_p$} (H);
\draw[-{Implies},double distance=1.5pt,shorten >=15pt,shorten <=15pt] (H) to node[scale=.7] [right] {$HG\widetilde{p}$} (J);
\draw[-{Implies},double distance=1.5pt,shorten >=15pt,shorten <=15pt] (S1) to node[scale=.7] [right] {$\xi_{Gp}$} (S');
\draw[-{Implies},double distance=1.5pt,shorten >=15pt,shorten <=15pt] (S') to node[scale=.7] [right] {$H\widetilde{Gp}$} (S);
\draw[-{Implies},double distance=1.5pt,shorten >=15pt,shorten <=15pt] (G1) to node[scale=.7] [right] {$\xi_{Gp}$} (G);
\end{tikzpicture}
\end{center}
On the other hand $H_\xi(G_\psi\,(p,\,\widetilde{p})\,)=(HG,\,H\psi\cdot\xi G)\,(p,\,\widetilde{p})$ so the associated 3-cell is, using  (\ref{equ:cc2}), 
\begin{center}
\begin{tikzpicture}[scale=0.9]
\node (A) at (-2.5,5) {$RHGF$}; 
\node (B) at (0,5) {$RHGF'$};
\node (C) at (2.5,6) {$RHGF$};
\node (D) at (5,6) {$RHGF'$};
\node (E) at (-2.5,1) {$HGFS$};
\node (F') at (-2.5,3) {$HGTF$}; 
\node (a) at (0,3) {$HGTF'$};
\node (b) at (2.5,4) {$HQGF$};
\node (c) at (5,4) {$HQGF'$};
\node (d) at (2.5,2) {$HGTF$};
\node (e) at (5,2) {$HGTF'$};
\node (f) at (0,1) {$HGF'S$};
\node (g) at (2.5,0) {$HGFS$};
\node (h) at (5,0) {$HGF'S$.};
\node at (1.25,3) {$=$};

\draw[->] (b) to node[scale=.7](G) {} (c);
\draw[->] (d) to node[scale=.7](H) {} (e);
\draw[->] (g) to node[scale=.7](J) {} (h);
\draw[->] (C) to node[scale=.7](G1) {} (D);
\draw[->] (A) to node[scale=.7](S1) {} (B);
\draw[->] (E) to node[scale=.7](S) {} (f);
\draw[->] (F') to node[scale=.7](S') {} (a);
\draw[->] (A) to (F');
\draw[->] (F') to (E); 
\draw[->] (B) to (a); 
\draw[->] (C) to (b); 
\draw[->] (D) to (c); 
\draw[->] (b) to (d);
\draw[->] (c) to (e); 
\draw[->] (d) to (g);
\draw[->] (e) to (h); 
\draw[->] (a) to (f);

\node (j) at (-2,5) {};
\node (k) at (-2,3) {};
\draw[-{Implies},double distance=1.5pt,shorten >=15pt,shorten <=15pt] (j) to node[scale=.7] [right] {$(H\psi\cdot\xi G)_p$} (k);
\draw[-{Implies},double distance=1.5pt,shorten >=15pt,shorten <=15pt] (S') to node[scale=.7] [right] {$HG\widetilde{p}$} (S);

\draw[-{Implies},double distance=1.5pt,shorten >=15pt,shorten <=15pt] (G1) to node[scale=.7] [right] {$\xi G_{p}$} (G);
\draw[-{Implies},double distance=1.5pt,shorten >=15pt,shorten <=15pt] (G) to node[scale=.7] [right] {$H\psi_p$} (H);
\draw[-{Implies},double distance=1.5pt,shorten >=15pt,shorten <=15pt] (H) to node[scale=.7] [right] {$HG\widetilde{p}$} (J);
\end{tikzpicture} 
\end{center}
But $\xi G_p=\xi_{Gp}$, since $\catK$ is a Gray-category, and so the required equality holds.

For part (ii), by definition, the 3-cell component of $H_\xi\,(q,\,\widetilde{q})\,)F_\phi$ is: 
\begin{center}
\begin{tikzpicture}[scale=0.9]
\node (A) at (-2.5,5) {$RHGF$}; 
\node (B) at (0,5) {$RHG'F$};
\node (F') at (-2.5,3) {$HGTF$}; 
\node (a) at (0,3) {$HG'TF$};
\node (E) at (-2.5,1) {$HGFS$};
\node (f) at (0,1) {$HG'FS$};

\node at (1.25,3) {$=$};

\node (C) at (2.5,6) {$RHGF$};
\node (D) at (5,6) {$RHGF'$};
\node (b) at (2.5,4) {$HQGF$};
\node (c) at (5,4) {$HQG'F$};
\node (d) at (2.5,2) {$HGTF$};
\node (e) at (5,2) {$HG'TF$};
\node (g) at (2.5,0) {$HGFS$};
\node (h) at (5,0) {$HG'FS$};

\node at (6.25,3) {$=$};

\node (A') at (7.5,5) {$RHGF$}; 
\node (B') at (10,5) {$RHGF'$};
\node (F'') at (7.5,3) {$HQGF$}; 
\node (a') at (10,3) {$HQG'F$};
\node (E') at (7.5,1) {$HGFS$};
\node (f') at (10,1) {$HG'FS$.};

\draw[->] (A) to node(S1) {} (B);
\draw[->] (A) to (F');
\draw[->] (F') to (E);  
\draw[->] (B) to (a); 
\draw[->] (F') to node(S') {} (a);
\draw[->] (E) to node(S) {} (f);
\draw[->] (a) to (f);
\draw[->] (F') to (E);

\draw[->] (b) to node(G) {} (c);
\draw[->] (d) to node(H) {} (e);
\draw[->] (g) to node(J) {} (h);
\draw[->] (C) to node(G1) {} (D);
\draw[->] (C) to (b); 
\draw[->] (D) to (c); 
\draw[->] (b) to (d);
\draw[->] (c) to (e); 
\draw[->] (d) to (g);
\draw[->] (e) to (h); 

\draw[->] (A') to node(S1') {} (B');
\draw[->] (A') to (F'');
\draw[->] (F'') to (E'); 
\draw[->] (B') to (a'); 
\draw[->] (F'') to node(S'') {} (a');
\draw[->] (E') to node(S2) {} (f');
\draw[->] (a') to (f');
\draw[->] (F'') to (E');

\draw[-{Implies},double distance=1.5pt,shorten >=15pt,shorten <=15pt] (S1) to node[scale=.7] [right] {$\widetilde{Hq}F$} (S');
\draw[-{Implies},double distance=1.5pt,shorten >=15pt,shorten <=15pt] (S') to node[scale=.7] [right] {${Hq_\phi}^{-1}$} (S);

\draw[-{Implies},double distance=1.5pt,shorten >=15pt,shorten <=15pt] (G1) to node[scale=.7] [right] {$(\xi_q)F$} (G);
\draw[-{Implies},double distance=1.5pt,shorten >=15pt,shorten <=15pt] (G) to node[scale=.7] [right] {$H\widetilde{q}F$} (H);
\draw[-{Implies},double distance=1.5pt,shorten >=15pt,shorten <=15pt] (H) to node[scale=.7] [right] {${Hq_\phi}^{-1}$} (J);

\draw[-{Implies},double distance=1.5pt,shorten >=15pt,shorten <=15pt] (S1') to node[scale=.7] [right] {$\xi_{qF}$} (S'');
\draw[-{Implies},double distance=1.5pt,shorten >=15pt,shorten <=15pt] (S'') to node[scale=.7] [right] {$H(\widetilde{qF})$} (S2);
\end{tikzpicture}
\end{center}
Finally, part (iii) is completely analogous to the first one using the inverse of \ref{equ:cc3} instead of \ref{equ:cc2}. 
\end{proof}

The combination of Lemmas~\ref{lemma:strict-hom-cats}, ~\ref{lemma:comp-ass}, ~\ref{lemma:cub-comp} and  ~\ref{lemma:ass-2-3-cells} proves Theorem~\ref{thm:psm-gray}. 

\section{Liftings to pseudoalgebras}
\label{sec:lift}

We now recall from~\cite[Section~7]{MarmolejoF:dislp} and~\cite[Section~6]{LackS:cohapm}
the definition of the Gray-category $\Lift(\catK)$ of pseudomonads in~$\catK$ and liftings to pseudoalgebras. In \cite{MarmolejoF:dislp} this was written as $\Psm(\catK)$, but we prefer to use that notation for the Gray-category introduced in Section~\ref{sec:psem}, since it seems the natural generalization of the 2-category of monads defined by Street in \cite{StreetR:fortm}. We will then show that $\Lift(\catK)$ is equivalent to $\Psm(\catK)$, which will be used in Section~\ref{sec:psel} for
our results on pseudodistributive laws.

The 0-cells of $\Lift(\catK)$ are  pseudomonads $(X,S)$ in $\catK$.  For 0-cells $(X,S)$ and $(Y,T)$,
a 1-cell $(F, \lift{F}) : (X,S) \rightarrow (Y, T)$ consists of  a 1-cell $F : X \rightarrow Y$ 
in $\catK$ and a Gray-transformation
 $\lift{F} : \psSalg \rightarrow \psTalg$ making the following diagram commute
\begin{equation*}
\label{equ:lift}
\xymatrix{
\psSalg \ar[r]^{\lift{F}} \ar[d]_{U} & \psTalg \ar[d]^{U} \\
X  \ar[r]_{F}                     &  Y }
\end{equation*} 
where, using implicitly the Yoneda lemma for Gray-categories, we write $X$ and $Y$ instead of $\catK(X,\,-)$ and $\catK(Y,\,-)$. We refer to $\lift{F}$ as a 
{\em lifting} of $F$ to pseudoalgebras. Analogous terminology will be used
for the 2- and 3-cells introduced below. 

\begin{lemma}
\label{lift-morph-lemma}
Let $(F, \phi) : (X,S) \rightarrow (Y, T)$  be a pseudomonad morphism. Then, there exists a lifting $\lift{F} : \psSalg \rightarrow \psTalg$ of $F : X \rightarrow Y$. 
\end{lemma}

\begin{proof}
Let us consider a fixed $I \in \catK$. First, let us observe that if 
$A$ is an $I$-indexed  pseudo-$S$-algebra, then $FA$  is naturally an $I$-indexed  
pseudo-$\mt$-algebra,  with structure map  given by the composite
\[
\xymatrix{
\mth A \ar[r]^{\phi_A} & \mhs A \ar[r]^{Fa} & FA }
\]
and associativity and unit 3-cells provided by the pasting diagrams 
\[
\begin{xy}
(-10,30)*+{\mtth A}="4";
(-10,0)*+{\mth A}="3";
(10,30)*+{\mths A}="2";
(10,0)*+{\mhs A}="1"; 
(10,15)*+{\mhss A}="0"; 
(-1,18)="5";
(-1,14)="6";
(24,25)="11";
(24,20)="12";
(24,9)="13";
(24,4)="14";
(40,30)*+{\mth A}="8";
(40,15)*+{\mhs A}="9";
(40,0)*+{F A \, ,}="10";
{\ar^{\mt \phi_A} "4";"2"};
{\ar_{\phi_A} "3";"1"};
{\ar_{\mut_{(FA)}} "4";"3"};
{\ar^{\phi_{\ms A}} "2";"0"};
{\ar^{F \mus_A} "0";"1"};
{\ar^{\mth a} "2";"8"};
{\ar_{\mhs a} "0";"9"};
{\ar_{F a} "1";"10"};
{\ar^{\phi_A} "8";"9"};
{\ar^{Fa} "9";"10"};
{\ar@{=>}^{\; \scriptstyle{\phi_a}} "11";"12"};
{\ar@{=>}^{\; \scriptstyle{F \bar{a}}} "13";"14"};
{\ar@{=>}^{\; \scriptstyle{\bar{\phi}_A}} "5";"6"};
\end{xy}
\hspace{1.5cm} 
\begin{xy}
(-15,30)*+{FA}="4";
(8,30)*+{\mth A}="2";
(8,0)*+{F A \, .}="1";
(2,23)="5";
(-3,23)="6";
(2,9)="8";
(-3,9)="9";
(8,15)*+{\mhs A}="7";
{\ar@/_1pc/_{F \es_A} "4";"7"};
{\ar@/_2pc/_{1_{FA}} "4";"1"};
{\ar^(.45){\es_{FA}} "4";"2"};
{\ar^{\phi_A} "2";"7"};
{\ar^{Fa} "7";"1"};
{\ar@{=>}^(.4){{\scriptstyle{\tilde{\phi}_A} \; }} "6";"5"};
{\ar@{=>}^(.4){{\scriptstyle{F \tilde{a}} \; }} "9";"8"};
\end{xy} 
\]
The coherence condition~\eqref{equ:alg1} for  $FA$ follows by an application of the coherence condition~\eqref{equ:m1} for $F$ and the coherence condition~\eqref{equ:alg1} for $A$. The  coherence 
condition~\eqref{equ:alg2} for $FA$ follows by applying the 
coherence condition~\eqref{equ:m2} for  $F$ and the coherence 
condition~\eqref{equ:alg2} for $A$.  Secondly, we observe that if $f : A \rightarrow B$ is
a pseudo-$\ms$-algebra morphism, then $Ff : FA \rightarrow FB$ is naturally a  
pseudo-$\mt$-algebra morphism, as we have the following pasting diagram:
\[
\begin{xy}
(10,30)*+{\mth A}="2";
(10,0)*+{FA}="1"; 
(10,15)*+{\mhs A}="0"; 
(24,25)="11";
(24,20)="12";
(24,9)="13";
(24,4)="14";
(40,30)*+{\mth B}="8";
(40,15)*+{\mhs B}="9";
(40,0)*+{F B \, .}="10";
{\ar_{\phi_{A}} "2";"0"};
{\ar_{F a} "0";"1"};
{\ar^{\mth f} "2";"8"};
{\ar_{\mhs f} "0";"9"};
{\ar_{Ff} "1";"10"};
{\ar^{\phi_B} "8";"9"};
{\ar^{Fb} "9";"10"};
{\ar@{=>}^{\; \scriptstyle{\phi_f}} "11";"12"};
{\ar@{=>}^{\; \scriptstyle{F \bar{f} }} "13";"14"};
\end{xy}
\] 
The coherence conditions~\eqref{equ:mor1} and~\eqref{equ:mor2} follow immediately by
the axioms for a Gray-category. 
Finally, if $\alpha : f \rightarrow g$ is a pseudo-$\ms$-algebra
$2$-cell, the required pseudo-$\mt$-algebra 2-cell is given by 
$F \alpha : F f \rightarrow F g$.  We have thus defined the components of
a Gray-natural
transformation $\lift{F} : \psSalg \rightarrow \psTalg$, which is clearly a lifting of $F : X \rightarrow Y$. 
\end{proof}

Given 1-cells $(F, \lift{F}) : (X,S) \rightarrow (Y, T)$ and $(F', \lift{F'}) : (X,S) \rightarrow (Y, T)$, a 2-cell~$(p, \lift{p}) : (F, \lift{F}) \rightarrow (F', \lift{F'})$ in $\Lift(\catK)$ consists of a 2-cell $p : F \rightarrow F'$ in $\catK$
and  a Gray-modification $\lift{p} : \lift{F} \rightarrow \lift{F'}$ such that
the following diagram commutes
\begin{equation*}
 \label{equ:liftp}
\xymatrix{
U  \lift{F}  \ar[rr]^{U \lift{p}} \ar@{=}[d]
& & U \lift{F'}  \ar@{=}[d] \\
F  U  \ar[rr]_{p U} & &  F' U .} 
\end{equation*}
The vertical arrows are the identities, which hold by the assumption that $\lift{F}$ and $\lift{F'}$
are liftings of $F$ and $F'$, respectively. 

\begin{lemma}
\label{lift-transf-lemma}
Let $(p, \tilde{p}) : (F, \phi) \rightarrow (F', \phi')$  be a pseudomonad transformation. 
Then, there exists a lifting $\lift{p} : \lift{F} \rightarrow \lift{F'}$ of
 $p : F \rightarrow F'$, where~$\lift{F}$ and~$\lift{F'}$ are the liftings of~$F$ and~$F'$ 
 associated to the pseudomonad morphisms~$(F, \phi)$ and~$(F', \phi')$, respectively, defined as in Lemma~\ref{lift-morph-lemma}.
\end{lemma}

\begin{proof}
Let $I \in \catK$.  We need to define a pseudonatural transformation $\lift{p} :
\lift{F}_I \rightarrow \lift{F'}_I$. We
define the component of $\lift{p}$ associated to an~$I$-indexed pseudo-$S$-algebra $A$ to be 
the~$I$-indexed pseudo-$T$-algebra morphism given by $p_A : FA \rightarrow F'A$ and the 2-cell  
\[
\begin{xy}
(10,30)*+{\mth A}="2";
(10,0)*+{FA}="1"; 
(10,15)*+{\mhs A}="0"; 
(24,25)="11";
(24,20)="12";
(24,9)="13";
(24,4)="14";
(40,30)*+{\mtk A}="8";
 (40,15)*+{\mks A}="9";
(40,0)*+{F' A \, .}="10";
{\ar_{\phi_{A}} "2";"0"};
{\ar_{F a} "0";"1"};
{\ar^{\mt p_A} "2";"8"};
{\ar_{p_{\ms A}} "0";"9"};
{\ar_{p_A} "1";"10"};
{\ar^{\phi'_A} "8";"9"};
{\ar^{F'a} "9";"10"};
{\ar@{=>}^{\; \scriptstyle{\bar{p}_A}} "11";"12"};
{\ar@{=>}^{\; \scriptstyle{p_{a}^{-1}}} "13";"14"};
\end{xy}
\]
To prove the  condition~\eqref{equ:mor1} for the pseudoalgebra morphism $p_A$, we 
apply the axioms for a Gray-category and then condition~\eqref{equ:t1} for the pseudomonad transformation $p$. To establish
condition~\eqref{equ:mor2}, it is sufficient to apply the coherence condition~\eqref{equ:t2} 
for the pseudomonad transformation $p$, and then the axioms for a Gray-category.  
By definition,~$\lift{p}$ is a lifting of $p$ as required. 
\end{proof}

Finally, for 2-cells $(p, \lift{p})$ and $(q, \lift{q})$, a 
3-cell $\alpha : (p, \lift{p}) \rightarrow (q, \lift{q})$ consists of a 3-cell and $\alpha : p \rightarrow q$ and a Gray-perturbation  $\lift{\alpha} : \lift{p} \rightarrow \lift{q}$ making the following diagram commute 
\begin{equation*}
\label{equ:liftalpha}
\xymatrix{
U \lift{p} \ar[rr]^{U \lift{\alpha}} \ar@{=}[d]
& & U \lift{q} \ar@{=}[d] \\
p U \ar[rr]_{\alpha U} & & qU \, .} 
\end{equation*}
As before, the vertical arrows are the identities that are part of the assumption that~$\lift{p}$ 
and~$\lift{q}$ are liftings of $p$ and $q$, respectively. Composition and identities of
$\Lift(\catK)$ are defined in the evident way, using those of $\catK$ and $\Gray$.

\begin{lemma}
\label{lift-mod-lemma}
Given a pseudomonad modification $\alpha : (p, \tilde{p})  \rightarrow (q, \tilde{q})$
we can define a lifting $\lift{\alpha} : \lift{p} \rightarrow \lift{q}$ of $\alpha$ as the 
Gray-perturbation whose components are, for a pseudo-$\ms$-algebra $A$, the
3-cells~$\alpha_A : p_A \rightarrow q_A$.
\end{lemma}

\begin{proof} 
It suffices to check that,  these $3$-cells are a pseudo-$\mt$-algebra 
$2$-cells. To prove this, apply the axioms for a Gray-category and the coherence 
axiom~\eqref{equ:md}. 
\end{proof}

\medskip

We use these results to define a Gray-functor 
\[
\Phi:\Psm(\catK)\rightarrow\Lift(\catK) \, .
\]
On objects, $\Phi$ acts as the identity. For two pseudomonads $(X,S)$ and $(Y,T)$ in $\catK$, the hom-2-functors
\[
\Phi_{(X,S),\,(Y,T)}:\Psm(\catK)(\,(X,S),\,(Y,T)\,)\longrightarrow\Lift(\catK)(\,(X,S),\,(Y,T)\,) 
\]
are defined sending a pseudomonad morphism, pseudomonad transformation and pseudomonad modification to the
associated liftings, using Lemmas~\ref{lift-morph-lemma}, \ref{lift-transf-lemma} and \ref{lift-mod-lemma}, respectively,
Here, the Gray-functoriality of~$\Phi$ is standard verification, which we omit for brevity, limiting ourselves to highlight that this includes checking that~$\Phi$ preserves composition strictly, {\em i.e.}
that the lifting associated to the composite of two pseudomonad morphisms is equal (rather than
just equivalent by invertible 2-cells) to the composite of the liftings obtained from the pseudomonad morphisms. Theorem~\ref{thm:equiv-lift-psm} states that the 
construction of $\Psm(\cat{K})$ given in Section~\ref{sec:psem} is equivalent to the one by Marmolejo in~\cite{MarmolejoF:dislp}.

\begin{theorem}
\label{thm:equiv-lift-psm}
The Gray-functor $\Phi \colon \Psm(\catK) \to \Lift(\catK)$ is a triequivalence. 
\end{theorem}

\begin{proof}
Since $\Phi$ is clearly bijective on objects, it suffices to prove that locally it is a biequivalence. Let us begin by considering a lifting $\lift{F} : \psSalg \rightarrow \psTalg$
of a 1-cell~$F : X \rightarrow Y$.  By the definition of a lifting, the following diagram of 
2-categories and 2-functors commutes:
\begin{equation}
\label{lift-diagr}
\begin{gathered}
\xymatrix{
\psSalg(X) \ar[rr]^{\lift{F}_X} \ar[d]_{U_X} & &  \psTalg(X) \ar[d]^{U_X} \\
\catK(X,X)  \ar[rr]_{\catK(X,F)}                     &  &  \catK(X,Y) \, .}
\end{gathered}
\end{equation}
Let us now observe that $S : X \rightarrow X$ can be regarded as an $X$-indexed pseudo-$S$-algebra,
with structure map given by the 2-cell $\mus : \mss \rightarrow \ms$. By the commutativity of the
diagram~(\ref{lift-diagr}), this pseudo-$S$-algebra is mapped by the 2-functor~$\lift{F}_X$ into a pseudo-$T$-algebra with underlying 1-cell  $\mhs : X \rightarrow Y$, with structure map 
a 2-cell of the form~$\phi_0 : \mths \rightarrow \mhs$, and invertible 3-cells fitting in the diagrams
\[
\begin{xy}
(-12,8)*+{\mtths }="4";
(-12,-8)*+{\mths  }="3";
(12,8)*+{\mths }="2";
(12,-8)*+{\mhs }="1"; 
(-1,2)="5";
(-1,-3)="6";
{\ar^{\mt \phi_0} "4";"2"};
{\ar_{h'} "3";"1"};
{\ar_{\mut_{\mhs }} "4";"3"};
{\ar^{\phi_0} "2";"1"};
{\ar@{=>}^{\ \scriptstyle{\bar{\phi_0}}} "5";"6"};
\end{xy}
\hspace{1.5cm}
\begin{xy}
(-16,8)*+{\mhs}="4";
(8,8)*+{\mths}="2";
(8,-8)*+{\mhs\,.}="1";
(3,1)="5";
(-3,1)="6";
{\ar@/_/_{1_{\mhs}} "4";"1"};
{\ar^(.45){\et_{\mhs }} "4";"2"};
{\ar^{h'} "2";"1"};
{\ar@{=>}^(.4){{\scriptstyle{\tilde{\phi_0}}}} "6";"5"};
\end{xy}
\]
The desired pseudomonad morphism $(F, \phi) : (X,S) \rightarrow (Y,T)$ is then obtained by
letting~$\phi : \mth \rightarrow \mhs$ be the composite 
\[
\xymatrixcolsep{36pt}
\xymatrix{ 
\mth  \ar[r]^(.42){\mth s} & \mths  \ar[r]^{\phi_0} \ar[r] & \mhs .}
\]
The appropriate 3-cells are provided by the following pasting diagrams
\[
\begin{xy}
(10,0)*+{\mths}="1"; 
(10,30)*+{\mtths}="2";
(-15,0)*+{\mth}="3";
(-15,30)*+{\mtth}="4";
(-5,18)="5";
(-5,14)="6";
(30,30)*+{\mths}="8";
(55,0)*+{\mhs}="10";
(24,18)="11";
(24,14)="12";
(24,9)="13";
(24,4)="14";
(55,30)*+{\mthss}="15";
(55,15)*+{\mhss}="16";
(44,23)*+{}="17";
(44,18)*+{}="18";
{\ar^{\mtth \es} "4";"2"};
{\ar_{\mth \es} "3";"1"};
{\ar_{\mut F} "4";"3"};
{\ar^{\mut \mhs} "2";"1"};
{\ar^{\mt \phi_0} "2";"8"};
{\ar_{\phi_0} "1";"10"};
{\ar@/_1pc/_(.45){\phi_0} "8";"10"};
{\ar^{\mth \es \ms} "8";"15"};
{\ar^{\phi_0 \ms} "15";"16"};
{\ar^{\mh \mus} "16";"10"};
{\ar@{=>}^{\; \bar{\phi_0}} "11";"12"};
{\ar@{=>}^{\; \mut_{\mh \es}} "5";"6"};
{\ar@{=>}_{\; \gamma} "17";"18"};
\end{xy}
\hspace{.5cm} 
\begin{xy}
(-15,30)*+{\mh}="0";
(-15,15)*+{\mhs}="4";
(8,30)*+{\mth}="2";
(8,0)*+{\mhs\,,}="1";
(0,23)="5";
(-5,23)="6";
(2,7)="8";
(-3,7)="9";
(8,15)*+{\mths}="7";
{\ar_{F \es} "0";"4"};
{\ar^{\et \mh } "0";"2"};
{\ar@/_1.5pc/_{1_{\mhs}} "4";"1"};
{\ar^(.45){\et \mhs} "4";"7"};
{\ar^{\mth \es} "2";"7"};
{\ar^{\phi_0} "7";"1"};
{\ar@{=>}^(.4){{\scriptstyle{\et_{\mh \es}} \; }} "6";"5"};
{\ar@{=>}^(.4){{\scriptstyle{\tilde{\phi_0}} \; }} "9";"8"};
\end{xy} 
\]
where $\gamma$  is the inverse to the  2-cell obtained from the 
following pasting of invertible 2-cells:
\begin{equation*}
\begin{xy}
(15,20)*+{\mths}="1";
(15,0)*+{\mthss}="5";
(40,0)*+{\mths}="6";
(15,-20)*+{\mhss}="9"; 
(40,-20)*+{\mhs\,.}="10";
(24,-7)*+{}="15";
(24,-12)*+{}="21";
{\ar@/^1pc/^{1_{\mths}} "1";"6"};
{\ar^{\phi_0} "6";"10"};
{\ar_{\phi_0 \ms} "5";"9"};
{\ar_{\mth \es \ms} "1";"5"};
{\ar^{\mth \mus } "5";"6"};
{\ar_{\mh \mus } "9";"10"};
{\ar@{=>}^{\; F \alpha \; } "15";"21"};
(24,12)*+{}="30";
(24,7)*+{}="31";
{\ar@{=>}^<<<{\mth \rho} "30";"31"};
(9,-2)*+{}="36";
(9,-7)*+{}="37";
\end{xy}
\end{equation*}
Let us now consider a lifting $(p, \lift{p}) : (F, \lift{F}) 
\rightarrow (F', \lift{F'})$ of a 2-cell $p : F \rightarrow F'$. 
We can define a pseudomonad transformation 
$p: (F, \phi) \rightarrow (F', \phi')$  by considering the following pasting diagram:
\[
\begin{xy}
(10,30)*+{\mth}="2";
(10,0)*+{\mhs}="1"; 
(10,15)*+{\mths}="0"; 
(24,25)="11";
(24,20)="12";
(24,9)="13";
(24,4)="14";
(40,30)*+{\mtk}="8";
 (40,15)*+{\mtks}="9";
(40,0)*+{\mks\,,}="10";
{\ar_{\mth \es} "2";"0"};
{\ar_{\phi_o} "0";"1"};
{\ar^{\mt p} "2";"8"};
{\ar_{\mt p \ms} "0";"9"};
{\ar_{p \ms} "1";"10"};
{\ar^{\mtk \es} "8";"9"};
{\ar^{\phi'_0} "9";"10"};
{\ar@{=>}^{Tp_u^{-1}}  "11";"12"};
{\ar@{=>}^{\; \scriptstyle{\bar{p \ms}}} "13";"14"};
\end{xy}
\]
in which the bottom 3-cell is part of the structure making $p {\ms} : \mhs  \rightarrow \mks$ into 
a pseudoalgebra morphism. Finally, if $(\alpha, \lift{\alpha}) : (p, \lift{p}) \rightarrow
(q, \lift{q})$ is a lifting of a 3-cell~$\alpha : p \rightarrow  q$, 
then $\alpha : p \rightarrow q$ is a pseudomonad modification.
These definitions determine a 2-functor
\[
\Psi_{(X,S), (Y,T)} :  \Lift(\catK)  \big( (\bx, \ms), (\by, \mt) \big) 
 \longrightarrow  \Psm(\catK) \big( (\bx, \ms), (\by, \mt) \big) 
\]
which provides the required quasi-inverse to $\Phi_{(X,S), (Y, T)}$. We omit the construction of
the required invertible pseudonatural transformations $\eta : 1 \rightarrow \Psi\Phi$
and $\varepsilon : \Phi\Psi \rightarrow 1$, since this is not difficult.
\end{proof}

\section{Pseudodistributive laws}
\label{sec:psel}

\begin{defn}
Let $(X, S)$ and $(X,T)$ be pseudomonads in $\catK$. A \myemph{pseudodistributive law} of $T$ over $S$ consists of a 2-cell $\di : \mst \rightarrow~\mts$ and invertible 3-cells
\[
\begin{xy}
(-10,30)*+{\msst}="4";
(-10,0)*+{\mst }="3";
(10,30)*+{\msts}="2";
(10,0)*+{\mts}="1"; 
(10,15)*+{\mtss}="0"; 
(-1,18)="5";
(-1,14)="6";
{\ar^{\ms \di} "4";"2"};
{\ar_{\di} "3";"1"};
{\ar_{\mus \mt} "4";"3"};
{\ar^{\di \ms} "2";"0"};
{\ar^{\mt \mus} "0";"1"};
{\ar@{=>}^{\; \scriptstyle{\bar{\mus}}} "5";"6"};
\end{xy}
\hspace{2cm} 
\begin{xy}
(-8,30)*+{\mt}="4";
(8,30)*+{\mst}="2";
(8,0)*+{\mts}="1";
(4,20)="5";
(-1,20)="6";
{\ar@/_/_{\mt \es} "4";"1"};
{\ar^(.45){\es \mt} "4";"2"};
{\ar^{\di} "2";"1"};
{\ar@{=>}^(.4){{\scriptstyle{\bar{\es}}}} "6";"5"};
\end{xy} 
\]
\[
\begin{xy}
(-10,30)*+{\mstt}="4";
(-10,0)*+{\mtts}="3";
(10,30)*+{\mst}="2";
(10,0)*+{\mts}="1"; 
(-10,15)*+{\mtst}="0"; 
(-1,18)="5";
(-1,14)="6";
{\ar^{\ms \mut} "4";"2"};
{\ar_{\mut \ms} "3";"1"};
{\ar_{\di \mt} "4";"0"};
{\ar^{\di} "2";"1"};
{\ar_{\mt \di} "0";"3"};
{\ar@{=>}^{\; \scriptstyle{\bar{\mut}}} "5";"6"};
\end{xy}
\hspace{2cm}
\begin{xy}
(-8,0)*+{\ms}="4";
(8,30)*+{\mst}="2";
(8,0)*+{\mts}="1";
(2,9)="5";
(2,14)="6";
{\ar_{\et \ms} "4";"1"};
{\ar@/^/^(.45){\ms \et} "4";"2"};
{\ar^{\di} "2";"1"};
{\ar@{=>}^(.4){\; {\scriptstyle{\bar{\et}}}} "6";"5"};
\end{xy} 
\]
satisfying the coherence conditions (C1)-(C8) stated in Appendix~\ref{app:coh}. 
\end{defn}

\begin{rmk} For the convenience
of the reader, Table~\ref{tab:comp-coh-axioms} describes the correspondence between 
the presentation of the coherence conditions for pseudodistributive laws here and
in \cite{TanakaM:psedl, MarmolejoF:dislp}. In the table, each row lists 
different formulations of the same  axiom.

\begin{table}[hbt]
\begin{tabular}{|c|c|c|} 
\hline
Appendix \ref{app:coh} & Marmolejo \cite{MarmolejoF:dislp} & Tanaka \cite{TanakaM:psedl}\\ \hline
(C1) & (coh 4) & (T6) \\
(C2) & (coh 2) & (T2) \\ 
(C3) & (coh 5) & (T9) \\
(C4) & (coh 3) & (T8) \\
(C5) & (coh 1) & (T1) \\
(C6) & (coh 6) & (T10) \\
(C7) & (coh 9) & (T7) \\
(C8) & (coh 7) & (T5) \\
(C9) & - & (T3) \\
(C10) & (coh 8) & (T4) \\
\hline
\end{tabular}
\caption{Comparison of coherence conditions.}
\label{tab:comp-coh-axioms} 
\end{table} 
\end{rmk}

\begin{rmk} \label{thm:explain-coherence}
Our development in Section~\ref{sec:lift} allows us to give a clear explanation for the coherence conditions for pseudodistributive laws, summarised in Table~\ref{tab:coherence-axioms}. 

\begin{table}[htb]
\begin{tabular}{|c|c|} 
\hline
Axiom & Coherence condition \\ \hline
(C1) and (C2) & $(T, d) : (X,S) \rightarrow (X,S)$ is a pseudomonad  morphism \\ 
(C3) and (C4) & $(\mut, \bar{\mut}) : (\mt, \di)^2  \rightarrow (\mt, \di)$ is a
pseudomonad transformation \\
(C5) and~(C6) &  $(\et, \bar{\et}) :  (X, 1_X) \rightarrow (\mt, \di)$ is a pseudomonad  transformation \\
(C7) & $\alpha$ is pseudomonad modification \\
(C8) & $\rho$ is a pseudomonad modification \\
(C9) & $\lambda$ is a pseudomonad modification \\
\hline
\end{tabular}
\caption{Coherence axioms for pseudodistributive laws.}
\label{tab:coherence-axioms} 
\end{table} 

The coherence axioms, (C9) and (C10) of Appendix~\ref{app:coh} have been
shown to be derivable from the others in~\cite[Theorem 2.3 and Proposition 4.2]{MarmolejoF:cohplr}.
Indeed, axiom (C9) is a particular case of a provable coherence condition for a
pseudomonad morphism and follows from (C1) and (C2) 
({\em cf.} Proposition~\ref{thm:marwood}). By duality, one can see 
that axiom~(C10) is a particular case of a provable coherence condition for a 
pseudomonad op-morphism and  follows from (C7) and (C8).
\end{rmk}


%
%
%

The explanation of the axioms for a pseudodistributive law in Remark~\ref{thm:explain-coherence} proves the following 
straightforward, but important, proposition. 

\begin{prop}
\label{prop:distr-law-as-obj}
The objects of $\Psm(\Psm(\catK))$ are exactly pseudodistributive laws in~$\catK$. 
\end{prop}

\begin{proof}
An object of $\Psm(\Psm(\catK))$ consists of an object $(X,S)$ of $\Psm(\catK)$, {\em i.e.} a pseudomonad
in $\cat{K}$,  together with a pseudomonad $(T, d) \colon (X,S) \to (X,S)$ on it in~$\Psm(\catK)$, which is
exactly a pseudodistributive law by Remark~\ref{thm:explain-coherence}.
\end{proof}

We can now give a new simple proof  of Marmolejo's fundamental result asserting the
equivalence between a pseudodistributive law of a pseudomonad~$T$ over
a pseudomonad~$S$ and a lifting of the pseudomonad~$T$ to the 2-categories of
pseudoalgebras for~$S$~\cite{MarmolejoF:dislp}.

\begin{theorem}
\label{thm:distr-iff-lift}
Let $\cat{K}$ be a Gray-category, $(X, S)$ and $(X,T)$ be pseudomonads in $\catK$. A pseudodistributive law $d\colon ST\rightarrow TS$ is equivalent to a lifting of $T$ to pseudo-$S$-algebras. 
\end{theorem}

\begin{proof} By Theorem~\ref{thm:psm-gray}, $\Psm(\catK)$ is a Gray-category and therefore we can consider the Gray-category $\Psm(\Psm(\catK))$.  

Next, observe that that $\Psm(-)$ preserves triequivalences between Gray-categories, i.e.\ given a triequivalence of Gray-categories $\Phi \colon \catK \to \catK'$, then it is possible
to define a triequivalence $\Psm(\Phi) \colon \Psm(\catK) \to \Psm(\catK')$. 
The construction of $\Psm(\Phi)$ is evident, and the verification that it
is a triequivalence is a long, but routine, calculation. For example, to prove essential surjectivity, we need to show that for every pseudomonad~$(X', T')$ in $\catK'$, there is a pseudomonad~$(X,T)$ in~$\catK$ that is 
mapped by~$\Psm(\Phi)$ to a pseudomonad that is biequivalent to~$(X', T')$ in~$\Psm(\catK')$.
For this, one defines~$(X,T)$ using the essential surjectivity of~$\Phi$, carefully inserting
 coherence isomorphisms that are part of the given triequivalence where appropriate. 

Applying this fact to the triequivalence of Theorem~\ref{thm:equiv-lift-psm}, we get a triequivalence:
\[
\Psm(\Psm(\catK))\simeq \Psm(\Lift(\catK)) \, .
\]
An object on the left hand side is
exactly a pseudodistributive law by Proposition~\ref{prop:distr-law-as-obj}. Similarly, an object on the
right hand side consists exactly of a pseudomonad~$(X,S)$ in~$\cat{K}$ and a pseudomonad 
$T \colon X \to X$ with a lifting $\hat{T} \colon \psSalg \to \psSalg$. 
\end{proof}

We conclude the paper by outlining how duality can be applied as in \cite[Section~4]{StreetR:fortm}
to obtain an equivalence between
pseudodistributive laws and extensions to Kleisli objects. Fix a Gray-category $\cat{K}$ and let $(X,T)$ 
be a pseudomonad in it. By definition, a \emph{right pseudo-$T$-module} in $\cat{K}$ is a left $T$-module 
in~$\cat{K}^{\op}$. We then have a Gray-functor
\begin{equation}
\label{equ:right-modules}
\mathrm{Mod}_T \colon \cat{K}^{\op} \to \Gray \, .
\end{equation}
Assuming the evident definition of a lifting 
to 2-categories of right pseudomodules, we have the following
corollary of Theorem~\ref{thm:distr-iff-lift}.

%
%

\begin{cor}
\label{cor:distr-leftmod}
Let $(X, S)$ and $(X,T)$ be pseudomonads in $\catK$. A pseudodistributive law~$d:ST\rightarrow TS$ is equivalent to a lifting of $S$ to right pseudo-$T$-modules.  \qed
\end{cor}

The equivalence of Corollary~\ref{cor:distr-leftmod} becomes more familiar under the assumption 
that~$\catK$ has
Kleisli objects. Recall that a \emph{Kleisli object} for a pseudomonad $(X,T)$ in $\cat{K}$ is an 0-cell $X_T \in \cat{K}$ and a right pseudo-$T$-module
$J_T \colon X \to X_T$, which is universal in the sense that the 2-functor
\[
 \cat{K}(X_T, I)  \to \mathrm{Mod}_T(I) \, , 
\]
induced by composition with $J_T$, is an equivalence of 2-categories, thus making the Gray-functor in~\eqref{equ:right-modules} representable. Now, 
a pseudodistributive law $d:ST\rightarrow TS$ is equivalent to a lifting of $S$ to right pseudo-$T$-modules, as in 
\begin{center}
\begin{tikzpicture}
\node (TF') at (3.5,2) {$\textrm{Mod}_T$};
\node (TF) at (0,2) {$\textrm{Mod}_T$};
\node (F'S) at (3.5,0) {$\catK(-,\,X) \, .$};
\node (FS) at (0,0) {$\catK(-,\,X)$};
\draw[->] (TF) to node[above] {$\lift{S}$} (TF');
\draw[->] (TF) to node[left] {$U$} (FS);
\draw[->] (FS) to node[below] {$S\circ-$} (F'S);
\draw[->] (TF') to node[right] {$U$} (F'S);
\end{tikzpicture}
\end{center}
This, in turn, is equivalent to 
\begin{center}
\begin{tikzpicture}
\node (TF') at (2,2) {$X_T$};
\node (TF) at (0,2) {$X_T$};
\node (F'S) at (2,0) {$X \, ,$};
\node (FS) at (0,0) {$X$};

\draw[->] (TF) to node[above] {$\lift{S}$} (TF');
\draw[->] (FS) to node[left] {$J_T$} (TF);
\draw[->] (FS) to node[below] {$S$} (F'S);
\draw[->] (F'S) to node[right] {$J_T$} (TF');
\end{tikzpicture}
\end{center}
which describes an extension of $S$ to the Kleisli object of $T$.

\begin{rmk} \label{rem:kleisli} 
We conclude the paper by briefly discussing the question of whether the Gray-category $\Gray$ has Kleisli objects. Given a 2-category $X$ and pseudomonad $T \colon X \to X$, there are two reasonable 
options to be the Kleisli object for $T$, mirroring the one-dimensional situation. In both
cases, the objects are the same objects as those of~$X$, but they have different hom-categories of morphisms. The first option is to define the hom-category between two objects $x, y \in X$ to be $X(x, Ty)$. With this definition we only get a bicategory,
 not a 2-category, and so we step outside~$\Gray$. The second option (which we will call $X_T$), is to take the hom-category of morphisms between $x$ and $y$ to consist of pseudoalgebras morphisms from $Tx$ to $Ty$ (considered as free algebras). 
This is a 2-category and so one could try to show that it is a Kleisli object for $\Gray$. In order to do this, 
one should prove that, for any 2-category $I$, there is an equivalence as in~\eqref{equ:right-modules}.
However, the construction taking  a $I$-indexed right  pseudo-$T$-module to a 2-functor 
$X_T \rightarrow I$ is only a pseudofunctor and not a strict 2-functor, thus leading again outside $\Gray$. The reason for this is that we need to use the pseudonaturality of the module action $\lambda$ and other coherence isomorphisms. 
Because of this, it seems that $\Gray$ does not have Kleisli objects. We suspect that, once it is defined what it means for a tricategory to have Kleisli objects, it should be possible to show that the tricategory $\mathbf{2\text{-}Cat}_{\mathrm{psd}}$ of 2-categories, pseudofunctors, pseudonatural transformations and modifications has Kleisli objects. The same should hold also for $\mathbf{Bicat}$, the tricategory 
of bicategories, pseudofunctors, pseudonatural transformation and modifications. We leave the investigation of these problems to future research.

\end{rmk}

\appendix

\section{Coherence conditions for pseudodistributive laws}
\label{app:coh}
We limit ourselves to drawing the boundaries
of these diagrams and explain in text which 3-cells should be inserted in them, except for the
3-cells coming from the structure of a Gray-category of $\catK$.

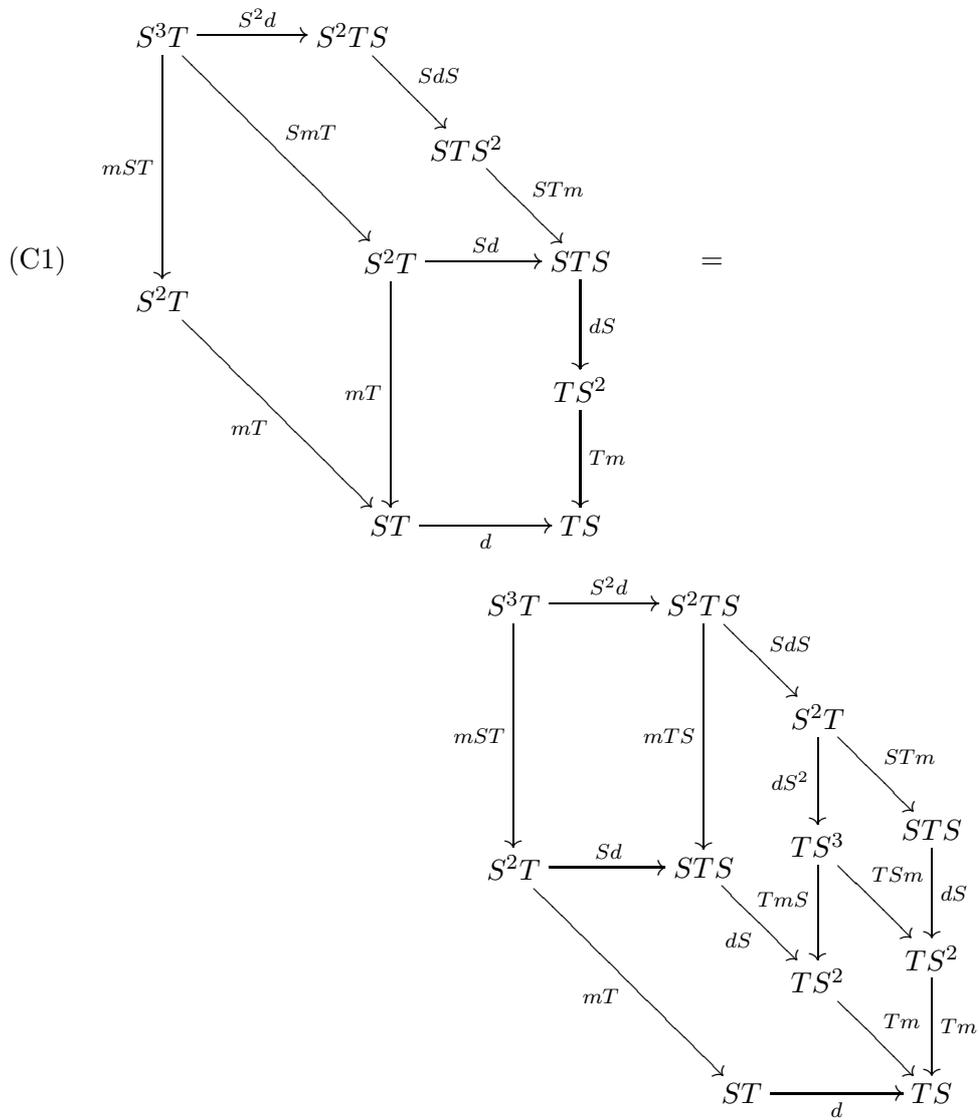
\begin{multline}
\tag{C1}
\begin{xy} 
(0,30)*+{\mssst}="1";
(25,30)*+{\mssts}="2";
(40,15)*+{\mstss}="6";
(0,-5)*+{\msst}="3";
(30,0)*+{\msst}="9";
(55,0)*+{\msts}="10";
(30,-35)*+{\mst}="11";
(55,-35)*+{\mts}="12";
(55,-17)*+{\mtss}="13";
(17,12)*+{}="31";
(23,6)*+{}="32";
(28,-7)*+{}="35";
(28,-14)*+{}="36";
(7,0)*+{}="37";
(7,-6)*+{}="38";
{\ar^{\mss \di} "1";"2"};
{\ar_{\mus \mst} "1";"3"};
{\ar^{\ms \mus \mt} "1";"9"};
{\ar^{\ms \di \ms} "2";"6"};
{\ar_{\mus \mt} "3";"11"};
{\ar^{\mst \mus} "6";"10"};
{\ar^{\ms \di} "9";"10"};
{\ar_{\mus \mt} "9";"11"};
{\ar^{\di \ms} "10";"13"};
{\ar^{\mt \mus} "13";"12"};
{\ar_{\di} "11";"12"};
\end{xy} \qquad = \\ 
\begin{xy} 
(0,30)*+{\mssst}="1";
(25,30)*+{\mssts}="2";
(0,-5)*+{\msst}="3";
(25,-5)*+{\msts}="4";
(40,15)*+{\msst}="6";
(12,6)*+{};
(40,-20)*+{\mtss}="8";
(40,-2)*+{\mtsss}="14";
(55,0)*+{\msts}="10";
(30,-35)*+{\mst}="11";
(55,-35)*+{\mts}="12";
(55,-17)*+{\mtss}="13";
(7,0)*+{}="17";
(7,-6)*+{}="18";
(15,-8)*+{}="30";
(20,-14)*+{}="31";
(33,-2)*+{}="32";
(33,-7)*+{}="33";
{\ar^{\mss \di} "1";"2"};
{\ar_{\mus \mst} "1";"3"};
{\ar_{\mus \mts} "2";"4"};
{\ar^{\ms \di \ms} "2";"6"};
{\ar^{\ms \di} "3";"4"};
{\ar_{\mus \mt} "3";"11"};
{\ar_{\di \ms} "4";"8"};
{\ar_{\di \mss} "6";"14"};
{\ar_(.4){\mt \mus \ms} "14";"8"};
{\ar^(.4){\mts \mus} "14";"13"};
{\ar^{\mst \mus} "6";"10"};
{\ar^{\mt \mus} "8";"12"};
{\ar^{\di \ms} "10";"13"};
{\ar^{\mt \mus} "13";"12"};
{\ar_{\di} "11";"12"};
\end{xy}
\end{multline} 
In~(C1), the left-hand side pasting is obtained using $\ms \bar{m}$,
$\bar{m}$, and the associativity 3-cell of the pseudomonad $\ms$; the right-hand side pasting is obtained using  the associativity  3-cell of the pseudomonad $\ms$ and $\bar{m}$.

\begin{equation}
\tag{C2}
\begin{xy}
(0,15)*+{\mst}="1";
(15,0)*+{\msst}="5";
(40,0)*+{\msts}="6";
(40,-18)*+{\mtss}="7";
(15,-35)*+{\mst}="9"; 
(40,-35)*+{\mts}="10";
(28,-7)*+{}="15";
(28,-14)*+{}="21";
{\ar@/^1pc/^{\mst \es} "1";"6"};
{\ar@/_1pc/_{1_{\mst}} "1";"9"};
{\ar^{\di \ms} "6";"7"};
{\ar^{\mt \mus} "7";"10"};
{\ar^{\mus \mt} "5";"9"};
{\ar^{\ms \es \mt} "1";"5"};
{\ar^{\ms \di} "5";"6"};
{\ar_{\di} "9";"10"};
(20,9)*+{}="30";
(20,4)*+{}="31";
(9,-2)*+{}="36";
(9,-7)*+{}="37";
\end{xy}
\quad = \quad  
\begin{xy} 
(0,15)*+{\mst}="1";
(0,-3)*+{\mts}="8";
(40,0)*+{\msts}="6";
(40,-18)*+{\mtss}="7";
(40,-35)*+{\mts}="10";
{\ar@/^1pc/^{\mst \es} "1";"6"};
{\ar_{\di} "1";"8"};
{\ar^{\mts \es} "8";"7"};
{\ar^{\di \ms} "6";"7"};
{\ar^{\mt \mus} "7";"10"};
{\ar@/_1pc/_{\mathit{Id}} "8";"10"};
\end{xy}
\end{equation}
In~(C2), the left-hand side pasting is obtained using $\ms \bar{\es}$,
$\bar{m}$, and the left unit  3-cell of the pseudomonad $\ms$; the
right-hand side pasting is obtained using the left unit  3-cell of the pseudomonad $\ms$.

\begin{equation}
\tag{C3}
\begin{xy}
(0,15)*+{\mss}="1";
(25,15)*+{\msst}="2";
(0,-25)*+{\ms}="16";
(40,0)*+{\msts}="6";
(40,-20)*+{\mtss}="10";
(40,-40)*+{\mts}="11";
{\ar^{\mss \et} "1";"2"};
{\ar_{\mus} "1";"16"};
{\ar^{\ms \di} "2";"6"};
{\ar@/_1pc/_{\ms \et \ms} "1";"6"};
{\ar@/_2pc/_(.65){\et \mss} "1";"10"};
{\ar@/_1pc/_{\et \ms} "16";"11"};
{\ar^{\mt \mus} "10";"11"};
{\ar^{\di \ms} "6";"10"};
(20,10)*+{}="30";
(14,10)*+{}="31";
\end{xy}
 = 
\begin{xy}
(0,15)*+{\mss}="1";
(25,15)*+{\msst}="2";
(0,-25)*+{\ms}="8";
(40,0)*+{\msts}="6";
(40,-20)*+{\mtss}="10";
(40,-40)*+{\mts}="11";
(25,-25)*+{\mst}="19";
(12,8)*+{}="15";
(12,4)*+{}="16";
(-13,0)*+{}="17";
(-13,-6)*+{}="18";
{\ar^{\mss \et} "1";"2"};
{\ar^{\ms \di} "2";"6"};
{\ar^{\di \ms} "6";"10"};
{\ar^{\mt \mus } "10";"11"};
{\ar^{\di} "19";"11"};
{\ar_{\mus} "1";"8"};
{\ar^{\ms \et} "8";"19"};
{\ar_{\mus \mt} "2";"19"};
(20,-30)*+{}="30";
(14,-30)*+{}="31";
(33,-2)*+{}="32";
(33,-7)*+{}="33";
{\ar@/_1pc/_{\et \ms} "8";"11"};
\end{xy}
\end{equation}
For~(C3), the left-hand side pasting is obtained using $\ms \bar{\et}$,
$\bar{\et} \ms$; the 
right-hand side is obtained using $\bar{\mus}$ and $\bar{\mus}$.

\begin{equation}
\tag{C4}
\begin{xy}
(0,15)*+{1_X}="1";
(25,15)*+{\mt}="2";
(0,-5)*+{\ms}="8";
(25,-5)*+{\mst}="19";
(40,-20)*+{\mts}="10";
(7,0)*+{}="17";
(7,-6)*+{}="18";
(15,-8)*+{}="30";
(20,-14)*+{}="31";
(33,-2)*+{}="34";
(33,-7)*+{}="35";
{\ar^{\et} "1";"2"};
{\ar^{\es \mt} "2";"19"};
{\ar_{\di} "19";"10"};
{\ar_{\es} "1";"8"};
{\ar^{\ms \et} "8";"19"};
{\ar@/^1pc/^{\mt \es} "2";"10"};
{\ar@/_/_{\et \ms} "8";"10"};
\end{xy}
 = 
\begin{xy}
(0,15)*+{1_X}="1";
(25,15)*+{\mt}="2";
(0,-5)*+{\ms}="8";
(40,-20)*+{\mts}="10";
{\ar^{\et} "1";"2"};
{\ar@/^1pc/^{\mt \es} "2";"10"};
{\ar_{\es} "1";"8"};
{\ar@/_/_{\et \ms} "8";"10"};
\end{xy}
\end{equation}
For~(C4), the left-hand side pasting is obtained using $\bar{\es}$ and $\bar{\et}$; the right-hand side is obtained from pseudonaturality of $t$.

\begin{multline} 
\tag{C5}
\begin{xy} 0;/r.20pc/:
(0,30)*+{\msstt}="1";
(30,30)*+{\msst}="2";
(15,-5)*+{\mtsst}="3";
(15,15)*+{\mstst}="7";
(30,0)*+{\mstts}="9";
(60,0)*+{\msts}="10";
(30,-40)*+{\mtt \mss}="11";
(60,-40)*+{\mtss}="12";
(30,-20)*+{\mtsts}="13";
(0,-30)*+{\mstt}="31";
(15,-45)*+{\mtst}="34";
(30,-60)*+{\mtts}="32";
(60,-60)*+{\mts}="33";
(28,-7)*+{}="35";
(28,-14)*+{}="36";
(7,0)*+{}="37";
(7,-6)*+{}="38";
{\ar^{\mss \mut} "1";"2"};
{\ar_{\mus \mtt} "1";"31"};
{\ar^{\ms \di \mt} "1";"7"};
{\ar^{\mst \di} "7";"9"};
{\ar^{\ms \di} "2";"10"};
{\ar^{\ms \mut \ms} "9";"10"};
{\ar^{\di \mts} "9";"13"};
{\ar^{\mt \di \ms} "13";"11"};
{\ar^{\di \ms} "10";"12"};
{\ar^{\mut \mss} "11";"12"};
{\ar_{\di \mst} "7";"3"};
{\ar_{\mt \mus \mt} "3";"34"};
{\ar_(.6){\mts \di} "3";"13"};
{\ar^{\mtt \mus} "11";"32"};
{\ar^{\mt \mus} "12";"33"};
{\ar_{\di \mt} "31";"34"};
{\ar_{\mt \di} "34";"32"};
{\ar_{\mut \ms} "32";"33"};
\end{xy} \qquad = \\
\begin{xy} 0;/r.20pc/:
(0,15)*+{\msstt}="1";
(25,15)*+{\msst}="2";
(25,-20)*+{\mst}="4";
(55,-15)*+{\msts}="10";
(55,-30)*+{\mtss}="12";
(0,-20)*+{\mstt}="31";
(30,-50)*+{\mtts}="32";
(55,-50)*+{\mts}="33";
(15,-35)*+{\mtst}="34";
(28,-7)*+{}="35";
(28,-14)*+{}="36";
(7,0)*+{}="37";
(7,-6)*+{}="38";
{\ar^{\mss \mut} "1";"2"};
{\ar^{\ms \di} "2";"10"};
{\ar^{\di \ms} "10";"12"};
{\ar^{\mt \mus} "12";"33"};
{\ar_{\mus \mtt} "1";"31"};
{\ar_{\di \mt} "31";"34"};
{\ar^{\ms \mut} "31";"4"};
{\ar_{\mt \di} "34";"32"};
{\ar_{\mut \ms} "32";"33"};
{\ar_{\mus \mt} "2";"4"};
{\ar^(.4){\di} "4";"33"};
\end{xy}
\end{multline} 
For~(C5), the left-hand side pasting is obtained using 
$\ms \bar{\mut}$, $\bar{\mut} \ms$ and $\bar{\mus} \mt$;
the right-hand side pasting is obtained using
$\bar{\mus}$ and $\bar{\mut}$.

\begin{equation}
\tag{C6}
\begin{xy}
(0,20)*+{\mtt}="1";
(25,20)*+{\mt}="2";
(0,0)*+{\mstt}="8";
(25,0)*+{\mst}="19";
(10,-15)*+{\mtst}="9"; 
(20,-30)*+{\mtts}="12"; 
(45,-30)*+{\mts}="13";
{\ar^{\mut} "1";"2"};
{\ar^{\ms \mut} "8";"19"};
{\ar^{\mut \ms} "12";"13"};
{\ar_{\es \mtt} "1";"8"};
{\ar^{\es \mt} "2";"19"};
{\ar_{\di \mt} "8";"9"};
{\ar^{\mt \di} "9";"12"};
{\ar_{\di} "19";"13"};
{\ar@/^2pc/^{\mt \es} "2";"13"};
\end{xy}
= 
\begin{xy}
(0,20)*+{\mtt}="1";
(25,20)*+{\mt}="2";
(0,0)*+{\mstt}="8";
(10,-15)*+{\mtst}="9"; 
(20,-30)*+{\mtts}="12"; 
(45,-30)*+{\mts}="13";
{\ar_{\es \mtt} "1";"8"};
{\ar_{\mt \di} "9";"12"};
{\ar^{\mut} "1";"2"};
{\ar@/^2pc/^{\mt \es} "2";"13"};
{\ar@/^2pc/^{\mtt \es} "1";"12"};
{\ar_{\di \mt} "8";"9"};
{\ar_{\mut \ms}  "12";"13"};
{\ar@/^.5pc/^(.65){\mt \es \mt} "1";"9"};
\end{xy}
\end{equation}
In~(C6), the left-hand side pasting is obtained using $\bar{\es}$ and $\bar{\mut}$; the 
right-hand side pasting is obtained using $\bar{\es} \mt$. 



\medskip

\begin{equation}
\tag{C7}
\begin{xy}
(0,24)*+{\msttt}="1";
(40,24)*+{\mst}="2";
(0,8)*+{\mtstt}="3";
(0,-8)*+{\mttst}="5";
(0,-24)*+{\mttts}="7";
(40,-24)*+{\mts}="8";
(20,32)*+{\mstt}="9";
(20,16)*+{\msst}="10";
(20,-16)*+{\mtst}="12";
(20,-32)*+{\mtts}="13";
{\ar_{\di \mtt} "1";"3"};
{\ar_{\mt \di \mt} "3";"5"};
{\ar_{\mtt \di} "5";"7"};
{\ar^{\di} "2";"8"};
{\ar_{\di \mt} "10";"12"};
{\ar_{\mt \di} "12";"13"};
{\ar@/^/^{\mst \mut} "1";"9"};
{\ar@/^/^{\ms \mut} "9";"2"};
{\ar@/_/_{\ms \mut \mt} "1";"10"};
{\ar@/_/_{\mut \mst} "5";"12"};
{\ar@/_/_{\mut \mts} "7";"13"};
{\ar@/_/_{\ms \mut} "10";"2"};
{\ar@/_/_{\mut \ms} "13";"8"};
\end{xy}
\quad = \quad 
\begin{xy}
(0,24)*+{\msttt}="1";
(40,24)*+{\mst}="2";
(0,8)*+{\mtstt}="3";
(0,-8)*+{\mttst}="5";
(0,-24)*+{\mttts}="7";
(40,-24)*+{\mts}="8";
(20,32)*+{\mstt}="9";
(20,16)*+{\mtst}="10";
(20,-16)*+{\mtts}="12";
(20,-32)*+{\mtts}="13";
{\ar_{\di \mtt} "1";"3"};
{\ar_{\mt \di \mt} "3";"5"};
{\ar_{\mtt \di} "5";"7"};
{\ar^{\di} "2";"8"};
{\ar_{\di \mt} "9";"10"};
{\ar_{\mt \di} "10";"12"};
{\ar@/^/^{\mst \mut} "1";"9"};
{\ar@/^/^{\ms \mut} "9";"2"};
{\ar@/^/^{\mts \mut} "3";"10"};
{\ar@/^/^{\mt \mut \ms} "7";"12"};
{\ar@/_/_{\mut \mts} "7";"13"};
{\ar@/^/^{\mut \ms} "12";"8"};
{\ar@/_/_{\mut \ms} "13";"8"};
\end{xy}
\end{equation}
For~(C7), the left-hand side pasting is obtained using the associativity 3-cell of the pseudomonad $T$, 
$\bar{\mut}$ and $\bar{\mut} \mt$; the right-hand side pasting is obtained using
$T \bar{\mut}$, $\bar{\mut}$ and the associativity 3-cell of the pseudomonad $T$.

\begin{equation}
\tag{C8}
\begin{xy}
(0,24)*+{\mst}="1";
(40,24)*+{\mst}="2";
(0,-24)*+{\mts}="7";
(40,-24)*+{\mts}="8";
(20,16)*+{\mstt}="10";
(20,-8)*+{\mtst}="12";
(20,-32)*+{\mtts}="13";
{\ar_{\di} "1";"7"};
{\ar^{\di} "2";"8"};
{\ar^{\di \mt} "10";"12"};
{\ar^{\mt \di} "12";"13"};
{\ar@/^2pc/^{1_{\mst}} "1";"2"};
{\ar@/_/^{\ms \et \mt} "1";"10"};
{\ar@/_/_(.65){\et \mst} "1";"12"};
{\ar@/_/_{\et \mts} "7";"13"};
{\ar@/_/^{\ms \mut} "10";"2"};
{\ar@/_/_{\mut \ms} "13";"8"};
\end{xy}
\qquad = \quad 
\begin{xy}
(0,24)*+{\mst}="1";
(0,-24)*+{\mts}="7";
(40,-24)*+{\mts}="8";
(20,-32)*+{\mtts}="13";
{\ar_{\di} "1";"7"};
{\ar@/^2pc/^{1_{\mts}} "7";"8"};
{\ar@/_/_{\et \mts} "7";"13"};
{\ar@/_/_{\mut \ms} "13";"8"};
\end{xy}
\end{equation}
For~(C8), the left-hand side pasting is obtained using the right unit 3-cell of the pseudomonad~$T$, 
$\bar{\mut}$, $\bar{\et} \ms$; the right-hand side pasting is  the right unit 3-cell of the pseudomonad~$T$.

\begin{equation}
\tag{C9}
\begin{xy}
(45,-20)*+{\mst}="1";
(65,-20)*+{\mts}="2";
(5,20)*+{\mst}="3";
(30,20)*+{\mts}="4";
(5,0)*+{\msst}="5";
(30,0)*+{\msts}="6";
(45,-8)*+{\mtss}="7";
(15,14)*+{}="8";
(15,8)*+{}="9";
(36,-5)*+{}="10";
(44,-11)*+{}="11";
(38,8)*+{}="12";
(38,2)*+{}="13";
(62,0)*+{}="14";
(62,-6)*+{}="15";
{\ar_{\di} "1";"2"};
{\ar^{\di} "3";"4"};
{\ar_{\es \mst} "3";"5"};
{\ar_{\mus \mt} "5";"1"};
{\ar^(.35){\mt \mus} "7";"2"};
{\ar^(.4){\di \ms} "6";"7"};
{\ar_{\es \mts} "4";"6"};
{\ar^{\ms \di} "5";"6"};
{\ar@/^/^(.55){\mt \es \ms}  "4";"7"};
{\ar@/^2pc/^{1_{\mts}} "4";"2"};
\end{xy} = 
\begin{xy}
(25,-20)*+{\mst}="1";
(45,-20)*+{\mts}="2";
(5,20)*+{\mst}="3";
(5,0)*+{\msst}="5";
(20,6)*+{}="6";
(20,-2)*+{}="7";
{\ar_{\di} "1";"2"};
{\ar_{\es \mst} "3";"5"};
{\ar_{\mus \mt} "5";"1"};
{\ar@/^1pc/^{1_{\mst}} "3";"1"};
\end{xy} 
\end{equation}
For~(C9), the left-hand side pasting is obtained using the
right unit 3-cell of the pseudomonad~$S$,  $\bar{\es}S$ and $\bar{\mus}$;
 the right-hand side pasting is obtained
using the right unit 3-cell of the pseudomonad~$S$.

\begin{equation}
\tag{C10}
\begin{xy}
(0,24)*+{\mst}="1";
(40,24)*+{\mst}="2";
(40,-24)*+{\mts}="8";
(20,32)*+{\mstt}="9";
{\ar^{\di} "2";"8"};
{\ar@/^/^{\mst \et} "1";"9"};
{\ar@/^/^{\ms \mut} "9";"2"};
{\ar@/_2pc/_{1_{\mst}} "1";"2"};
\end{xy}
\qquad = \quad 
\begin{xy}
(0,24)*+{\mst}="1";
(40,24)*+{\mst}="2";
(0,-24)*+{\mts}="7";
(40,-24)*+{\mts}="8";
(20,32)*+{\mstt}="9";
(20,8)*+{\mtst}="11";
(20,-16)*+{\mtts}="12";
{\ar_{\di} "1";"7"};
{\ar^{\di} "2";"8"};
{\ar_{\di \mt} "9";"11"};
{\ar_{\mt \di} "11";"12"};
{\ar@/^/^{\mst \et} "1";"9"};
{\ar@/^/^{\ms \mut} "9";"2"};
{\ar@/^/^(.65){\mts \et} "7";"11"};
{\ar@/^/_{\mt \et \ms} "7";"12"};
{\ar@/_2pc/_{1_{\mts}} "7";"8"};
{\ar@/^/_{\mut \ms} "12";"8"};
\end{xy}
\end{equation}
For~(C10), the left-hand side pasting uses the left unit 3-cell of the pseudomonad $T$. The right-hand
side pasting is obtained using $\bar{\mut}$ and  the left unit 3-cell of the pseudomonad~$T$.

\section{Some technical proofs}
\label{app:coherence-comp}
\label{app:well-defineness-Fphi}
\label{app:coherence-interch}

\subsection*{Coherence diagrams for $(GF,\,G\phi\cdot\psi F)$}
We show only equation in~\eqref{equ:m2}. 
Using the coherence diagram (\ref{equ:m2}) for $(G,\,\psi)$ and for $(F,\,\phi)$, it suffices to prove that the following two diagrams are equal: 
\begin{center}
\begin{tikzpicture}[scale=0.95]
\node (b) at (3,6) {$QGTF$};
\node (d) at (3,4) {$GT^2F$};
\node (f) at (3,0) {$GTF$};
\node (QGFS) at (6,6) {$QGFS$};
\node (GTFS) at (6,4) {$GTFS$};
\node (GFS2) at (6,2) {$GFS^2$};
\node (GFS) at (6,0) {$GFS$};
\node (QGF) at (-2,8) {$QGF$};
\node (GTF) at (-2,5) {$GTF$}; 
\draw[->] (QGF) [bend left=10] to node[scale=.7]  [above,xshift=3, yshift=0.1] {$QGtF$} (b);
\draw[->] (QGF) [bend left=25] to node[scale=.7]  (QGFs) [above] {$QGFs$} (QGFS);
\draw[->] (QGF) to node[scale=.7] [left] {$\psi F$} (GTF);
\draw[->] (GTF) [bend right=25] to node[scale=.7]  [left] {$1_{GTF}$} (f);
\draw[->] (GTF) to node[scale=.7]  [above] {$GTtF$} (d);

\draw[->] (b) to node[scale=.7] (G'f) [left] {$\psi TF$} (d);
\draw[->] (d) to node[scale=.7] (Gf2) [left]{$G\mut F$} (f);
\draw[->] (b) to node[scale=.7] (QGphi) [above]{$QG\phi$} (QGFS);
\draw[->] (QGFS) to node[scale=.7] (psiFS) [right]{$\psi FS$} (GTFS);

\draw[->] (d) to node[scale=.7] (gF2) [below]{$GT\phi$} (GTFS);
\draw[->] (GTFS) to node[scale=.7] (G'f2) [right] {$G\phi S$} (GFS2);
\draw[->] (GFS2) to node[scale=.7] (G'f2) [right] {$GF\mus $} (GFS);
\draw[->] (f) to node[scale=.7] (gF''2) [below] {$G\phi$} (GFS);
\draw[-{Implies},double distance=1.5pt,shorten >=40pt,shorten <=40pt] (gF2) to node[scale=.7] [left] {$G\bar{\phi}$} (gF''2);
\draw[-{Implies},double distance=1.5pt,shorten >=20pt,shorten <=20pt] (QGphi) to node[scale=.7] [left] {${\psi_\phi}$} (gF2);

\node (psi1) at (0.5,6.5) {};
\node (psi2) at (0.5,5.5) {};
\draw[-{Implies},double distance=1.5pt] (psi1) to node[scale=.7] [right] {$\psi_{tF}$} (psi2);

\node (eta1) at (0.5,3.5) {};
\node (eta2) at (0.5,2.5) {};
\draw[-{Implies},double distance=1.5pt] (eta1) to node[scale=.7] [right] {$G\eta_TF$} (eta2);

\node (phi1) at (3,7.5) {};
\node (phi2) at (3,6.5) {};
\draw[-{Implies},double distance=1.5pt] (phi1) to node[scale=.7] [right] {$QG\widetilde{\phi}$} (phi2);
\end{tikzpicture} 
\end{center}
\begin{center}
\begin{tikzpicture}[scale=0.95]
\node (QGF) at (1,8) {$QGF$};
\node (GTF) at (1,6) {$GTF$};
\node (QGFS) at (6,7) {$QGFS$};
\node (GTFS) at (6,4) {$GTFS$};
\node (GFS2) at (6,2) {$GFS^2$};
\node (GFS) at (6,0) {$GFS$.};
\node (GT2F) at (3,4) {$GT^2F$}; 
\node (GFS') at (3,0) {$GTF$};

\path[->] 
(QGF) edge [bend left=15] node[scale=.7]  (QGFs) [above] {$QGFs$} (QGFS)
	edge node[scale=.7, left] {$\psi F$} (GTF)
(GTF) edge [bend left=15] node[scale=.7, above, xshift=0.5cm] (GTFs) {$GTFs$} (GTFS)
	edge [bend right=25] node[scale=.7, left] {$1_GTF$} (GFS')
	edge node [scale=.7, left,yshift=-0.3cm] {$GTtF$} (GT2F)
(GFS') edge node[scale=.7, below] {$G\psi$} (GFS)
(QGFS) edge node[right, scale=.7,] {$\psi FS$} (GTFS)
(GTFS) edge node[scale=.7] [right] {$G\phi S$} (GFS2)
(GT2F) edge node[scale=.7] [below] {$GT\phi$} (GTFS)
	edge node[scale=.7] [right] {$G\mut F$} (GFS')
(GFS2) edge node[scale=.7] [right] {$GF\mus $} (GFS);
\node (psis1) at (3.5,7.3) {};
\node (psis2) at (3.5,6.3) {};
\draw[-{Implies},double distance=1.5pt] (psis1) to node[scale=.7] [right] {$\psi F_s$} (psis2);

\node (phis1) at (3.5,5.3) {};
\node (phis2) at (3.5,4.3) {};
\draw[-{Implies},double distance=1.5pt] (phis1) to node[scale=.7] [right] {$GT\widetilde{\phi}$} (phis2);

\node (eta1) at (4.5,2.5) {};
\node (eta2) at (4.5,1.5) {};
\draw[-{Implies},double distance=1.5pt] (eta1) to node[scale=.7] [right] {$G\bar{\phi}$} (eta2);

\node (e1) at (1.8,3.5) {};
\node (e2) at (1.8,2.5) {};
\draw[-{Implies},double distance=1.5pt] (e1) to node[scale=.7] [right] {$G\eta_TF$} (e2);

\end{tikzpicture}
\end{center}
This equality holds using (\ref{equ:cc1}) and (\ref{equ:cc2}).

\subsection*{$F_\phi$ is well-defined}

Given a pseudomonad transformation $(q,\,\bar{q}):(G,\,\psi)\rightarrow(G',\,\psi')$ in $P_{\catK}(\,(Y,\,T),\,(Z,\,Q)\,)$ we want to show that $(qF,\,\overline{qF})$ is a pseudomonad transformation as well. We will show just 
equation~\eqref{equ:t1}, since~\eqref{equ:t2} can be proved similarly. 
The required equality follows from equation~\eqref{equ:t1} for $q$ and the equation below, which can be proved using  (\ref{equ:cc1}) twice.
\begin{center}
\begin{tikzpicture}[scale=0.85]
\node at (10,5) {$=$};
\node (b) at (0,11) {$QGTF$};
\node (d) at (0,8) {$GT^2F$};
\node (f) at (0,2) {$GTF$};
\node (QGFS) at (3,9) {$QGFS$};
\node (GTFS) at (3,6) {$GTFS$};
\node (GFS2) at (3,3) {$GFS^2$};
\node (GFS) at (3,0) {$GFS$};
\node (QG'TF) at (4,11) {$QG'TF$};
\node (QG'FS) at (7,9) {$QG'FS$}; 
\node (G'TFS) at (7,6) {$G'TFS$};
\node (G'FS2) at (7,3) {$G'FS^2$};
\node (G'FS) at (7,0) {$G'FS$}; 

\draw[->] (b) to node[scale=.7] (G'f) [left] {$\psi TF$} (d);
\draw[->] (d) to node[scale=.7] (Gf2) [left]{$G\mut F$} (f);
\draw[->] (b) to node[scale=.7] (QGphi) [above, yshift=0.1cm]{$QG\phi$} (QGFS);
\draw[->] (QGFS) to node[scale=.7] (psiFS) [right]{$\psi FS$} (GTFS);

\draw[->] (d) to node[scale=.7] (gF2) [below, yshift=-0.3cm]{$GT\phi$} (GTFS);
\draw[->] (GTFS) to node[scale=.7] (G'f2) [right] {$G\phi S$} (GFS2);
\draw[->] (GFS2) to node[scale=.7] (G'f2) [right] {$GF\mus $} (GFS);
\draw[->] (f) to node[scale=.7] (gF''2) [below, yshift=-0.2cm] {$G\phi$} (GFS);

\path[->] 
(b) edge node[scale=.7] (QqTF) [above] {$QqTF$} (QG'TF)
(QG'TF) edge node[scale=.7] [right, xshift=0.1cm, yshift=0.1cm] {$QG'\phi$} (QG'FS)
(QGFS) edge node[scale=.7] (QqFS) [above] {$QqFS$} (QG'FS)
(QG'FS) edge node[scale=.7] [right] {$\psi'FS$} (G'TFS)
(GTFS) edge node[scale=.7] (qTFS) [above] {$qTFS$} (G'TFS)
(G'TFS) edge node[scale=.7] [right] {$G'\phi S$} (G'FS2)
(GFS2) edge node[scale=.7] (qFS2) [above] {$qFS^2$} (G'FS2)
(G'FS2) edge node[scale=.7] [right] {$G'F\mus $} (G'FS)
(GFS) edge node[scale=.7] (qFS) [below] {$qFS$} (G'FS);

\draw[-{Implies},double distance=1.5pt,shorten >=50pt,shorten <=50pt] (gF2) to node[scale=.7] [left] {$G\bar{\phi}$} (gF''2);
\draw[-{Implies},double distance=1.5pt,shorten >=30pt,shorten <=30pt] (QGphi) to node[scale=.7] [left] {${\psi_\phi}$} (gF2);

\node (x1) at (2,11) {};
\node (y1) at (5,9) {};
\draw[-{Implies},double distance=1.5pt,shorten >=30pt,shorten <=30pt] (x1) to node[scale=.7] [right, xshift=0.3cm, yshift=0.2cm] {${Qq_\phi}^{-1}$} (y1);

\draw[-{Implies},double distance=1.5pt,shorten >=25pt,shorten <=20pt] (QqFS) to node[scale=.7] [right] {$\bar{q}FS$} (qTFS);
\draw[-{Implies},double distance=1.5pt,shorten >=25pt,shorten <=20pt] (qTFS) to node[scale=.7] [right] {${q_{\phi S}}^{-1}$} (qFS2);
\draw[-{Implies},double distance=1.5pt,shorten >=30pt,shorten <=30pt] (qFS2) to node[scale=.7] [right] {${q_{F\mus }}^{-1}$} (qFS);
\node at (18,0) {};
\end{tikzpicture} 
\begin{tikzpicture}[scale=0.85]
\node at (2,0) {};
\node at (8,8) {$=$};
\node (rb) at (15,11) {$QG'TF$};
\node (rd) at (15,8) {$G'T^2F$};
\node (rf) at (15,2) {$G'TF$};
\node (rQGFS) at (18,9) {$QG'FS$};
\node (rGTFS) at (18,6) {$G'TFS$};
\node (rGFS2) at (18,3) {$G'FS^2$};
\node (rGFS) at (18,0) {$G'FS$};
\node (rGTF) at (11,2) {$GTF$};
\node (rh) at (14,0) {$GFS$};
\node (rQG'TF) at (11,11) {$QG'TF$};
\node (rGT2F) at (11,8) {$GT^2F$}; 

\draw[->] (rb) to node[scale=.7] (rG'f) [left] {$\psi' TF$} (rd);
\draw[->] (rd) to node[scale=.7] (rGf2) [left]{$G'\mut F$} (rf);
\draw[->] (rb) to node[scale=.7] (rQGphi) [above, yshift=0.1cm, xshift=0.2cm]{$QG'\phi$} (rQGFS);
\draw[->] (rQGFS) to node[scale=.7] (rpsiFS) [right]{$\psi' FS$} (rGTFS);

\draw[->] (rd) to node[scale=.7] (rgF2) [below, yshift=-0.3cm]{$G'T\phi$} (rGTFS);
\draw[->] (rGTFS) to node[scale=.7] (rG'f2) [right] {$G'\phi S$} (rGFS2);
\draw[->] (rGFS2) to node[scale=.7] (rG'f2) [right] {$G'F\mus $} (rGFS);
\draw[->] (rf) to node[scale=.7] (rgF''2) [above, yshift=0.2cm] {$G'\phi$} (rGFS);

\path[->] 
(rQG'TF) edge node[scale=.7] (rQqTF) [above] {$QqTF$} (rb)
	edge node[scale=.7] [left, xshift=-0.2cm, yshift=-0.2cm] {$\psi TF$} (rGT2F)
(rGT2F) edge node[scale=.7] (rqT2F) [below] {$qT^2F$} (rd)
	edge node[scale=.7] [left, xshift=-0.2cm, yshift=-0.2cm] {$G\mut F$} (rGTF)	
	
(rGTF) edge node[scale=.7] (rqTF) [above] {$qTF$} (rf)
	edge node[scale=.7] [left, xshift=-0.2cm, yshift=-0.2cm] {$G\phi$} (rh)
(rh) edge node[scale=.7] (rqFS) [below] {$qFS$} (rGFS);

\draw[-{Implies},double distance=1.5pt,shorten >=40pt,shorten <=40pt] (rgF2) to node[scale=.7] [left] {$G'\bar{\phi}$} (rgF''2);
\draw[-{Implies},double distance=1.5pt,shorten >=30pt,shorten <=30pt] (rQGphi) to node[scale=.7] [left] {${\psi'_\phi}$} (rgF2);

\draw[-{Implies},double distance=1.5pt,shorten >=25pt,shorten <=30pt] (rqTF) to node[scale=.7] [right, xshift=0.5cm] {${q_\phi}^{-1}$} (rqFS);

\draw[-{Implies},double distance=1.5pt,shorten >=25pt,shorten <=25pt] (rQqTF) to node[scale=.7] [left] {$\bar{q}TF$} (rqT2F);
\draw[-{Implies},double distance=1.5pt,shorten >=45pt,shorten <=45pt] (rqT2F) to node[scale=.7] [left] {$q_{\mut F}$} (rqTF);
\end{tikzpicture}
\end{center}

Let 
\[
(q,\,\bar{q}),\,(q',\,\bar{q'}):(G,\,\psi)\rightarrow(G',\,\psi')
\] 
be pseudomonads transformations in $P_{\catK} (\,(Y,\,T),\,(Z,\,Q)\,)$. Given a pseudomonad modification $\beta:(q,\,\bar{q})\rightarrow(q',\,\bar{q'})$ we want to show that $\beta F$ is a pseudomonad modification from $(qF,\,\bar{qF})$ to $(q'F,\,\bar{q'F})$. So we need to show the following equation 

\begin{center}
\begin{tikzpicture}
\node (QGF) at (0,6) {$QGF$}; 
\node (QG'F) at (3,6) {$QG'F$};
\node (GTF) at (0,3) {$GTF$};
\node (G'TF) at (3,3) {$G'TF$};
\node (GFS) at (0,0) {$GFS$};
\node (G'FS) at (3,0) {$G'FS$};

\path[->]
(QGF) edge [bend left=20] node[scale=.7] (Qq'F) [above] {$QqF$} (QG'F)
	edge node[scale=.7] [left] {$\psi F$} (GTF)
	edge [bend right=20] node[scale=.7] (QqF) [below] {$Qq'F$} (QG'F)
(QG'F) edge node[scale=.7] [right] {$\psi'F$} (G'TF)
(GTF) edge [bend right=20] node[scale=.7] (q'TF) [above] {$q'TF$} (G'TF)
	edge node[scale=.7][left] {$G\phi$} (GFS)
(G'TF) edge node[scale=.7] [right] {$G'\phi$} (G'FS)
(GFS) edge [bend right=20] node[scale=.7] (q'FS) [below] {$q'FS$} (G'FS);

\draw[-{Implies},double distance=1.5pt,shorten >=5pt,shorten <=5pt] (Qq'F) to node[scale=.7] [right] {$Q\beta F$} (QqF);
\draw[-{Implies},double distance=1.5pt,shorten >=20pt,shorten <=20pt] (QqF) to node[scale=.7] [right] {$\bar{q'}F$} (q'TF);
\draw[-{Implies},double distance=1.5pt,shorten >=30pt,shorten <=30pt] (q'TF) to node[scale=.7] [right] {${q'_\phi}^{-1}$} (q'FS);

\node at (4.5,3) {$=$};

\node (QGFR) at (6,6) {$QGF$}; 
\node (QG'FR) at (9,6) {$QG'F$};
\node (GTFR) at (6,3) {$GTF$};
\node (G'TFR) at (9,3) {$G'TF$};
\node (GFSR) at (6,0) {$GFS$};
\node (G'FSR) at (9,0) {$G'FS$.};

\path[->]
(QGFR) edge [bend left=20] node[scale=.7] (QqFR) [above] {$QqF$} (QG'FR)
	edge node[scale=.7] [left] {$\psi F$} (GTFR)
(QG'FR) edge node[scale=.7] [right] {$\psi'F$} (G'TFR)
(GTFR) edge [bend left=20] node[scale=.7] (q'TFR) [above] {$qTF$} (G'TFR)
	edge node[scale=.7][left] {$G\phi$} (GFSR)
(G'TFR) edge node[scale=.7] [right] {$G'\phi$} (G'FSR)
(GFSR) edge [bend right=20] node[scale=.7] (qFSR) [below] {$q'FS$} (G'FSR)
	edge [bend left=20] node[scale=.7] (q'FSR) [above] {$qFS$} (G'FSR);

\draw[-{Implies},double distance=1.5pt,shorten >=5pt,shorten <=5pt] (q'FSR) to node[scale=.7] [right] {$\beta FS$} (qFSR);
\draw[-{Implies},double distance=1.5pt,shorten >=25pt,shorten <=25pt] (QqFR) to node[scale=.7] [right] {$\bar{q}F$} (q'TFR);
\draw[-{Implies},double distance=1.5pt,shorten >=25pt,shorten <=25pt] (q'TFR) to node[scale=.7] [right] {${q_\phi}^{-1}$} (q'FSR);
\end{tikzpicture}
\end{center}
This can be shown to hold using the coherence axiom for $\beta$ and (\ref{equ:cc3}).

\subsection*{Coherence for $q_p$}
Given $(p,\bar{p}):(F,\,\phi)\rightarrow(F',\,\phi')$ and $(q,\bar{q}):(G,\,\psi)\rightarrow(G',\,\psi')$ two 2-cells in $P_\catK$ we want to prove that $q_p:G'p\cdot qF\rightarrow qF'\cdot Gp$ is a 3-cell in $P_\catK$. First of all, it is useful to note that
\begin{center}
\begin{tikzpicture}[scale=0.85]
\node at (7,2) {$=$};
\node (GFS) at (0,0) {$GFS$};
\node (GTF) at (0,3) {$GTF$};
\node (GF'S) at (3,-1) {$GF'S$};
\node (G'FS) at (3,1) {$G'FS$};
\node (GTF') at (3,4) {$G'TF$};
\node (G'F'S) at (6,0) {$G'F'S$};
\node (G'TF') at (6,3) {$G'TF'$};

\path[->]
(GFS) edge [bend right=15] node[scale=.7] [below] {$GpS$} (GF'S)
	edge [bend left=15] node[scale=.7] (GpS) [above] {$qFS$} (G'FS)
(GTF) edge [bend left=15] node[scale=.7] (GTp) [above] {$qTF$} (GTF')
	edge node[scale=.7] [left] {$G\phi$} (GFS)
(GF'S) edge [bend right=15] node[scale=.7] [below] {$qF'S$} (G'F'S)
(GTF') edge [bend left=15] node[scale=.7] (qTF') [above] {$G'Tp$} (G'TF')
	edge node[scale=.7] [left] {$G'\phi$} (G'FS)
(G'TF') edge node[scale=.7] [right] {$G'\phi'$} (G'F'S)
(G'FS) edge [bend left=15] node[scale=.7] (qF'S) [above] {$G'pS$} (G'F'S);

\draw[-{Implies},double distance=1.5pt,shorten >=10pt,shorten <=10pt] (G'FS) to node[scale=.7] [right] {${q_{pS}}^{-1}$} (GF'S);
\draw[-{Implies},double distance=1.5pt,shorten >=20pt,shorten <=20pt] (GTp) to node[scale=.7] [right] {${q_\phi}^{-1}$} (GpS);
\draw[-{Implies},double distance=1.5pt,shorten >=20pt,shorten <=20pt] (qTF') to node[scale=.7] [right] {$G'\bar{p}$} (qF'S);
\end{tikzpicture}
\begin{tikzpicture}[scale=0.85]
\node (GTF) at (0,3) {$GFS$};
\node (QGF) at (0,6) {$GTF$};
\node (GTF') at (3,2) {$GF'S$};
\node (QGF') at (3,5) {$GTF'$};
\node (QG'F) at (3,7) {$G'TF$};
\node (G'TF') at (6,3) {$G'F'S$.};
\node (QG'F') at (6,6) {$G'TF'$};

\path[->]
(GTF) edge [bend right=15] node[scale=.7] (GTp) [below] {$GpS$} (GTF')
(QGF) edge [bend left=15] node[scale=.7] [above] {$qTF$} (QG'F)
	edge [bend right=15] node[scale=.7] (QGp) [below] {$GTp$} (QGF')
	edge node[scale=.7] [left] {$G\phi$} (GTF)
(GTF') edge [bend right=15] node[scale=.7] (qTF') [below] {$qF'S$} (G'TF')
(QGF') edge [bend right=15] node[scale=.7] (QqF') [below] {$qTF'$} (QG'F')
	edge node[scale=.7] [left] {$G\phi'$} (GTF')
(QG'F) edge [bend left=15] node[scale=.7] [above] {$G'Tp$} (QG'F')
(QG'F') edge node[scale=.7] [right] {$G'\phi'$} (G'TF');


\node (a) at (2.5,7) {};
\node (b) at (2.5,5) {};
\draw[-{Implies},double distance=1.5pt,shorten >=10pt,shorten <=10pt] (a) to node[scale=.7] [right] {${q_{Tp}}^{-1}={(qT)_p}^{-1}$} (b);
\draw[-{Implies},double distance=1.5pt,shorten >=20pt,shorten <=20pt] (QGp) to node[scale=.7] [right] {$G\bar{p}$} (GTp);
\draw[-{Implies},double distance=1.5pt,shorten >=20pt,shorten <=20pt] (QqF') to node[scale=.7] [right] {${q_{\phi'}}^{-1}$} (qTF');
\end{tikzpicture}
\end{center}

This equality is true since both diagrams are equal to the following one, using (\ref{equ:cc3}) for the one on the left-hand side and (\ref{equ:cc1}) for the one on the
right-hand side.  
\begin{center}
\begin{tikzpicture}
\node (a) at (0,3) {$GTF$};
\node (b) at (3,3) {$G'TF$};
\node (c) at (0,0) {$GF'S$};
\node (d) at (3,0) {$G'F'S$};
\draw[->] (a) to node[scale=.7] (f) [above] {$qTF$} (b);
\draw[->] (a) to node[scale=.7] [left] {$G(pS\cdot\phi)$} (c);
\draw[->] (b) [bend left] to node[scale=.7] (Gf) [right]{$G'(\phi'\cdot Tp)$} (d);
\draw[->] (b) [bend right] to node[scale=.7] (Gf') [left, yshift=-1cm]{$G'(pS\cdot\phi)$} (d);
\draw[->] (c) to node[scale=.7] (g'F) [below]{$g'F$} (d);
\draw[-{Implies},double distance=1.5pt,shorten >=40pt,shorten <=40pt] (b) to node[scale=.7] [left, xshift=-0.2cm, yshift=0.5cm] {${q_{(pS\cdot\phi)}}^{-1}$} (c);
\node (x') at (2.5,1.5) {};
\node (x) at (3.5,1.5) {};
\draw[-{Implies},double distance=1.5pt,shorten >=3pt,shorten <=3pt] (x) to node[scale=.7] [below, yshift=-0.2cm] {$G'\bar{p}$} (x');
\end{tikzpicture}
\end{center}
%
%
%
%
The proof can be concluded using the invertibility of the 3-cells involved and (\ref{equ:cc1}). 

%

\bibliographystyle{plain}

\begin{thebibliography}{10}

\bibitem{BarrM:toptt}
M.~Barr and C.~Wells.
\newblock {\em Toposes, triples, and theories}.
\newblock Springer, 1985.

\bibitem{BeckJ:disl}
J.~Beck.
\newblock Distributive laws.
\newblock In {\em Seminar on Triples and Categorical Homology Theory (ETH,
  Z\"urich, 1966/67)}, pages 119--140. Springer, 1969.
  
  

\bibitem{BungeM:cohera}
M.~Bunge.
\newblock Coherent extensions and relational algebras.
\newblock {\em Transactions of the American Mathematical Society}, 197:355--390,
  1974.
  

\bibitem{CattaniG:proomb}
G.~L. Cattani and G.~Winskel.
\newblock Profunctors, open maps, and bisimulation.
\newblock {\em Mathematical Structures in Computer Science}, 15(3):553--614,
  2005.

\bibitem{ChengE:psedl}
E.~Cheng, M.~Hyland, and J.~Power.
\newblock pseudodistributive laws.
\newblock {\em Electronic Notes in Theoretical Computer Science}, 83, 2003.

\bibitem{CurienP:opecdl}
P.-L. Curien.
\newblock Operads, clones, and distributive laws.
\newblock  Proc. of the International Conference on Operads and Universal Algebra, Nankai Series in Pure, Applied Mathematics and Theoretical Physics, Vol. 9, World Scientific, Singapore, 25-50 (2012).

\bibitem{DostalM:2-un-alg}
M. Dost\'al.
\newblock Two-dimensional universal algebra.
\newblock PhD Thesis, Czech Technical University in Prague, 2018.  
\newblock \url{https://math.feld.cvut.cz/dostamat/research/phd-thesis.pdf}.

\bibitem{EilenbergS:cloc}
S.~Eilenberg and G.~M. Kelly.
\newblock Closed categories.
\newblock In {\em Proceedings of La Jolla Conference on Categorical Algebra},
  pages 421--562. Springer-Verlag, 1966.

\bibitem{FioreM:carcbg}
M.~Fiore, N.~Gambino, M.~Hyland, and G.~Winskel.
\newblock The cartesian closed bicategory of generalised species of structures.
\newblock {\em Journal of the London Mathematical Society}, 77(2):203--220,
  2008.
  
\bibitem{FioreM:relpsm}
M.~Fiore, N.~Gambino, M.~Hyland, and G.~Winskel.
\newblock Relative pseudomonads, Kleisli bicategories, and substitution monoidal structures.
\newblock {\em Selecta Mathematica New Series}, 24: 2791, 2018. 


\bibitem{GambinoN:opebaf}
N.~Gambino and A.~Joyal.
\newblock On operads, bimodules and analytic functors.
\newblock {\em Memoirs of the American Mathematical Society}, 289(1184), (v) + 110pp, 2017.

\bibitem{GarnerR:polvpd}
R.~Garner.
\newblock Polycategories via pseudodistributive laws.
\newblock {\em Advances in Mathematics}, 218(3):781--827, 2008.

\bibitem{GarnerR:lowdsf}
R.~Garner and N.~Gurski.
\newblock The low-dimensional structures formed by tricategories.
\newblock {\em Mathematical Proceedings of the Cambridge Philosophical
  Society}, 146(3):551--589, 2009.

\bibitem{GordonR:coht}
R.~Gordon, A.~J. Power, and R.~Street.
\newblock Coherence for tricategories.
\newblock {\em Memoirs of the American Mathematical Society}, 117(558), 1995.

\bibitem{GurskiN:algtt}
N.~Gurski.
\newblock {\em Coherence in Three-Dimensional Category Theory}.
\newblock Cambridge University Press, 2013.

\bibitem{KellyG:cohtla}
G.~M. Kelly.
\newblock Coherence theorems for lax algebras and for distributive laws.
\newblock In G.~M. Kelly and R.~H. Street, editors, {\em Category Seminar
  (Proc. Sem., Sydney, 1972/1973)}, volume 420 of {\em Lecture Notes in
  Mathematics}, pages 281--375. Springer, 1974.

\bibitem{KellyG:bascec}
G.~M. Kelly.
\newblock {\em Basic {C}oncepts of {E}nriched {C}ategory {T}heory}.
\newblock Cambridge University Press, 1982.

\bibitem{KellyG:revetc}
G.~M. Kelly and R.H. Street.
\newblock Review of the elements of $2$-categories.
\newblock In G.~M. Kelly and R.~H. Street, editors, {\em Category Seminar
  (Proc. Sem. Sydney 1972/1973)}, volume 420 of {\em Lecture Notes in
  Mathematics}, pages 75--103. Springer, 1974.

\bibitem{LackS:cohapm}
S.~Lack.
\newblock A coherent approach to pseudomonads.
\newblock {\em Advances in Mathematics}, 152:179--202, 2000.

\bibitem{LackS:bicntg}
S.~Lack.
\newblock Bicat is not triequivalent to {G}ray.
\newblock {\em Theory and Applications of Categories}, 18(1):1--3, 2007.

\bibitem{LackS:fortm}
S.~Lack and R.~Street.
\newblock The formal theory of monads. {II}.
\newblock {\em Journal of Pure and Applied Algebra}, 175(1-3):243--265, 2002.


\bibitem{MacLaneS:catwm}
S.~Mac Lane.
\newblock {\em Categories for the working mathematician.}
\newblock{2nd ed.}, 
\newblock Springer, 1998.


\bibitem{MarmolejoF:docwsf}
F.~Marmolejo.
\newblock Doctrines whose structure forms a fully faithful adjoint string.
\newblock {\em Theory and applications of categories}, 3(2):24--44, 1997.

\bibitem{MarmolejoF:dislp}
F.~Marmolejo.
\newblock Distributive laws for pseudomonads.
\newblock {\em Theory and Applications of Categories}, 5(5):91--147, 1999.

\bibitem{MarmolejoF:dislpII}
F.~Marmolejo.
\newblock Distributive laws for pseudomonads {II}.
\newblock {\em Journal of Pure and Applied Algebra}, 195(1-2):169--182, 2004.


\bibitem{MarmolejoF:cohplr}
F.~Marmolejo and R.~J. Wood.
\newblock Coherence for pseudodistributive laws revisited.
\newblock {\em Theory and Applications of Categories}, 20(6):74--84, 2008.

\bibitem{AJPow90}
A.~J.~Power.
\newblock A $2$-Categorical Pasting Theorem.
\newblock {\em Journal of Algebra}, volume 129, pages 439-445, 1990.

\bibitem{StreetR:fortm}
R.~Street.
\newblock The formal theory of monads.
\newblock {\em Journal of Pure and Applied Algebra}, 2(2):149--168, 1972.

\bibitem{TanakaM:psedl}
M.~Tanaka.
\newblock {\em Pseudodistributive laws and a unified framework for variable
  binding}.
\newblock PhD thesis, Laboratory for the Foundations of Computer Science,
  School of Informatics, University of Edinburgh, 2004.
\newblock Available as LFCS Technical Report ECS-LFCS-04-438.

\bibitem{TanakaM:psedlav}
M.~Tanaka and J.~Power.
\newblock Pseudodistributive laws and axiomatics for variable binding.
\newblock {\em Higher-order and Symbolic Computation}, 19(2-3):305--337, 2006.

\bibitem{TanakaM:unicsb}
M.~Tanaka and J.~Power.
\newblock A unified category-theoretic semantics for binding signatures in
  substractural logics.
\newblock {\em Journal of Logic and Computation}, 16:5--25, 2006.

\bibitem{WalkerC:dislva}
C.~Walker. 
\newblock Distributive laws via admissibility.
\newblock {\em Applied Categorical Structures}, 27(6) 567--617, 2019.


\end{thebibliography}

\end{document}